\documentclass[twoside]{report}
\usepackage{epsf,latexsym,rawfonts}

\mathchardef\tnode="020E 

\def\arc{
  \hbox{\kern -0.15em
  \vbox{\hrule width 2.5em height 0.6ex depth -0.5 ex}
  \kern -0.33em}}

\def\darc{
  \rlap{\lower0.2ex\arc}{\raise0.2ex\arc}}

\def\farc{
  \rlap{\rlap{\lower0.45ex\arc}{\raise0.45ex\arc}}
       {\rlap{\lower0.15ex\arc}{\raise0.15ex\arc}}}

\def\tarc{
  \rlap{\rlap{\lower0.4ex\arc}{\raise0.4ex\arc}}
       {\arc}}

\def\stroke#1{
  \kern 0.05em
  \rlap\arc{{\textstyle{#1}}\atop\phantom\arc}
  \kern -0.22em}

\def\dstroke#1{
  \kern 0.05em
  \rlap\darc{{\textstyle{#1}}\atop\phantom\darc}
  \kern -0.22em}

\def\centerscript#1{
  \setbox0=\hbox{$\tnode$}
  \hbox to \wd0{\hss$\scriptstyle{#1}$\hss}}

\def\node{
  \def\super{}
  \def\sub{}
  \futurelet\next\dolabellednode}

  \let\sp=^
  \let\sb=_

  \def\dolabellednode{%
    \ifx\next\sb\let\next\getsub
    \else
      \ifx\next\sp\let\next\getsuper
      \else\let\next\donode
      \fi
    \fi
    \next}

  \def\getsub_#1{\def\sub{#1}\futurelet\next\dolabellednode}
  \def\getsuper^#1{\def\super{#1}\futurelet\next\dolabellednode}

  \def\donode{%
    \rlap{$\mathop{\phantom\tnode}\limits_{\centerscript{\sub}}
    ^{\centerscript{\super}}$}\tnode}

\def\varcdn{
  \kern -0.03em\vbox{\kern -0.5ex
  \hbox to \wd0{\hss\vrule width 0.04em depth 5.8ex\hss}
  \kern -0.3ex
  \hbox{$\tnode$}}}

\newtheorem{thm}{Theorem}[section]
\newtheorem{lem}[thm]{Lemma}
\newtheorem{res}[thm]{Result}
\newtheorem{cor}[thm]{Corollary}
\newtheorem{pr}[thm]{Proposition}

\newcommand{\qed}{\hfill\mbox{$\Box$}}

\newenvironment{pf}{\medskip {\bf Proof.\ }}{\qed\medskip}

\renewcommand{\O}{\Omega}
\newcommand{\Om}{\Omega}
\newcommand{\sm}{\setminus}
\newcommand{\Up}{\Upsilon}
\newcommand{\ep}{\varepsilon}
\newcommand{\G}{\Gamma}
\renewcommand{\L}{\Lambda}
\newcommand{\V}{\rm{V}}
\newcommand{\Aut}{{\rm Aut}}
\newcommand{\Th}{\Theta}
\newcommand{\Si}{\Sigma}
\newcommand{\Ga}{\Gamma}
\newcommand{\E}{{\rm E}}
\newcommand{\D}{\Delta}
\newcommand{\cL}{{\cal L}}
\newcommand{\cP}{{\cal P}}
\newcommand{\cH}{{\cal H}}
\newcommand{\cS}{{\cal S}}
\newcommand{\cM}{{\cal M}}
\newcommand{\cB}{{\cal B}}
\newcommand{\cE}{{\cal E}}
\newcommand{\cD}{{\cal D}}
\newcommand{\cX}{{\cal X}}
\newcommand{\cY}{{\cal Y}}
\newcommand{\cO}{{\cal O}}
\newcommand{\cG}{{\cal G}}
\newcommand{\cU}{{\cal U}}

\newcommand{\<}{\langle}
\renewcommand{\>}{\rangle}

\newcommand{\McL}{{\it McL}}
\newcommand{\Co}{{\it Co}}
\newcommand{\Suz}{{\it Suz}}
\newcommand{\BM}{{\it BM}}
\newcommand{\HS}{{\it HS}}
\newcommand{\HJ}{J_2}
\newcommand{\He}{{\it He}}
\newcommand{\Ru}{{\it Ru}}

\begin{document}
\pagenumbering{roman}
\begin{titlepage}
\bigskip
\begin{center}
{\huge Local Characterizations of\\ Geometries\\}
\vspace{2cm}
by\\
\vspace{2cm}
{\large Dmitrii V. Pasechnik\\}
\vspace{4cm}
This thesis is presented for the degree of\\
\bigskip
{\large Doctor of Philosophy\\}
\bigskip
of the\\
\bigskip
{\large University of Western Australia\\}
\vspace{3cm}
Department of Mathematics\\
{\large 1994\\}
\end{center}
\end{titlepage}

\begin{abstract}
One of the best  ways to understand the nature of finite simple groups
is   through  geometries  associated  with  them.  This  approach  for
classical and exceptional groups of Lie type has been quite successful
and has led to  the  deveopment of the concept of buildings and polar
spaces. The latter have been characterized by simple systems of axioms
with a combinatorial-geometric flavour.  It has been observed recently
that  geometries  similar to  buildings can be associated with  finite
sporadic simple  groups.  However, most of the known characterizations
of  such geometries for sporadic groups require additional assumptions
of a group-theoretic nature.  One  aim of  this thesis  is  to present
characterizations of  geometries for the sporadic groups $J_2$, $\Suz$,
$\McL$, $\Co_3$, $Fi_{22}$, $Fi_{23}$, $Fi_{24}$ and $\He$,  which are
in the same spirit  as the  characterizations  of buildings  and polar
spaces mentioned  above, in particular without  any assumption on  the
automorphism groups of the  geometries.  A by-product of these results
for $J_2$, $\Suz$ and $\He$  is a proof that certain presentations for
those groups are faithful.

Most of this work may be  viewed as a  contribution to  the  theory of
graphs  with  prescribed  neighbourhood.  The  result  on  graphs   of
$(+)$-points  of  GF(3)-orthogonal  spaces,  which  is also  used  for
characterization  of  geometries  for $Fi_{23}$  and $Fi_{24}$, may be
considered (along with the results for the groups $Fi_{22}$, $Fi_{23}$
and  $Fi_{24}$) as a generalization of a well-known theorem on locally
co\-triangular graphs.

Several of the geometries to be characterized are extensions of  polar
spaces.  Hyperovals  of  polar spaces  are natural  generalizations of
hyperovals of projective planes  of  even  order and play an important
role   in   investigations  of   extensions   of  polar   spaces.    A
computer-aided  enumeration   of  hyperovals  was   performed.   As  a
by-product  new  extended  generalized   quadrangles   were  found  as
hyperovals of the polar spaces $Q_5^+(4)$ and $H_5(4)$.
\end{abstract}

\tableofcontents
\chapter*{Preface}
\addcontentsline{toc}{chapter}{Preface}
\paragraph{Acknowledgments.}
\addcontentsline{toc}{section}{Acknowledgments}
The advice  and help of several people played an essential role in the
production of  this thesis. Especially, I thank  my  supervisor Cheryl
Praeger for her guidance, support, enthusiasm and tolerance during the
course  of  my study, and  for sharing with  me her vast experience of
scientific writing.

I started the research presented  in the thesis  while at Institute of
System Studies of Russian Academy of Sciences. I thank several present
or  former staff members of this  Institute,  in  particular Alexander
Ivanov, Sergey Shpectorov,  Sergey Tsaranov and my Diploma  supervisor
Igor  Faradjev, for their generosity with  their time and for numerous
extremely helpful discussions.

I thank Andries Brouwer, Hans Cuypers  and Antonio  Pasini for helpful
remarks,  discussions  and  suggestions.   I thank  several  anonymous
referees for their  careful reading of  the papers incorporated in the
thesis and  editors  of  the  corresponding  journals, namely  Francis
Buekenhout,  Arjeh  Cohen  and  Alexander  Ivanov,  for  their  prompt
responses.   I thank  Alice  Niemeyer for  editing  Chapter~9 and  for
sharing  her  valuable experience  in the use of modern computers, and
Roman Bogoyev for his readiness in  addressing  computer  software and
hardware related enquires.

Invitations to several  international  conferences (and the assistance
in covering travel expenses by various sources,  mainly the Department
of Mathematics of UWA) made it possible  for me to exchange ideas with
leading scholars in the field of research of this thesis.

I  acknowledge  the  support  of  an  Overseas  Postgraduate  Research
Scholarship of Australia and a  University Research Studentship of the
University of Western Australia during the course of my PhD.
\bigskip

\paragraph{Publications.}
\addcontentsline{toc}{section}{Publications}
The main part  of this  thesis consists  of  several research  papers,
which are  currently  at  various stages  of publication, as indicated
below.   The  references  made in  the  individual  papers  have  been
combined into a single  unified list  which is given at the end of the
thesis.

\begin{enumerate}
\item Chapter~\ref{chap:o3n} \cite{Pa:3tr}. 
On some locally 3-transposition graphs,
{\em in Proc. of Second Deinze Conf. ``Finite geometries and
combinatorics", LMS Lect. Notes Ser. 191, Cambridge Univ. Press, 1993,
319--326}.
\item Chapter~\ref{chap:EGQ33} \cite{Pa:EGQ33}.
The triangular extensions of a generalized quadrangle of order
$(3,3)$.  
{\em Bull. Belg. Math. Soc.}, 2:509--518, 1995.
\item Chapter~\ref{chap:suz} \cite{Pa:suz}.
Geometric characterization of graphs from the Suzuki chain.
{\em Eur. J. Comb. 14(1993), 491--499}.
\item Chapter~\ref{chap:fi2n} \cite{Pa:fi}.
Geometric characterization of the sporadic groups $Fi_{22}$,
$Fi_{23}$ and $Fi_{24}$. 
{\em J. Combin. Th. (A)}, 68:100--114, 1994.
\item Chapter~\ref{chap:epolsp} \cite{Pa:epolsp}.
Extending polar spaces of rank at least 3.  
{\em J. Comb. Th. (A)}, 72:232--242, 1995.
\item Chapter~\ref{chap:hex} \cite{CKP}.
Multiple extensions of generalized hexagons related to the simple
groups $\McL$ and $\Co_3$, (jointly with Hans~Cuypers and A.I.~Kasi\-ko\-va). 
{\em Proc. London Math. Soc.} (2) 54:16--24, 1996 
and {\em Research Report,
Department of Mathematics, University of Western Australia, July 1994,
No. 24.}
\item Chapter~\ref{chap:EGO} \cite{Pa:he}.
Extended generalized octagons and the group $\He$.
{\em Geom. Ded.}, 56:85--101, 1995.
\end{enumerate}

\chapter{Introduction}
\label{chap:intro}
\pagenumbering{arabic}

First,  in Section~\ref{intro:defno}  we  give  some  definitions  and
notation. This  is intended to  be  introductory  in nature, and to be
available for use as a  reference section when reading the rest of the
thesis. It may be skipped over providing the reader is familiar enough
with the terminology. Next in  Section~\ref{intro:tech} we outline the
approach  used  in  the   main  body   of  the  thesis  (that  is,  in
Chapters~\ref{chap:o3n}--\ref{chap:EGO}), aiming  at demonstrating the
common steps taken  there.  For  each  step we discuss the differences
and similarities among the chapters.

Section~\ref{intro:ex}  contains  a  model  problem, which illustrates
these steps in detail and also serves to demonstrate  the  flavour  of
results in the rest of the thesis.

Section~\ref{intro:sums} contains  summaries  of  the  results  in the
chapters comprising  the  main  body, including a  brief history and a
discussion of the ``completeness'' of the results in the sense that we
discuss whether one  might hope for stronger results. Interconnections
between  the chapters are pointed out  and related  results beyond the
scope of this thesis are surveyed.

\section{Definitions and notation}
\label{intro:defno}
Only rather general definitions are given here. More specific
definitions can be found in the chapters in which they are used.  The
main emphasis in this section is on definitions and facts regarding
diagram geometries. Although we do not use this general language in
the papers comprising the main part of this thesis (that is
Chapters~\ref{chap:o3n}--\ref{chap:EGO}), it helps to establish a link
between these papers. This general approach and language will be used
again in the final review Chapter \ref{chap:review}.

\subsection{Groups}
We assume knowledge of basic facts about permutation groups, classical
groups and finitely presented groups. Our notation on isomorphism types
of groups, extensions and so on will follow that of \cite{Atl}.
A standard reference is Aschbacher~\cite{As:fin}.
We give certain specific definitions here.

Let $G$ be a transitive permutation group acting on a set $\O$. 
The {\em rank} of $G$ is the number of orbits of the stabilizer
$G(\omega)$ of a point $\omega\in\O$ on $\O$. A rank 2 group is called
{\em doubly transitive}.

A $3$-{\em transposition group} is a group generated by a conjugacy
class $C$ of involutions (called 3-transpositions) such that for any
$a,b\in C$ either $(ab)^2=1$ or $(ab)^3=1$. 

We shall deal mainly with groups which arise as automorphism groups of
various incidence structures: graphs, diagram geometries, etc.

\subsection{Graphs}
Our notation is fairly standard, and mostly may be found in, e.g.,
\cite{BCN}.

We consider undirected graphs without loops
and multiple edges.
Given a graph $\G$, let us denote  the  set  of  vertices  by
$\V=\V\G$,  the  set  of edges by $\E=\E\G$. Let $X\subseteq
\V\G$.
We denote by $\langle X\rangle =\langle X\rangle _{\G}$ the subgraph $\Xi$
of $\G$ {\it induced} by $X\/$
(i.e. $\V\Xi=X,\/\E\Xi=\{ (u,v)\in \E\G| u,v\in X\}$).
Given two graphs $\G$ and $\D$, the graph $\G\cup\D$
(resp. the graph $\G\cap\D$) is
the graph with vertex set $\V\G\cup \V\D$ (resp.
$\V\G\cap \V\D$) and edge set $\E\G\cup \E\D$
(resp. $\E\G\cap \E\D$). Given $v\in \V\G$, we denote
$\G_i(v)=\langle \{x\in \V\G\mid x $ at distance $i$ from $v\}\rangle $,
and often write $\G(v)=\G_1(v)$.
Furthermore, $\G(X)=\bigcap_{x\in X} \G(x)$.
To simplify the notation we use
$\G(v_1,\dots,v_k)$ instead of $\G(\{v_1,\dots,v_k\})$ and
$u\in\G_i(\dots)$ instead of $u\in \V\G_i(\dots)$.

As usual, $k=k(\G)=|\G(x)|$, and $d(\G)$ is the
maximal distance between a pair of vertices occurring in $\G$.
For $w\in\G_i(u)$ we write $c_i=|\G_{i-1}(u)\cap\G(w)|$,
$b_i=|\G_{i+1}(u)\cap\G(w)|$. If $k=b_0$ is independent of the
choice of $x\in\V\G$, then $\G$ is said to be {\em regular}.
If for a connected graph $\G$
the numbers $c_i$ and $b_i$ are independent of the 
choice of vertices $u$, $w$ at distance $i$ for
$i\in\{0,\dots,d(\G)\}$ then $\G$ is called {\em distance regular}.
The parameters $c_i$ and $b_i$ are usually presented, for a given
graph $\G$ with $d=d(\G)$, as an {\em intersection array}:
$$\{b_0,\dots,b_{d-1};c_1,\dots,c_d\}.$$
We denote $\lambda=k-b_1-1$ and $\mu=c_2$.
A distance-regular graph of diameter 2 is called {\em strongly regular}.
 
If $\D$ is a (proper) subgraph of
$\G$ we denote this fact by $\D\subseteq\G$ (resp.
$\D\subset\G$). 

We denote  the complete $n$-vertex graph by $K_n$,  the complete
bipartite graph with parts of size $n$ and $m$
by $K_{n,m}$, the circuit of length $n$ by
$C_n$, and the empty graph by $\emptyset$.

Let $\G$, $\D$ be two graphs. We say that $\G$ is
{\it locally} $\D$ if $\G(v)\cong\D$ for each $v\in \V\G$.
More generally, if $\cD$ is a class of graphs then $\G$ is
locally $\cD$ if for each $v\in\V\G$ the subgraph $\G(v)$ is isomorphic
to a member of $\cD$.

Let $\G,\ \overline\G$ be two graphs. We say that $\G$ is a
{\it cover} of $\overline\G$ if there exists a mapping $\varphi$
from $\V\G$ to $\V\overline\G$ which maps edges to edges.
We shall be particularly interested in covers satisfying the
additional assumption that for
any $v \in \V\G$ the restriction of $\varphi$ to $\G(v)$ is an
isomorphism onto $\overline\G(\varphi(v))$. 

The automorphism group of a graph $\G$ will be denoted by $\Aut(\G)$.
A distance-regular graph $\G$ is called {\em distance-transitive} if the
rank of $\Aut(\G)$ equals $d(\G)+1$. If $d(\G)=2$ then $\G$ is called a
{\em rank} 3 graph.

For a 3-transposition group, the graph $\G$
with the set of 3-transpositions
as the vertex set, such that
two vertices are adjacent if the corresponding  involutions
commute, is called a {\em $3$-transposition graph}.
In most important classes of 3-transposition groups, $\G$ is a rank
3 graph.

\subsection{Incidence systems}  
Here we follow \cite{CHP}.

An {\em incidence system} $\G$ is a pair $\G(\cP,\cB)$ of sets $\cP$
and $\cB$, where the elements of $\cB$ are subsets of $\cP$ of size at
least 2. The elements of $\cP$ and $\cB$ are usually called {\em
points} and {\em blocks}, respectively, and we write $\cP(\G)=\cP$ and
$\cB(\G)=\cB$.  Incidence between points and blocks is given by
inclusion.  The {\em point graph} of $\G$ is the graph with vertex set
$\cP$ such that two vertices $p$ and $q$ are {\em adjacent} (notation
$p\perp q$) if there is a block containing both of them.  For any
$X\subseteq\cP$ we denote $X^\perp=\{p\in\cP\mid p\perp x {\rm\ for\
all\ } x\in X\}$.  The {\em subsystem} of $\G$ induced by the set $X$
is the incidence system $\G(X\cap\cP,\{B\cap X\mid B\in\cB, |B\cap
X|>1\})$.  We call $\G$ {\em connected} if its point graph is
connected.  The {\em incidence graph} of $\G$ is the bipartite graph
with vertex set $\cP\cup\cB$, such that vertices are adjacent if they
are incident.

Let $\G$ be an incidence system.
Given $X\subset\cP(\G)$ satisfying $|B-X|>1$ for all $B\in\cB(\G)$, 
we refer to the incidence system 
$$\G_X=(\G(X), \{B-X\mid X\subset B\in\cB(\G)\})$$
as the {\em residue} of $X$, or as the 
{\em local} system in $X$. Let $\cD$ be a class of incidence systems. 
$\G$ is called an {\em extension} (respectively {\em $k$-fold extension})
of $\cD$, or ({\em $k$-fold}) {\em extended} $\cD$, if, for any
$x\in\cP(\G)$ (respectively, for any $X\subset\cP(\G)$, $|X|=k$,
such that there
exists $B\in\cB(\G)$, $X\subset B$), the residue $\G_x$ (respectively $\G_X$)
is isomorphic to a member of $\cD$. If $\cD=\{\D\}$,
it is customary to refer to extensions of $\cD$ as extensions of $\D$
(or extended $\D$). 
If $\G$ is an extension of $\cD$ then 
the connected components of $\G$ are extensions of $\cD$,
also. Hence, unless otherwise stated, we assume our extensions to be connected.

\subsection{Some examples of incidence systems}
Here some well-known examples of incidence systems are given.

An  incidence system  is  called  a {\em  generalized  $d$-gon} if its
incidence graph has diameter  $d>2$ and  girth $2d$.  It  follows that
there is at most  one block  on any pair of  points, so the blocks  of
$\G$ are usually called  {\em lines}, and the term {\em  collinearity}
is  used instead of {\em  adjacency}.  Given $\G$, we define the  dual
system $\G^*$, in fact also  a  $2d$-gon, to be  the incidence  system
whose points (respectively lines) are  the lines (respectively points)
of $\G$,  and incidence is by inverse inclusion.   We say that $\G$ is
{\em nondegenerate} if for each point there exists a point at distance
$d$  from it, and the same holds in $\G^*$.  A line of $\G$ is said to
be {\em thick} if it contains more than 2 points, otherwise it is said
to  be {\em thin}. We say that $\G$ is {\em regular} of order  $(s,t)$
if each point is on exactly $t+1$ lines and each line contains exactly
$s+1$ points.   The numbers $s$ and $t$ are called {\em parameters} in
this case.

If $d=3$ then $\G$ is a projective plane.
Is customary to refer to a generalized 4-gon as a {\em generalized quadrangle}
(GQ, for short, or GQ$(s,t)$). Similarly, we refer
to generalized 8-gons as {\em generalized octagons} (respectively GO and
GO$(s,t)$).

A {\em polar space} $\G$ is the incidence system of isotropic points 
(that is, 1-dimensional subspaces) and
totally isotropic lines (that is, 2-dimensional subspaces)
of a vector space $V$, 
equipped with a nondegenerate symmetric,
alternating, or Hermitian form. Given a point $P$ and a line $\ell$ not
on $P$, $P$ is collinear to one or all points of $\ell$ (the
so-called {\em one-or-all axiom}). If both alternatives happen, then $\G$
is called {\em nondegenerate}. Note that if the second alternative never
happens then $\G$ is a GQ.
The {\em rank} of $\G$ is the (projective) dimension of a maximal
totally isotropic subspace of $\G$.

In fact, the one-or-all axiom, together with certain nondegerancy
conditions, can be taken as a definition of (nondegenerate) polar
spaces, see e.g. \cite{BuSh,BC}.

An {\em affine polar space} is the incidence system of points, lines,
etc. of a polar space with a geometric hyperplane removed. This is a
natural generalization of the construction of an affine space from a
projective space.
A {\em standard quotient} of an affine polar space $A$ is the incidence
system of (possibly nonmaximal) equivalence classes of points, lines,
etc. at mutual distance 3 in the point graph of $A$.
See \cite{CoSh,CuPa} for further discussion and an axiomatic
characterization of those objects.

\subsection{Diagram geometries}
Here we follow mainly \cite{Pas:book,Bue:geo:spo,Cam:PGPS}, see also
\cite{As:spo}. 

Our definition of geometries is inductive, as follows.
\begin{description}
\item[$(i)$] A {\em geometry of rank $1$} is just a set with at least
two elements (namely, the null graph with at least 2 vertices).
\item[$(ii)$] A {\em geometry of rank $n\ge 2$} is a pair $(\G,I)$,
where $\G$ is a connected $n$-partite graph and $I$
is its $n$-partition such that for any vertex $x$ of $\G$ the subgraph
$\G(x)$ is a geometry of rank $n-1$ with respect to the partition
of $V(\G(x))$ induced by $\G$ and $I$.
\end{description}
The classes of $I$ are usually called {\em types}.
There is in fact no problem in identifying $(\G,I)$ with the graph $\G$,
since the partition can be recovered from $\G$, and we shall denote
$(\G,I)$ just by $\G$. 
The vertices of $\G$ are called {\em elements}.
It is convenient to introduce the 
{\em type function $t$} from the set of elements to the set of types.
Elements are {\em incident} if they are adjacent in $\G$.
A set $X$ of pairwise incident elements of $\G$ is called a {\em flag}.
Given a flag $X$, $|X|<n$, the geometry $\G(X)$ is called the 
{\em residue} of $X$. Also, the restriction of $\Theta$ onto $\G(X)$ 
(i.e. $I-t(X)$) is called the {\em cotype} of $X$. Then,
$|I-t(X)|$ is called {\em corank} of $X$.

Let $X$ be a flag of cotype $\{\alpha\}$. $\G$ is said to admit the {\em order}
$a$ in $\alpha$ if $a=|\G(X)|-1$ does not depend upon the
particular choice of $X$ of cotype $\{\alpha\}$.

Note that a rank 2 geometry is simply a connected bipartite
graph such that every vertex is adjacent to more than one vertex.
An incidence system, as defined above, is a geometry of rank
2 (with types called points and blocks)
satisfying the additional property that if two vertices $x$, $y$ of
$\G$ correspond to {\em blocks} of the incidence system then 
$\G(x)\not=\G(y)$.
 
A {\em diagram} $\cD$ of $\G$ is a (possibly directed)
graph with vertex set $\Theta$, such that 
each edge/nonedge $(u,v)$ is labeled by the class of
isomorphism types of the residues of flags of cotype $\{u,v\}$.
Also, the vertices of $\cD$ are sometimes labeled by the names of types
and/or their orders. The idea of associating a diagram with $\G$ arises
from the observation that often properties of rank 2 residues of $\G$
greatly influence properties of $\G$ itself.

\subsubsection{Some rank 2 diagrams}
To be able to draw and understand diagrams, a common convention has been
developed for most important rank 2 geometries/incidence systems.

A {\em generalized $n$-gon} is denoted by 
$\node\stroke{(n)}\node$.
For $n=2,3,4,6,8$ the symbol $(n)$ is usually omitted and
$0,1,2,3,4$ ordinary parallel bonds respectively are drawn instead.
For instance, to denote $GQ(s,t)$ we draw 
$\node_s\darc\node_t$.
An explanation is needed for the case $n=2$. Indeed, a generalized digon
is not an incidence system as defined above. We define a
{\em generalized digon} to be
a rank 2 geometry whose incidence graph is complete bipartite.
\medskip

A {\em circular space}, which is the incidence system of elements and
(unordered) pairs of elements of a set, is denoted by
$\node\stroke{\subset}\node$.

\subsubsection{Some rank $n$ diagrams}
\label{sssubs:rkn}
To illustrate the convention introduced above, we give examples of
diagram for geometries of rank bigger than 2.
If a diagram looks like a path, we say that it is {\em linear}.
Throughout this thesis, mainly geometries with linear diagrams are
considered.

A projective space of (projective) dimension $n$ over GF($q$) has the
following diagram.
$$\node_q\arc\node_q\cdots\node_q\arc\node_q\arc\node_q,$$
the diagram has $n$ nodes, node number $i$ from the left depicts the
$i$-varieties.


A $k$-fold extension $\G$ of a class $\cX$ of incidence
systems has the following diagram.
\begin{equation}
\label{intro:d1}
\node\arc\node\cdots\node\arc\node\stroke{\subset}\node\stroke{\cX}\node
\end{equation}
-with $k+2$ nodes,
where the left type corresponds to points of $\G$, 
the right type corresponds to
blocks of $\G$ and the middle type(s) -- to pairs, triples,\dots, $k+1$-tuples
of points lying in blocks(s).
The following is a specialization of this example. 
Let $\cX=\Si_1$,\dots, $\Si_k$ be a chain of graphs such that
$\Si_{i+1}$ is locally $\Si_i$ for each $i=1,\dots,k-1$.
Then $\D=\Si_k$, $k>1$, gives rise to a geometry with the following diagram.
$$\node\arc\node\cdots\node\arc\node\stroke{\cX}\node$$
with $k+1$ nodes, where, counting from the left-hand side, elements of type
$i$ are complete $i$-vertex subgraphs of $\D$, for $1\le i\le k+1$.

For $\G$ with diagram (\ref{intro:d1}) we use the abbreviation
$c^k.\cX$, or $c.\cX$ for $k=1$.

The diagram
$$\node^1\arc\node^2\cdots\node^{k-1}\arc\node^k\stroke{\subset}
\node^{k+1}_s\arc\node^{k+2}_s\cdots\cdots\node^{k+n-2}_s\arc\node^{k+n-1}_s
\darc\node^{k+n}_t$$
for a $k$-fold extension of a polar space of rank $n$ 
(over GF($s$), providing $n>2$, or with lines of size $s+1$ otherwise)
will be abbreviated as $c^k.C_n(s,t)$, or $c.C_2(s,t)$ for $k=1$.
Here the labels over the nodes are just the numbers of the types.

\subsubsection{Morphisms, coset geometries}
Let $\G$, $\G'$ be two geometries over the same set of types $I$,
$|I|=n$, and let $t$, $t'$ be the type functions of $\G$, $\G'$, respectively.
A (type-preserving) {\em morphism} $f:\G'\rightarrow\G$ is a covering of
graphs such that $t'(x)=t(f(x))$.
A morphism is called a {\em $k$-cover} if $\G'(X)\cong\G(f(X))$ for
any flag $X$ of $\G'$ of corank $k$.
It is natural to call an $n$-covering an {\em isomorphism}, and
usually $(n-1)$-covers are called covers (without prefix). 

We say that an $m$-cover $f:\G'\rightarrow\G$  is {\em universal} if, for any
$m$-cover $g:\G''\rightarrow\G$, there is a unique morphism 
$h_g:\G'\rightarrow\G''$ such that $gh_g=f$.

One of the important results in this theory says the following.
\begin{thm}
Given a geometry $\G$, its
universal cover exists and is unique (up to isomorphism).\qed
\label{Univcover}
\end{thm}

An {\em automorphism group} $G$ of $\G$ is a group of isomorphisms of $\G$
onto itself, such a group $G$ will be a subgroup of the group $Aut\G$
of all the isomorphisms of $\G$ onto itself. 
It is {\em flag-transitive} if for any flags $X,X'$ of the
same cotype there exists $g\in G$ such that $X^g=X'$.

Let $C$ be a maximal flag of $\G$. For any $i\in I$, let $x_i\in X$ be
the element of $C$ of type $i$ and let $G_i$ be the stabilizer of $x_i$ in $G$.
For any $J\subseteq I$, let $G_J=\bigcap_{j\in J} G_j$.

Define a geometry $\G(G)$ as follows. The set of elements is the set
of all cosets of all of the $G_i$ in $G$ and two elements are incident if the
corresponding cosets have nonempty intersection.
The following two results form the basis for investigation of
flag-transitive geometries from the group-theoretic point of view.

\begin{pr}
Let $G$ be flag-transitive on $\G$. Then $\G(G)\cong\G$.\qed
\label{Cosetgeom}
\end{pr}

The group $G$ is an {\em amalgamated product} of $G_i$, $i\in I$,
with amalgamation given by $G_J$, where $J$ runs through the set 
$I^2$ of all the
unordered pairs of $I$. In other words, $G$ is a quotient of the group
$\hat{G}$ given by a
presentation ${\cal A}={\cal A}(G)$, 
where the generators are elements of $G_i$,
the relations are those holding in $G_i$,
$i\in I$, and for $J=\{i,j\}$ and $g_i\in G_i$, $g_i=g_j$, the
corresponding element in $G_j$, if and only if $g_i\in G_J$. Moreover,
the kernel $K$ of the homomorphism $\phi :\hat{G}\rightarrow G$ satisfies
$K\cap G_i=\{1\}$ for any $i\in I$. 

Note that the group $\hat{G}$ is usually called the {\em universal closure} of
${\cal A}$.
\begin{thm}
Let $G$ be flag-transitive on $\G$. 
Then the homomorphism $\phi :\hat{G}\rightarrow G$
induces the universal cover $f :\G(\hat{G})\rightarrow\G$.
\label{Univclosure}
\end{thm}
Note also that Theorem~\ref{Univclosure} is useful when we need a proof that a
presentation for $G$ is faithful. Namely, if we show that $\G=\G(G)$ does not
possess nontrivial covers, then it follows that the presentation ${\cal A}$
for $G$ is faithful.

\subsubsection{Point graphs}
Let $\G$ be a geometry with a linear diagram. Call the left-hand side type of
the elements {\em points}. The {\em point graph} $\Si$ of $\G$ is the graph
whose vertices are points, two vertices are adjacent if the corresponding
points are incident to a common element of type second (from the left) side.

The geometries we shall consider satisfy the following property.
\begin{itemize}
\item[(LL)] Each pair of points is incident to at most one element of the second
type.
\end{itemize}
If $\G$ is the geometry for an incidence system, as defined in
\ref{sssubs:rkn}, then obviously $\Si$ is just the point graph for this system.
In most cases we shall consider $\G$ can be recovered from $\Si$,
and we write $\G=\G(\Si)$ without any ambiguity.
Namely, the following holds.
\begin{itemize}
\item[(T)] Each point in any triple of pairwise adjacent points is incident
to an element of the third type.
\end{itemize}
In fact, (T) is equivalent to the property that, for any point $p$, the
subgraph $\Si(p)$ is isomorphic to the point graph of the residue $\G(p)$ of
$p$.
If (T) holds,
covers of $\G$ give rise to covers of $\Si$ inducing isomorphisms
on the neighbours of each vertex of $\Si$, and vise versa (that is,
the covers of $\Si$ inducing isomorphisms on $\Si(x)$ for any
$x\in\V\Si$ are covers of $\G$).

\section{Commentary on the approach} \label{intro:tech}
Here we briefly discuss our approach and techniques.
Suppose we wish to classify the geometries $\G$ with a
given diagram. As usual it is assumed that $\G$ satisfies (LL).
As well, we assume that $\G$ satisfies a condition which could be 
used to translate the problem into a classification problem regarding
graphs with prescribed neighbourhood. The most usual condition
is (T). However, note that
in Chapters~\ref{chap:hex}
and \ref{chap:EGO} a different type
of such condition is used, whereas in 
Chapter~\ref{chap:epolsp} 
we manage to {\em show} that (T) holds. 

So, here we shall slightly abuse our notation and identify 
$\G$ with its point graph.
Thus, our problem is translated into the problem of classifying
the graphs $\G$ which are locally $\cD$, where $\cD$ is a fixed 
class of graphs.
Mostly the problem is formulated (or can be easily shown to be formulated)
in such a way that $\cD$ consists of a
unique graph $\D$ (or at least that $\G(v)\cong\D$ for any $v\in\G$ and
any {\em particular} $\G$), 
with the exceptions of Chapters~\ref{chap:EGQ33} and
\ref{chap:EGO}. In particular in Chapter~\ref{chap:EGO} 
it takes a fair amount of 
effort to establish that $\cD$ is not too ``wild''.
So, here we assume that we already have that $\cD=\{\D\}$,
that is $\G$ is locally $\D$, and $\D$ is a {\em known} (connected) graph.
The following simple general fact is well-known. For  a particular
case (and a proof), see   Lemma~\ref{suz:L1}.

\begin{lem}
Let $u\in\V\G$, $v\in\G_2(u)$. The subgraph $\G(u,v)$ is locally $\cM$,
where $\cM$ is the family of the isomorphism types of the subgraphs
$\D(x,y)$, where $x$ and $y$ run through the set of pairs $(x,y)$
of nonadjacent vertices of $\D$.
\qed \label{intro:mugraphs}
\end{lem}
In all the cases we consider $\D$ has diameter 2, so $\cM$ does
not contain the empty graph. 

\paragraph{Digression.} Generally speaking, if the diameter of $\D$ 
is greater than 2, the problem of finding all the graphs
$\G$ might be very difficult.
It is worth mentioning two related results \cite{Wee1,Wee2}.
In \cite{Wee1} it is shown that if $\D$ is of girth at least 6 then there
always exists an infinite locally $\D$ graph $\G$.
In \cite{Wee2} it is shown that for ``sufficiently dense'' graphs $\D$ of
diameter 2 any locally $\D$ graph is finite. 
\medskip

Going back to our consideration, one of the main topics in this thesis
is the problem of determining all the possible $\G(u,v)$ using $\cM$. In
Chapters~\ref{chap:o3n}, 
\ref{chap:suz}, \ref{chap:fi2n} and 
\ref{chap:hex},  we manage to enumerate
all the possibilities for $\G(u,v)$ using ``local'' information from $\cM$ 
(here $\cM$ consists of the unique member). The key observation is that the
problem of defining all the locally $\cM$ {\em subgraphs}
of $\D$ is easier than the problem of defining all the locally $\cM$ 
{\em graphs}.

In Chapters~\ref{chap:EGQ33} and \ref{chap:epolsp} we were forced to use
a computer to enumerate the possibilities for $\G(u,v)$.
The technique used in Chapter~\ref{chap:EGO} is rather specific and
so we do not discuss it further here.

The next stage is to reconstruct the possibilities for the second layer
of $\G$. Here the approach is to pick a particular $\G(u,v)$ obtained at
the previous step and try to reconstruct $W=\G(v)\cap\G_2(u)$.
It turns out that the choice of appropriate $\G(u,v)$ for $w\in W$
is usually very restricted. Usually we are able to reconstruct
the subgraph induced on $\G_2(u)$ in terms of subgraphs
$\G(u,x)$, where $x\in\G_2(u)$. One important point is that in the cases
considered in Chapters~\ref{chap:o3n}, \ref{chap:EGQ33}, \ref{chap:fi2n} and
\ref{chap:epolsp} the subgraph $\G_2(u)$ appears as a proper cover of
a certain graph defined on the sets $\G(u,x)$, $x\in\G_2(u)$,
namely, in these cases $\G(u,x)=\G(u,x')$ for $x\not=x'\in\G_2(u)$.
The trick used there was to study the partial linear space $\cL$ on the
vertices of $\G$, where lines are defined as  the cocliques
$X\subset\G$ such
that $\G(x,x')=\G_X\not=\emptyset$ is independent of the particular
choice of a pair $x,x'\in X$. In these cases $\cL$ turned out to be
(generalized) Fischer spaces, objects which have been completely classified.

In most cases the graph $\G$ turned out to be of diameter 2, so this
step was the final one. Note that, if necessary, this step was performed
several times, until all the possibilities for $\G(x,y)$ were exhausted.

In the remaining cases $\G$ was either shown to be a cover of a locally
$\D$ graph of diameter 2 (Chapters~\ref{chap:suz}, \ref{chap:fi2n}), 
or reconstructed directly (Chapters~\ref{chap:hex}, \ref{chap:EGO}).

\section{An example} \label{intro:ex}
Here we consider an example illustrating the discussion in the previous
section. The following statement has been known at least 
since the work \cite{BH} by Buekenhout and Hubaut. 
\begin{pr}
There are exactly $2$ nonisomorphic $c.C_2(2,1)$-geometries $\G$
satisfying {\rm (LL)} and {\em (T)}, with $16$ and $20$ points,
respectively.
\end{pr}
\begin{figure}[bht]
\begin{center}
\leavevmode
\epsffile{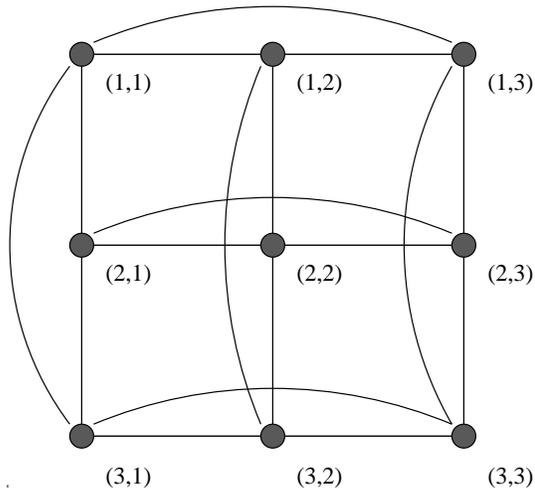}
\end{center}
\caption{The $3\times 3$-grid graph $\D$.}
\label{intro:fig1}
\end{figure}
We sketch the proof here. First, we observe that the problem is
equivalent to the classification of the locally $\D$ graphs $\G$, where
$\D$ is the $3\times 3$ grid graph, that is the graph whose vertices are
the pairs $(i,j)$, $i,j\in\{1,2,3\}$, such that 
two vertices $(i,j)$ and $(m,n)$
are adjacent if either $i=m$ or $j=n$, see Figure~\ref{intro:fig1}.
For such $\D$, the family $\cM$ (cf. Lemma~\ref{intro:mugraphs}) 
consists of the unique element, namely
the null graph on 2 vertices. Moreover, the subgraphs of $\D$
isomorphic to the 2-vertex null graph lie in one orbit of $G=Aut(\D)$.

Next, the set of all possible $\O=\G(u,v)$,
where $u\in\G$ and $v\in\G_2(u)$, consists of the
quadrangles and hexagons. Indeed, without loss of generality,
$x=(1,1)\in\O$ and $\O(x)=\{y,y'\}$, where $y=(1,2)$, $y'=(2,1)$.
Since $\O$ is locally $\cM$, the vertices $(1,3)$ and $(3,1)$ of
$\D=\G(u)$ do not belong to $\O$. Thus $\O(y)$ is either
$\{(1,1),(2,2)\}$ or $\{(1,1),(3,2)\}$. 

In the former case the
connected component of $\O$ containing $x$ is the quadrangle formed by
$x,y,y'$ and $(2,2)$. The only vertex of $\D$ nonadjacent to a vertex
in the quadrangle is $(3,3)$, so $\O$ is in fact connected.
Similarly in the latter case we find that $\O$ is the hexagon formed
by all the vertices of $\D$ excluding $(1,3)$, $(2,2)$ and $(3,1)$.

There are 9 quadrangles $Q$ in $\D$, for any $Q$ there exists $x\in\D$ such
that $Q=\D_2(x)$. Also, there are 6 hexagons $H$ in $\D$, for any $H$
there exists a 3-coclique $T$ such that $H=\D-T$.
Note that $Aut(\D)$ acts transitively on the quadrangles and on the hexagons.

Let $v\in\G_2(u)$. Then $\G(u,v)=\G(u,v')$ implies $v=v'$.
Indeed, let $x\in\G(u,v)$. Then, by inspection of $\G(x)$,
the common neighbourhood of $\G(u,v)\cap\G(x)$ in $\G(x)$ consists of
$u$ and $v$, so $v=v'$.

Assume first that $\G(u,v)$ is a hexagon.
Let $w\in\G(v)-\G(u)$. \\
Then $|\G(u,v,w)|=4$. Since the subgraph induced
on $\G(u,v,w)$ is not a quadrangle, $\G(u,w)$ is a hexagon, as well.
Let ${\cal H}$ denote the set of all the hexagons in $\D$. Let $\Si$ be
the graph with vertex set ${\cal H}$, such that two vertices $H$, $H'$ are
adjacent if $|H\cap H'|=4$. It is straightforward to see that $\Si\cong
K_{3,3}$. Now, since $|\G(v)-\G(u)|=3$, $\G_2(u)\cong\Si$. 
The remaining check that the graph just constructed is locally $\D$ is
straightforward. Finally, note that $\G(x,y)$ is a hexagon for any pair
of nonadjacent vertices of $\G$.

It remains to consider the case when $\G(u,v)$ is a quadrangle. 
As observed, in this case $\G(x,y)$ is a quadrangle for any pair of
vertices at distance 2. Counting in two ways the edges between $\G(u)$
and $\G_2(u)$, we see that $|\G_2(u)|=9$. So all the quadrangles in
$\G(u)$ appear as $\G(u,x)$ for some $x\in\G_2(u)$.
Let $w\in\G(v)-\G(u)$. Then $|\G(u,v,w)|=0$ or 2. 
Let $|\G(u,v,w)|=0$. The case $w\in\G_2(u)$ is impossible, since
intersection of any quadrangle of $\G(u)$ with $\G(u,v)$ is nontrivial.
So $w\in\G_3(u)$. 

Thus $|\G(u,v,w)|=2$ iff $w\in\G_2(u)$.  The graph $\Si$
defined on the set of all quadrangles of $\G(u)$, such that two vertices
are adjacent if the intersection of the corresponding quadrangles is of size
2, is isomorphic to $\D$.  Since $|\G(v)\cap\G_2(u)|=4$,
$\G_2(u)\cong\D$. 

It remains to show that $\G(w)=\G_2(u)$. As already noted, each
vertex $x\in\G_2(u)$ is adjacent to a unique vertex $w_x\in\G_3(u)$.
By inspection of $\G(v)$,  
we have that $w_x=w_v$ for any $x\in\G(v)\cap\G_2(u)$.
Since $\G_2(u)$ is connected, our claim holds.
It is straightforward to check that the graph just constructed is indeed
locally $\D$.

\section{The results} \label{intro:sums}
Here we give summaries for each of the chapters comprising the main body
of the thesis. We give statements of results, a brief historical
background of the problems, and a discussion on whether stronger
results might be possible.
Moreover, we point out several connections between these chapters and
survey related results of the author which have not been included
in this thesis.
A more general discussion of the place of these results
in a broader context will be given in Chapter~\ref{chap:review}.
We remind the reader that the information on where the contents of the
Chapters~\ref{chap:o3n}-\ref{chap:EGO} are published is given in the
Preface.

\paragraph{Chapter~\ref{chap:o3n}} is a contribution to the theory of
$c^k.C_2(s,t)$-geometries for $s=2$ or $4$ and $t=2$.
We assume that the geometries satisfy condition (T). Essentially, 
we improve here results from \cite{Mei}, where similar results were
proved under the stronger assumption of flag-transitivity.

It is well-known that there are exactly 3 nonisomorphic
$c.C_2(2,2)$-geomet\-ries, with 28, 32, and 36 points, respectively, cf.
\cite{BH,CHP}. The situation with $c^k.C_2(2,2)$-geometries is less
clear. According to \cite{HaSh}, the $c^k.C_2$-geometry having 
$c.C_2$-residues isomorphic to the 28-point example is unique. Its
point graph is the 3-transposition graph for the symmetric group on
$2k+6$ symbols. 

The 32-point example is not extendable, cf. \cite{CHP}.

Denote by $\cX$ the 36-point $c.C_2(2,2)$-geometry. By \cite{HaSh}, a
$c^2.C_2$-geometry having at least one point residue isomorphic to $\cX$ is
unique. The point graphs of this geometry and of $\cX$ are members of
the family $\cO^-_k$ of 3-transposition graphs for the orthogonal
group $O^\epsilon_{k+4}(3)$ (for $k=1,2$), where $\epsilon$ is empty if
$k$ is odd, and otherwise is $`-'$ or $`+'$ depending on whether $4|k$ or
not. We characterize the geometries $\G(\cO_k^-)$ for $k\ge 4$.
We show that if a $c^k.C_2$-geometry $\G$ has
$c^3.C_2$-residues isomorphic to $\G(\cO_3^-)$ then $\G\cong\G(\cO_k^-)$.

According to \cite{Mei}, there exists a $c^4.C_2$-geometry $\cY$ whose
$c^3.C_2$-residues are  proper triple covers of $\G(\cO_3^-)$, and $\cY$
is a $2$-cover of $\G(\cO_4^-)$. It follows that there exist many
nonisomorphic $c^4.C_2$-geometries, all of them quotients of $\cY$, some
of whose $c^3.C_2$-residues are isomorphic to $\G(\cO_3^-)$ and some to
its proper triple cover.
Hence our characterization appears to be close to best possible, in
some sense.

The author has shown that $\G(\cO_3^-)$ and its proper triple cover are
the only $c^3.C_2$-geometries having a $c.C_2$-residue isomorphic to
$\cX$ [in preparation].
\medskip

At present there are 3 nonisomorphic $c.C_2(4,2)$-geometries known. They
have 126, 378, and 162 points, respectively. The first one, to be
referred to as $\cX$, has its point graph $\cO_1^+$ in the family $\cO^+_k$
of 3-transposition graphs for the orthogonal group
$O_{k+5}^\epsilon(3)$, where $\epsilon$ is empty if $k$ is even,
or is $`-'$ or $`+'$ depending on whether $4|(k+5)$ or not.
The second example is the proper triple cover of the first one, see e.g.
\cite{CHP,Mei}.
The third example is new and is described in Chapter~\ref{chap:epolsp}.
At present the problem of classifying $c.C_2(4,2)$-geometries seems to
be hopeless.
We characterize the geometries $\G(O_k^+)$.
Namely, we show that if a $c^k.C_2$-geometry has the $c.C_2$-residues
isomorphic to $\cX$ then $\G\cong\G(O_k^+)$. 
This result will be used in Chapter~\ref{chap:fi2n} to classify
$c^m.C_3(4,2)$-geometries for $m>1$.

Again, for the same reason as for $\G(O_k^-)$, that is, the existence of a
proper 2-cover of $\G(\cO_2^+)$, our result is close to the best
possible one.

\paragraph{Chapter~\ref{chap:EGQ33}} is concerned with a proof that the
$c.C_2(3,3)$-geometry satisfying (T) is unique. 
This was shown under the additional assumption of flag-transitivity in
\cite{BH}, see also \cite{Mei,DFGMP}.
Along with \cite{BlBr44}, where $c.C_2(3,1)$-geometries satisfying (T)
were classified, and the author's paper \cite{Pa:mclco3}, where
$c^k.C_2(3,9)$-geometries satisfying (T) were classified, this completes
the classification of the triangular extensions of classical GQ(3,$t$)'s.

The proof depends upon a computer enumeration of hyperovals in these 
GQ(3,3)'s.
A similar idea has been implemented in Chapter~\ref{chap:epolsp} and in
the author's preprint \cite{Pa:mclco3}.

\paragraph{Chapter~\ref{chap:suz}.}
The {\em Suzuki chain} \cite{Suz} is the following series
$G_1,\dots,G_7$ of finite
groups $$A_4<L_2(3)<U_3(3)<\HJ<G_2(4)<\Suz<\Co_1.$$ The group $G_i$ is
an automorphism
group of a geometry $\cG_i$ of rank $i+1$ with diagram $\cD_i$:
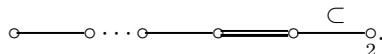
$$\node\arc\node\cdots\node\arc\node\darc\node\stroke{\subset}\node_2.$$
Let $\cG$ be a geometry with diagram $\cD_i$, for some 
$3<i<7$, satisfying (T) and such that its
residues with diagram $\cD_3$ are isomorphic to $\cG_3$.
We show that then $\cG\cong\cG_i$, or $i=6$ and
$\cG\cong 3\cG_6$,  where $3\cG_6$ is a proper triple cover of
$\cG_6$. 
In particular, this result gives the first characterization of
geometries for the sporadic groups $\HJ$ and $\Suz$ not involving
flag-transitivity assumptions. Later, Hans Cuypers \cite{Cuy:suz}
obtained a similar characterization of geometries for the Suzuki chain
groups which are dual to ours in certain sense: their point 
graphs are complements of the point graphs of our geometries, see
Chapter~\ref{chap:EGO} for a more detailed statement of this result
\cite{Cuy:suz}. 

L.~Soicher proved the result of Chapter~\ref{chap:suz}
under the extra assumption that the
automorphism group of $\G$ is flag-transitive. He was also able to
settle the case $i=7$ under this assumption.
Note that the assumption on the $\cG_3$-residues is important, since
there exist geometries with diagram $\cD_i$ obtained via
truncating the (infinite) Coxeter complex of the Coxeter group with
diagram
$$\node\arc\node\cdots\node\arc\node\darc\node\arc\node\arc\node.$$

It is worth mentioning that the classification of geometries
with diagram $\cD_3$ and $\cD_2$-residues isomorphic to $\cG_2$
and satisfying (T) follows from \cite{BFdFSh}.

Finally, note that it follows from the results of
Chapter~\ref{chap:suz} that the presentations for the groups
$\Aut(\HJ)$, $\Aut(G_2(4))$ and $\Aut(\Suz)$ given in \cite{Atl} are
faithful. This is the first computer-free proof of that fact.

\paragraph{Chapters~\ref{chap:fi2n} and \ref{chap:epolsp}.}
This is a contribution to the theory of $c^k.C_n(s,t)$-geomet\-ries in
the case $n>2$ and $t>1$. 
In particular, all such geometries for $k>1$ satisfying (T) and
$c.C_n(s,t)$-geometries have been classified.

It has been noted in \cite{BH} that if $\cG$ is a
$c.C_n(s,t)$-geometry, $n>2$, then $s=2$ or $4$. By
\cite{BH,Bue:extpolarsp}, all the examples with $s=2$ are standard
(possibly improper) quotients of affine polar spaces over GF(2).

The case $s=4$ remained open. With additional assumption on the
existence of a flag-transitive automorphism group, it was shown in
\cite{BH,DFGMP,Mei} that the only example existing is the
$c.C_3(4,2)$-geometry with automorphism group $Fi_{22}$.

In Chapter~\ref{chap:fi2n} we show that the only $c.C_3(4,2)$-geometry
satisfying (T) and an extra technical assumption on the set of common
neighbours of any 4 vertices of the point graph forming a 4-cycle is
the $Fi_{22}$-geometry. 

In Chapter~\ref{chap:epolsp} we show that the only
$c.C_3(4,t)$-geometry with $t>1$ is the $Fi_{22}$-geometry.
The proof depends upon a computer enumeration of hyperovals in rank 3
polar spaces over GF(4).

Let $\cG_k$ be a $c^k.C_n(s,t)$-geometry, $k>1$, $n>2$.
It follows from the classification of extended projective spaces
\cite{Hug65.1,Hug65.2} that 
$n=3$, $s=4$, $k\leq 3$. It was shown in \cite{Mei,vBW:fi:char} that
the only flag-transitive $\cG_k$'s are the geometries for the groups
$Fi_{23}$ ($k=2$), $Fi_{24}$ and $3\cdot Fi_{24}$ ($k=3$). 

In Chapter~\ref{chap:fi2n} we show that any $c^k.C_3(4,t)$-geometry,
$k>1$, $t=2$, satisfying (T), is one of the geometries for Fischer
groups just mentioned. 
It follows from the results of Chapter~\ref{chap:epolsp} that
the condition $t=2$ can be weakened to $t\geq 2$.

As a by-product, in Chapter~\ref{chap:epolsp} we give a construction
of a new $c.C_2(4,2)$-geometry. This is the first example of an
extension of a classical GQ$(s,t)$, $t>1$, not admitting a
flag-transitive automorphism group.

\paragraph{Chapter~\ref{chap:hex}.}
The groups $\McL$ and $\Co_3$ act flag-transitively on geometries
$\cG_k$ of rank $k$, where $k$ is 4 and 5 respectively, 
with the following diagram
$\cD_k$
$$\node\arc\node\cdots\node\stroke{\subset}\node\tarc\node,$$
and the $\cD_3$-residues are isomorphic to an extended generalized
hexagon $\cG_3$ (EGH, for short) with automorphism group $L_3(4)$.
Note that $\G=\cG_3$ satisfies the following property 
\begin{enumerate}
\item[$(*)$] $\{x_1,x_2,x_3\}$ is a clique  of the point graph not contained
in a circle, if and only if $x_2$ and $x_3$ are at distance 3 in the
point graph of $\G_{x_1}$.
\end{enumerate}
In fact, EGHs satisfying $(*)$ were completely classified by Cuypers
in \cite{Cuy:suz}.

We show that if $\cG$ is a geometry with diagram $\cD_k$, $k>3$,
satisfying (T), and its $\cD_3$-residues satisfy $(*)$, then either
$k=4$ and $\cG\cong\cG_4$, or $k=5$ and $\cG\cong 2\cG_5$, where
$2\cG_5$ is the double cover of $\cG_5$.

As well, we show that if we replace (T) by the assumption $k\le 5$, we
get one of the three geometries: $\cG_4$, $\cG_5$, or $2\cG_5$
(note that $\cG_5$ does not satisfy (T)).

R.~Weiss~\cite{We:MCL} proved a similar result under an additional
assumption that $\cG$ admits a flag-transitive automorphism group.
It follows from our result that the presentations for the groups $\McL$
and $\Co_3$ given in \cite{We:MCL} are faithful (whereas in
\cite{We:MCL} Weiss has used a coset enumeration performed on a
computer to prove this). 

Note that an assumption similar to $(*)$ appears to be a necessary one,
since, according to Pasini~\cite{Pa:covers}, extended generalized
hexagons have infinite universal covers.

The main result of the Chapter is one of the two known characterizations of
geometries for the groups $\McL$ and $\Co_3$ not involving
flag-transitivity assumptions. The other result is a classification of
$c^k.C_2(3,9)$-geometries given by the author in \cite{Pa:mclco3}. The
latter has the disadvantage that a computer is used at one stage, which
prompted another characterization to be included in the thesis.

\paragraph{Chapter~\ref{chap:EGO}.}
The group $\He$ acts flag-transitively on an extended generalized
octagon (EGO, for short) $\cH$, that is a geometry with a diagram
$$\node\stroke{\subset}\node\farc\node.$$
Note that $\cH$ satisfies $(*)$ (see the paragraph about
Chapter~\ref{chap:hex} above). 

We show that there exist exactly two nonisomorphic EGO satisfying
$(*)$, namely $\cH$ and its subgeometry admitting $2L_3(4)$ as a
flag-transitive automorphism group.

Unlike the rest of the results included in this thesis, the point
residues of the geometries in question form a rather ``wild'' class,
and it takes quite an effort to establish that $(*)$ brings this class
down to just two objects, namely the flag generalized octagons of the
generalized quadrangles $W(2)$ and $W(4)$.

As well, we show that these EGO do not admit any further extensions.

Again, R.~Weiss~\cite{We:he} proved a similar result under an
additional assumption that $\cG$ admits a flag-transitive automorphism
group.  It follows from our result that the presentation for the group
$\He$ given in \cite{We:he} is faithful (whereas in \cite{We:he} this
follows from a coset enumeration performed on a computer).

Again, an assumption similar to $(*)$ seems unavoidable, since EGO, as
well as EGH, have infinite universal covers, cf. \cite{Pa:covers}.

\chapter{On some locally 3-transposition graphs}
\label{chap:o3n} 
\paragraph{Abstract.}\footnote{A part of this research was
 completed when this author held a position at the 
Institute for System Analysis, Moscow.}
Let $\Si^\ep_n$ be the graph defined on the (+)-points of an 
 $n$-dimen\-sio\-nal
GF(3)-space carrying a nondegenerate symmetric bilinear form with
discriminant $\ep$, points are adjacent iff they are perpendicular. We prove
that if $\ep=1$, $n\geq 6$ (resp. $\ep=-1$, $n\geq 7$)
then $\Si^\ep_{n+1}$ is the unique connected locally
$\Si^\ep_n$ graph. One may view this result as a characterization of a
class of $c^k\cdot C_2$-geometries (or 3-transposition groups).
We briefly discuss an application  of the result to a 
 characterization of
Fischer's sporadic groups.
\section{Introduction and results}
The study of geometries of classical groups as point-line
systems with fixed local structure is very extensive. For instance, see
Tits~\cite{Ti:la}. We refer the reader to a paper of Cohen and Shult~\cite{CoSh}
 for a
brief survey on more recent results.

In this paper we use the aforementioned approach to characterize some
``nonclassical" geometries, namely 3-transposition graphs that arise 
 from
orthogonal GF(3)-groups.

For the concept of  a 
3-transposition graph
and a 3-transposition group, 
  see Fischer~\cite{Fi}. In
case of characteristic 2 groups such graphs come from classical geometries.
It is worthwhile to mention 
 the work of Hall and Shult~\cite{HaSh}
characterizing some class of 3-transposition graphs as {\it locally
cotriangular} ones. Note that the graphs $\Si^\ep_n$ considered in the
present paper are not locally cotriangular (with finitely 
many
exceptions for small $n$). 

There is a one-to-one correspondence between 3-transposition
graphs and {\it Fischer spaces} (see, e.g. Buekenhout \cite{Bue:fi}, Cuypers
\cite{Cuy:genfisp}, Weiss \cite{We1}, \cite{We2}), namely 3-transposition graphs
are complements to the collinearity graphs of Fischer spaces. It turns
out that it is possible to exploit this duality and classification of
Fischer spaces in the final part of the proof of our theorem~\ref{o3n:mainth}.
Note, however, that locally $\Si^\ep_n$ graphs not always arise from
Fischer spaces. E.g., there are at least two nonisomorphic
locally $\Si^-_6$-graphs, one of them has diameter 4, (i.e. it does not
correspond to any Fischer space) see e.g. \cite{BCN}. The author has 
 shown
[in preparation] that they are the only examples of such graphs.

Throughout the paper we consider undirected graphs without loops and multiple
 edges.
Given a graph $\Ga$, let us denote the set of vertices by
$\V=\V\Ga$,  the  set  of edges by $\E=\E\Ga$.
Given two graphs $\Ga,\/\Delta$, the graph $\Ga\cup\Delta$
(resp. the graph $\Ga\cap\Delta$) is
the graph with the vertex set $\V\Ga\cup \V\Delta$ (resp.
$\V\Ga\cap \V\Delta$) and the edge set $\E\Ga\cup \E\Delta$
(resp. $\E\Ga\cap \E\Delta$). Given $v\in \V\Ga$, we denote by
$\Ga_i(v)$ the subgraph induced by vertices at distance $i$ from $v$,
and $\Ga_1(v)=\Ga(v)$.
Furthermore, $\Ga(X)=\bigcap_{x\in X} \Ga(x)$.
To simplify the notation we use
$\Ga(v_1,\dots,v_k)$ instead of $\Ga(\{v_1,\dots,v_k\})$ and
$u\in\Ga_i(\dots)$ instead of $u\in \V\Ga_i(\dots)$.
As usual, $v=v(\Ga)=|\V\Ga|$,  $k=k(\Ga)=v(\Ga(x))$, where
$x\in\V\Ga$. Let $y\in\Ga_2(x)$. We denote
$\mu=\mu(x,y)=\mu(\Ga)=v(\Ga(x,y))$. Of course, 
 we use
$k$, $\mu$ if it makes sense, 
 i.e. if those numbers are independent on the
particular choice of the corresponding vertices.
If $\Delta$ is a (proper) subgraph of
$\Ga$ we denote this fact as $\Delta\subseteq\Ga$ (resp.
$\Delta\subset\Ga$).

We denote  the complete multipartite graph with the $m$ parts of equal size $n$
by $K_{n\times m}$.
$\Aut(\Ga)$ denotes the automorphism group of $\Ga$.
Our group-theoretic notation is as in \cite{Atl}.
Let $\Ga$, $\Delta$ be two graphs. We say that $\Ga$ is
{\it locally} $\Delta$ if $\Ga(v)\cong\Delta$ for any $v\in \V(\Ga)$.
Let $\Ga,\ \overline\Ga$ be two graphs. We say that $\Ga$ is a
{\it cover} of $\overline\Ga$ if there exists a mapping $\varphi$
from $\V\Ga$ to $\V\overline\Ga$ which maps edges to edges.
Suppose we have a chain of graphs $\Si_1,\dots,\Si_n$, such that
$\Si_i$ is locally $\Si_{i-1},\  i=2,\dots,n$. Then for any complete
$k$-vertex subgraph $\Up$ of $\Si_m$, $1<k<m\leq n$
we have $\Si_m(\V\Up)\cong\Si_{m-k}$.
We say that a graph $\Ga$ is a {\it triple} graph if for each nonedge
$(u,v)$ there exist a unique $w\in \V\Ga$ such that
$\Ga(u,v)=\Ga(u,w)=\Ga(v,w)$.

We slightly adopt notation and several basic facts from \cite{Mei}.
Let $T=T_n$ be an 
$n$-dimensional GF(3)-vector space carrying a nondegenerate
symmetric bilinear form (,) with discriminant $\ep$. We say that the point
$\<v\>\subset T$ (or a nonzero vector $v\in T$) is of type $(+)$, $(-)$ or
{\it isotropic} according to $(v,v)=1,\, -1,\, 0 $ respectively. Since the form
 is
constant on a point, the notation like $(p,q)$ for points $p$, $q$ will be
used freely. The orthogonal complement of $X\subseteq T$ in $T$ is denoted by
$X^\perp$.

Define the graph $\Si^\ep_n=\Ga(\V,\E)$ as follows. Let $\V$ be the
set of (+)-points of $T$. Define $\E=\{(u,v)\subset \V\times\V |
(u,v)=0\}$, i.e. the edges are pairs of perpendicular (+)-points.
Given $\Ga=\Si^\ep_n$, we denote by $T(\Ga)$ the underlying
GF(3)-space. Note that $\Si^\ep_n$, $n\geq 5$ if $\ep=1$, $n\geq 4$ if
$\ep=-1$, is a rank 3 graph with  automorphism group
$GO^\mu_n(3)$, where Witt defect $\mu$ is empty if $n$ is odd, otherwise
$\mu=(-\ep)^{n/2}$.
Since $T_n$ can be represented as the orthogonal direct sum of a (+)-point
and $T_{n-1}$, the graph $\Si^\ep_n$ is locally $\Si^\ep_{n-1}$.
We denote $\Si^1_n$ by $\Si^+_n$, and $\Si^{-1}_n$ by $\Si^-_n$.
Note that $\Si^\ep_n$ is a triple graph.
We refer the interested reader to \cite{Mei} for more detailed
information about $\Si^\ep_n$.
We will prove the following theorem. 

\begin{thm} Let $\Th=\Th_{n+1}$ be a connected locally $\Si^\ep_{n+1}$-graph,
where $\ep=1$, $n\geq 6$, or $\ep=-1$, $n\geq 7$.
Then $\Th$ is isomorphic to $\Si^\ep_{n+1}$.
\label{o3n:mainth}
\end{thm}

\subsubsection*{\bf Remark 1} 
One may view the graphs $\Th_{n+1}$ as the collinearity graphs
of 
 certain $c^k\cdot C_2(s,t)$-geometries ${\cal G}(\Th_{n+1})$
(here $k=n-(\ep+9)/2$), i.e. rank $k+2$ geometries
with diagram
$$\circ \!{\phantom{00000}\over \phantom{0}}
\!\!\circ\cdots\circ \!\!{\phantom{00000}\over \phantom{0}}
\!\!\circ \!\!{{\phantom{00} c \phantom{00}}\over \phantom{0}}
{\vphantom{\bigl\langle}}_s
\!\!\!\circ \!\! = \!\! = \!\! = \!\! = \!\! =
\!{\vphantom{\bigl\langle}}_t\!\!\!\circ$$
Conversely, the elements of the geometry may be viewed as $i$-cliques
of the graph with natural incidence, $i=1,2,\dots, k+1, s+k+1$.
Meixner has proved the following result~\cite{Mei}. 

\begin{res}
Let ${\cal G}$ be a residually connected flag-transitive
$c^k\cdot C_2$-geometry. Then if the $c^1\cdot C_2$-residues (resp.
$c^3\cdot C_2$-residues) of ${\cal G}$ are isomorphic to ${\cal G}(\Si^+_6)$
(resp. to ${\cal G}(\Si^-_7)$) then ${\cal G}={\cal G}(\Si^+_{k+5})$
(resp. ${\cal G}={\cal G}(\Si^-_{k+4}))$.
\label{o3n:Meires}
\end{res}

Our theorem  
implies that 
the flag-transitivity
assumption can be replaced to a geometric condition (X) from \cite{Mei}. (X)
states that in the collinearity graph of ${\cal G}$ each $i$-clique is the
shadow of some element, $i=1,2,\dots k+1,s+k+1$.

\subsubsection*{\bf Remark 2}
 The significance of $\Si^+_{i+3}$ as subgraphs of
 3-transposition graphs $\Delta_{2i}$ of Fischer's sporadic simple groups 
 $Fi_{2i}$ ($i=2,3,4$) is well-known.
Namely, for $\Ga=\Delta_{2i}$ the subgraph $\Ga(x,y)$, where $x$ and $y$ are two
 vertices at distance 2, is isomorphic to $\Si^+_{i+3}$. The author
 used theorem~\ref{o3n:mainth} 
 to prove the following result~\cite{Pa:fi}. 

\begin{res}
Any connected locally $\Delta_{22}$ (resp.\ $\Delta_{23}$) graph is isomorphic
to $\Delta_{23}$ (resp. to $\Delta_{24}$ or to its 3-fold antipodal cover).
\qed
\end{res}

\section{Proof of the theorem} 

\paragraph{Preliminaries} A proof of the following technical statement is omitted.

\begin{lem}
Let $\Ga=\Si^\ep_n$,  $n\geq 5$,
and let $a,b,c$ be isotropic 
points of $T=T(\Ga)$. \hfil\break
\indent {If $(a,b)=0$, $a\not=b$, $(a,c)\not=0$ then} \hfil\break
\indent {(i) $a^\perp\cap c^\perp\cap\Ga\cong\Si^{-\ep}_{n-2}$,} \hfil\break
\indent {(ii) $a^\perp\cap b^\perp\cap\Ga$ is not isomorphic to
$\Si^{-\ep}_{n-2}$.}    \hfil\break
\indent {(iii) Denote $\Up_p=\Ga\cap p^\perp$. If
$\Up_a\cap\Up_b$ contains $\{v\}\cup\Up_a(v)$ for some
$v\in\V\Up_a$ then $a=b$.}
\hfil\break
\indent {(iv) For any $u\in\V\Ga\setminus\V\Up_a$, the subgraph
$\Ga(u)\cap\Up_a$ is isomorphic to $\Si^{-\ep}_{n-2}$. }
\label{L4}
\end{lem}

\paragraph{Neighbourhood of two vertices at distance two}
We start with a simple general fact. Let $\Delta$ be a connected
graph satisfying the following property

\begin{description}
\item[(*)] For any $u\in \V\Delta$ and
$v\in \V\Delta\setminus (\V\Delta(u)\cup\{ u \})$
the subgraph $\Delta(u,v)$ is isomorphic to some $M_{\Delta}$, whose
isomorphism type is independent on the particular choice of $u$ and $v$.
\end{description}

\begin{lem}
Let $\Ga$ be locally $\Delta$ graph, where $\Delta$
satisfies {\rm (*)}. Then for  
 any $u\in \V\Ga$, $v\in\Ga_2(u)$
the graph $\Ga(u,v)$ is locally $M_{\Delta}$. \qed
\label{L6}
\end{lem}

Note that $\Ga=\Si^\ep_n$ satisfies (*), when $\ep=1$, $n\geq 5$
(resp. $\ep=-1$, $n\geq 4$). Indeed, the stabilizer
of $u\in \V\Ga$ in $\Aut(\Ga)$ acts transitively on
$\Ga_2(u)$. This implies (*).
Thus lemma~\ref{L6} 
 holds for locally $\Ga$ graphs. The next statement
characterizes locally $M_\Ga$-subgraphs of $\Ga$.

\begin{pr}
Let $\Ga=\Si^\ep_n$, and either $\ep=1$, $n\geq 6$,
or $\ep=-1$, $n\geq 7$. Let $\Om\subset\Ga$ be a locally $M_\Ga$ graph.
Then there exists a unique isotropic point $p\subset T(\Ga)$ such that
$\Om=\Ga\cap p^\perp$.
\label{P7}
\end{pr}

\begin{pf}
We proceed by induction on $n$. It is straightforward to check its basis.
We leave it to the reader, noticing that for
$\ep=1$, $n=6$ (resp. $\ep=-1$, $n=6$) it suffices to classify locally
$K_{3\times 2}\cup K_{3\times 2}\cup K_{3\times 2}$
(resp. $K_{3\times 4}\cup K_{3\times 4}\cup K_{3\times 4}$) subgraphs of
$\Ga$.

Now let us check the inductive step. Let $\Om$ be a connected component
of a locally $M_\Ga$ subgraph of $\Ga$, $v_1\in\V\Om$,
$v_2\in\Om(v_1)$. By the inductive hypothesis,
$\Om(v_i)=\Ga(v_i)\cap p^\perp_i$, where $p_i$ is an
isotropic point of $T(\Ga(v_i))$, hence of $T(\Ga)$ ($i=1,2$).
Since $p_2\subset v^\perp_1$, it defines a locally $M_{\Ga(v_1)}$
subgraph $\Up$ of $\Ga(v_1)$, and $\Up(v_2)=\Om(v_1,v_2)$.
Hence by lemma~\ref{L4}~(iii), 
 applied to $\Ga(v_1)$, $p_1=p_2$.
Therefore $\Om=\Ga\cap p^\perp_1$. Finally, it is easy to check that
$\Om$ is a unique connected component of the subgraph under consideration.
\end{pf}

\paragraph{Final part of the proof} Let $\Th$ be a connected locally
$\Ga=\Si^\ep_n$ graph. Here we assume either $\ep=1$, $n\geq 6$ or
$\ep=-1$, $n\geq 7$. Pick a vertex $u\in\V\Th$.

\begin{lem}
\hfil\break
\indent (i) $\mu(\Th)=\mu(\Si^\ep_{n+1})$, $v(\Th)=v(\Si^\ep_{n+1})$.\hfil\break
\indent (ii) $\Th$ is a triple graph.
\label{L8}
\end{lem}

\begin{pf}
(i). The first claim follows from proposition~\ref{P7}.
Indeed, counting in two ways
the edges between $\Th(u)$ and $\Th_2(u)$, we obtain the precise value of the
number of vertices of $\Th$. Let $v\in\Th_2(u)$, $w\in\Th(v)\setminus\Th(u)$.
By lemma~\ref{L4}~(iv) we obtain $\Th(u,v,w)\cong\Si^{-\ep}_{n-2}$, hence
 nonempty.
Thus the diameter of $\Th$ equals two, and we are done.

(ii). Assume that there exist three distinct vertices $v_i\in\Th_2(u)$,
$i=1,2,3$ such that $\Up=\Th(u,v_1)=\Th(u,v_2)=\Th(u,v_3)$.
It contradicts the fact that $\Th(w)$, where
$w\in\V\Up$, is a triple graph.
Observe that $|\Th_2(u)|$ is exactly twice the number of isotropic points
in $T(\Ga)$. Hence we have no choice determining the edges between $\Th(u)$
and $\Th_2(u)$. Now since $\Si^\ep_{n+1}$ is a triple graph, the same is
true for 
 $\Th$.
 \end{pf}

Let us denote $\Ga=\Th(u)$, $\Xi=\Th_2(u)$. Let $\Up$ be the graph defined
on isotropic points of $T(\Ga)$, two points are adjacent if they are
not perpendicular.

\begin{pr}
$\Xi$ is a two-fold cover of $\Up$.
\label{P9}
\end{pr}

\begin{pf}
Let $(v_1,v_2)\in\E\Xi$. By lemma~\ref{L4}~(iv) we have $\Th(u,v_1,v_2)\cong
\Si^{-\ep}_{n-2}$. Denote by $p_i$ the isotropic point of $T(\Ga)$
such that $p_i^\perp\cap\Ga=\Th(u,v_i)$ ($i=1,2$). By proposition~\ref{P7} such
 a point
$p_i$ exists and is unique ($i=1,2$). By lemma~\ref{L4}~(i),
(ii), $(p_1,p_2)\not=0$.

Conversely, assume $p_1$, $p_2$ are nonperpendicular isotropic points of
$T(\Ga)$.  Denote by $\Om_i$ the locally $M_\Ga$ subgraph of
$\Ga$, which corresponds to $p_i$ ($i=1,2$).
For each $\Om_i$ we have exactly two vertices
$v_{ij}\in\Th_2(u)$ such that $\Om_i=\Th(u,v_{ij})$ ($i,j=1,2$).
Considering the neighbourhood of $v_{11}$, we see that
both $v_{21}$ and $v_{22}$ cannot be ajacent to $v_{11}$.
On the other hand $k(\Xi)=k(\Up)$. Hence one of $v_{21}$ and $v_{22}$ must
be ajacent to $v_{11}$. We have shown that the mapping
$v\mapsto\Th(u,v)$, where $v\in\V\Xi$, is a covering from $\Up$ to $\Xi$.
\end{pf} 

The latter statement implies that $\Xi$ possesses 
 an involutory automorphism
$g_u$ which interchanges any $v,w\in\V\Xi$ such that $\Th(u,v)=\Th(u,w)$, and
fixes $\{u\}\cup\Th(u)$ pointwise.

Consider the subgroup $G_u$ of $\Aut(\Th)$ generated by $g_x$, $x\in\Th(u)$.
We have $G_u\cong\Aut(\Ga)$. Therefore $\Th$ is the collinearity graph of a
$c^k\cdot C_2$-geometry satisfying the conditions of
result~\ref{o3n:Meires}. 
Hence our result follows from it. However, we would like to give
a complete proof of theorem~\ref{o3n:mainth} here.

Our claim is that $\{g_v|v\in\V\Th\}$ is a class of 3-transpositions in
$\Aut(\Th)$. Indeed, clearly for any $x\in\Th(u)$ the involutions $g_u$ and
$g_x$ commute. Now let $y\in\Th_2(u)$. We must prove $\tau=(g_u g_y)^3=1$.
Note that $\tau$ belongs to the kernel of the
action of the stabilizer of every $v\in\Th(u,y)$ on $\Th(u)$. Therefore
it fixes every vertex of $\Th$.
Our claim is proved. The use of the classification of 3-transposition groups
given in \cite{Fi} completes the proof of theorem~\ref{o3n:mainth}.
\qed 

\subsubsection*{\bf Note} 
Hans Cuypers (personal communication) has suggested another idea how to
complete the proof, which is much more geometric. Namely, it may be
easily
shown that the partial linear space on $\V\Th$, whose lines are triples,
is an irreducible Fischer space (see introduction) 
 (or a locally polar
geometry with affine planes, see Cuypers and Pasini \cite{CuPa}). Then
the use of classification of these objects \cite{Cuy:genfisp} (resp.
\cite{CuPa}) completes the proof.

\chapter{The triangular extensions of a 
generalized quadrangle of order $(3,3)$}
\label{chap:EGQ33}
\paragraph{Abstract.}
We show that the triangular extension of a generalized quadrangle of order
(3,3) is unique. The proof depends upon certain computer calculations.
\section{Introduction and the result}
Extensions of finite generalized quadrangles (EGQ, for short), 
or, more generally, of polar
spaces, play a important role as incidence geometries admitting
interesting automorphism groups, such as sporadic simple, or some classes
of (extensions of) classical groups.  Buekenhout and Hubaut \cite{BH}
initiated the study of extensions of polar spaces from geometric point
of view by proving some characterization theorems, in particular they
classified locally polar spaces such that the lines of the residual
polar space are of size 3. They also classified locally polar spaces
admitting classical group acting on point residues,  later on these
results were generalized in a more general framework of flag-transitive
diagram geometries, see a survey \cite{PaYo} by Pasini and Yoshiara for
an extensive bibliography.  However, very few 
characterizations are known which do not assume group actions.
Blokhuis and Brouwer \cite{BlBr44} and P.Fisher \cite{Fi44} classified
EGQ(3,1), 
apart from the already mentioned result in \cite{BH} for polar lines of
size 3, and generalizing it result \cite{Bue:extpolarsp} by Buekenhout. 
The author \cite{Pa:fi,Pa:3tr} characterized extensions
of polar spaces related to some 3-transposition groups, including
Fischer's sporadic simple groups.  In \cite{Pa:mclco3} he proved the
uniqueness of EGQ(3,9) and classified its further extensions.

Here we shall be concerned with triangular EGQ(3,3). 
For basic definitions and general account
on EGQ see Cameron, Hughes and Pasini \cite{CHP}. A triangular EGQ$(s,t)$
may and will be viewed as a graph $\G$ such that the subgraph $\G(u)$
induced on the neighbourhood of any vertex $u$  is isomorphic to the
collinearity graph of a generalized quadrangle of
order $(s,t)$, or GQ$(s,t)$, for short. 
Concerning GQ$(s,t)$, a standard reference is \cite{PayTh}.
There are two nonisomorphic GQ(3,3),
dual to each other. One, usually denoted $W(3)$, is the point-line
system of the isotropic points and totally isotropic with respect to a
nondegenerate symplectic form lines of the 3-dimensional projective
space over GF(3) (PG(3,3), for short). The other one, usually denoted 
$Q_4(3)$, is defined similarly with a
nondegenerate symmetric bilinear form in place of the symplectic one
and PG(4,3) in place of PG(3,3). 

Let ${\cal U}_n$ be the graph defined on the nonisotropic points on a
$n$-dimensional GF(4)-space $T=T({\cal U}_n)$ carrying a nondegenerate hermitian form,
points are adjacent if they are perpendicular. Note that ${\cal U}_4$ is
isomorphic to the collinearity graph of $W(3)$, and ${\cal U}_{n+1}$ is
locally ${\cal U}_n$. 

We say that a graph $\G$ is locally ${\cal D}$ (or $\D$), 
where $\D$ is a family of graphs (resp. $\D$ a graph), 
if for any vertex $u$ of $\G$ the subgraph
$\G(u)$ is isomorphic to a member of ${\cal D}$ (resp. to $\D$).
In our case ${\cal D}$ consists of the collinearity graphs of $W(3)$ and
$Q_4(3)$.

\begin{thm}
Let $\G$ be a triangular {\rm EGQ(3,3)} (in other words, $\G$ is locally
${\cal D}$). Then $\G$ is isomorphic to ${\cal U}_5$.
\label{EGQ33}
\end{thm}

Under the additional
assumption that a classical group is induced on $\G(u)$ for
any vertex $u$, this statement was proved in \cite{BH}.
Later on, this was improved in \cite{DFGMP,Yo1} in the slightly more
general framework of the classification of flag-transitive
$c.C_2$-geometries, still involving a strong assumption on group action. 

\section{Proof}
Let $\G$ be a triangular EGQ(3,3).
Our approach is based on the observation made in
\cite{BH} that given a point $u$ of $\G$ and a point $v$ at distance 2
from $u$, its common neighbourhood $\G(u,v)$ (which will be often
called a {\em $\mu$-graph of $\G(u)$}), is a {\em hyperoval} (or a
{\em local subspace}, in the terminology of \cite{BH})
in $\D=\G(u)$, that is a subset $\Phi$ of points of $\D$ such that each
line of $\D$ meets either 0 or 2 points of $\Phi$. 
Hyperovals of GQ were studied by several authors, see e.g.
\cite{DFGP,Pas:02GQ}. However the results achieved are concerned mainly
with various extreme cases, and nothing like a classification of the
hyperovals in GQ, which is required in our approach, exists. 

So we classify the hyperovals of $Q_4(3)$ and $W(3)$, using a computer. 
Then we rule out most of the 
hyperovals of $Q_4(3)$, using some simple criteria. 
As an immediate corollary we have that $\G(x)\cong\G(y)$
for any distance two pair of points $x$, $y$ of $\G$ such that
$|\G(x,y)|$ is of certain size. It gives us an opportunity to
eliminate most of the hyperovals of $W(3)$.

The remaining ones are exactly the 45 hyperovals as appeared in ${\cal U}_5$.
Then we assume that $\G(u)=\D\cong W(3)$. 
We deduce that $\G$ is a strongly regular graph having the same parameters 
as ${\cal U}_5$. Then we establish that $\G(x)\cong\D$ for {\em any} $x\in\G$.
Moreover, we see that $\G$ has {\em quadruples}, that is, for each
nonadjacent pair of points $x$, $y$ there are exactly two other points
$z$, $z'$ such that $\G(x,y)=\G(x,y,z,z')$. This defines on $\G$ the
structure of a partial linear space with line size 4. 
One can then check that the latter partial
linear space such that the lines and the affine planes on any point form
a finite GQ, which were
classified in \cite{CuyEGQ}. The application of \cite{CuyEGQ,CoSh} 
completes the proof that $\G\cong{\cal U}_5$ (alternatively, we
demonstrate how to use the classification of generalized Fisher spaces
\cite{Cuy:genfisp,CuySh} to get the same result). 

The remaining case, where $\G(x)\cong Q_4(3)$ for any $x$,
is dealt with similarly. It turns out that this assumption leads to a
contradiction.

\subsection{Preliminaries}
The determination of all the hyperovals in a GQ(3,3) $\D$ was based on an
almost straightforward backtrack exhaustive search. 
It is natural to regard a hyperoval $\O$ as the subgraph induced by $\O$
in $\D$. Clearly if $\O$ is disconnected then each connected component
of $\O$ is a hyperoval, as well. So we look for the connected hyperovals
only, and then, if possible, glue components together. We note, however,
that all the hyperovals in $\D$ turn out to be connected.
The main way to
reduce the number of objects found by the search was the use of
the group $G=Aut(\D)$. Indeed, for the set of  orbits of $G$ on the
hyperovals of $\D$ it suffices to find a representative $R_k$ for each orbit
$O_k$.
Moreover, the following idea proved to be highly successful. 

Let $S$ be a graph which is a subgraph of $\O$ for any hyperoval
$\O$ of $\D$ (for instance, $S\cong K_2$).  Let ${\cal S}=\{S_j\}$ be a set of
representatives of the orbits of $S$ on the subgraphs of $\D$ isomorphic
to $S$. Then  a set ${\cal R}=\{R_k\}$ of representatives of $G$-orbits on the
hyperovals may be chosen in such a way that each $R_k\in {\cal R}$ contains some
$S_j\in {\cal S}$. 

So the problem is to find such $S$ that the set ${\cal S}$ is not huge
and, on the other hand, the number of different hyperovals containing
given $S_k\in {\cal S}$ is not huge, as well. Since any $\O$ is a
triangle-free graph of valence 4, we choose as $S$ the 4-{\em claw}, that is
the subgraph induced on the union of $\{x\}$ and $\O(x)$.
We find a set ${\cal S}$. Then for each $S_k\in {\cal S}$ we
perform the exhaustive search of the hyperovals containing $S_k$.
The resulting set ${\cal R}$ need not be a minimal one, that is, it may
contain several representatives for one $G$-orbit.
Thus, finally, we construct a minimal set ${\cal R'}$.

The computer calculations were carried out using {\sf GAP} system for
algebraic computations \cite{GAP} along with the package GRAPE
for computations in
graphs \cite{GRAPE}. The latter uses the package NAUTY for computations of
automorphisms and isomorphisms of graphs \cite{NAUTY}.
\begin{pr}
The hyperovals $\Phi$ of $\D\cong{\cal U}_4\cong W(3)$ are as follows.
\begin{enumerate}
\item $432$ of size $20$. The $20$ points outside $\Phi$ are collinear to 
$8$ points of $\Phi$.
\item $540$ of size $16$.
\item $720$ of size $12$. There are $2$ points outside $\Phi$ collinear to $6$
points of $\Phi$.
\item $45$ of size $8$. 
The $32$ points outside $\Phi$ are collinear to $2$ points of
$\Phi$. 
Let $v\in\Phi$ and $\Phi'$ be a hyperoval containing $v$. 
Then $\Phi(v)=\Phi'(v)$ implies $\Phi'=\Phi$. Let $\Phi'$ be a hyperoval of
type 2. Then $\Phi\cap\Phi'\not\cong K_2$. There is one-to-one correspondence
between the hyperovals of size 8 and the isotropic points of $T=T(\D)$
given by $\Phi=\D\cap p^\perp$, where $p\in T$ is an isotropic point.
\end{enumerate}
The group $Aut(\D)$ acts transitively on the hyperovals of each
type.
\label{lsW3} \qed
\end{pr}
Note that 1) above contradicts the first (technical) part of the statement of 
\cite[Proposition~8]{BH}. Note that the second part of the statement
remains valid, as we shall see later.
\begin{pr}
The hyperovals $\Phi$ of $\D\cong Q_4(3)$ are as follows.
\begin{enumerate}
\item $1080$ of size $14$. There are $8$ points outside $\Phi$ collinear to 
$6$ points of $\Phi$. There are only $4$ hyperovals intersecting 
$\Phi$ in $3K_2$.\label{hs14}
\item $360$ of size $18$. There are $12$ points outside $\Phi$ collinear
to $6$ points of $\Phi$. All the hyperovals intersecting $\Phi$ in $3K_2$
are of type \ref{hs14}.\label{hs18}
\item $324$ of size $20$. The points outside $\Phi$ are collinear
to $8$ points of $\Phi$. All the hyperovals intersecting $\Phi$ in $4K_2$
are of type \ref{hs18}. \label{hs20}
\item $135$ of size $16$. There are $16$ points outside $\Phi$ collinear to $4$
points of $\Phi$, and the remaining $8$ points are collinear to $8$
points of $\Phi$. The hyperovals intersecting $\Phi$ in $4K_2$ are of type
\ref{hs14} or \ref{hs18}. \label{hs16}
\item $216$ of size $10$. There are $20$ points outside $\Phi$ collinear to $2$
points of $\Phi$, and the remaining $10$ points are collinear to $4$
points of $\Phi$. There are $60$ (resp. $20$) hyperovals of type \ref{hs10}
(resp. of type \ref{hs14}) intersecting $\Phi$ in $K_2$.
\label{hs10} 
\item $270$ of size $12$. There are $24$ points outside $\Phi$ collinear
to $4$ points of $\Phi$, and the remaining $4$ points are not collinear to the
points of $\Phi$. There are exactly $24$ hyperovals of type \ref{hs12}
intersecting $\Phi$ in $2K_2$. \label{hs12} 
\end{enumerate}
The group $Aut(\D)$ acts transitively on the hyperovals of each
type.\qed
\label{hoQ43}
\end{pr}
\begin{lem}
In the notation of Proposition \ref{hoQ43}, only hyperovals of types
\ref{hs10} or \ref{hs12} may appear as $\mu$-graphs of EGQ.
\end{lem}\label{hoQ43:poss}
\begin{pf}
It follows from Proposition \ref{hoQ43}, 1, that any hyperoval $\Phi$ 
of type \ref{hs14} cannot appear as $\mu$-graphs. 
Indeed, there must be at least  $8$ hyperovals intersecting $\Phi$ in
$3K_2$, but there are only $4$ such ones.

Next, if $\Phi$ is of type \ref{hs18} there must be other $\mu$-graphs
intersecting $\Phi$ in $3K_2$, but all of them, by Proposition
\ref{hoQ43}, must be of type \ref{hs14}, a contradiction.

Similarly, we reject hyperovals of types \ref{hs20} and \ref{hs16}.
\end{pf}

Now the following statement is immediate.
\begin{cor}
Let $x$, $y$ be two points of $\,\G$ at distance $2$ such that\\
$|\G(x,y)|=20$. Then $\G(x)\cong\G(y)\cong W(3)$.\qed
\label{mu20}
\end{cor}

\subsection{A point of $W(3)$-type exists}
Here we assume that
$\G(u)=\D\cong W(3)$ for some $u\in\G$. 
\begin{lem}
In the notation of Proposition \ref{lsW3}, only hyperovals of types 2 or
4 may appear as $\mu$-graphs.
\label{2or4}
\end{lem}
\begin{pf}
First, we show that the hyperovals of type 1 cannot appear as
$\mu$-graphs. Let $\Phi=\G(u,v)$ be type 1 hyperoval. By Corollary
\ref{mu20}, we have $\G(v)\cong W(3)$.  
The following facts obtained by means of a computer
will be used.
 
1) There are two orbits $O_1$, $O_2$ in the action of the stabilizer $H$
of $\Phi$ in $G=Aut(\D)$ on the set of edges of $\Phi$ 
of lengths 30 and 10, respectively.

2) There are 25 hyperovals intersecting $\Phi$ in the disjoint union
of 4 copies of $K_2$. All of them are of type 1.
The group $H$ has two orbits $\O_1$, $\O_2$
on this set of hyperovals of lengths 20 and 5, respectively.
If $\Psi\in\O_1$ then it contains exactly one edge from $O_2$.
If $\Psi\in\O_2$ then all the edges from $\Psi\cap\Phi$ belong to $O_2$.
Given $e\in O_1$, there exist exactly two hyperovals $\Psi\in\O_1$
such that $e\subset\Psi$.

By Proposition \ref{lsW3}
there are 20 vertices in $\Xi=\G(v)\sm\G(u)$ such that for each
$x\in\Xi$ we have that $\G(x)\cap\Phi$ is the disjoint union of 4 copies
of $K_2$. Hence $\G(u,x)$ belongs to $\O=\O_1\cup\O_2$. So we have
constructed a mapping $\phi$ from $\Xi$ to $\O$. It is an injection,
since the
disjoint union of 4 copies of $K_2$ determines four of the lines of
$\D$, and $\D$ is a partial linear space.
We have to choose 20 of the 25 elements of $\O$. Since the line size of
GQ(3,3) is 4, for each edge $e$ of $\Phi$ there exist exactly two $\Psi\in\phi(\Xi)$
such that $e\subset\Psi$. Now it follows from 2) that $\O_1\subseteq\phi(\Xi)$.
Hence $\phi(\Xi)=\O_1$.  

Now the graph $\G_2(u)$ is isomorphic to a connected component of the
graph whose vertex set is the set $\Phi^G$ (i.e. the hyperovals of type
1), and the edge set is $\{\Phi^G,\Psi^G\}$, where $\Psi\in\O_1$. It is easy
to check either by computer or exploiting its $G$-invariance 
that the latter graph is connected. This is a
contradiction, since the latter graph has 432 vertices, whereas
$\G_2(u)$ has 54. Thus $\Phi$ cannot be of type 1.

Let $\Phi$ be of type 3. There exist a two-element subset $\Xi$ of 
$\G(v)\sm\G(u)$ such that $\G(x)\cap\Phi$ is the disjoint union of 3
copies of $K_2$ for any $x\in\Xi$. On the other hand, computer
calculations show that each hyperoval $\Psi$ such that $\Phi\cap\Psi$
is the disjoint union of 3 copies of $K_2$  is of type 1. This is the
contradiction.
\end{pf}

The following is well known.
\begin{lem} Let $\Xi$ be a {\rm GQ(3,3)}. Assume that $\Xi$ has a local 
hyperoval of type 4, that is, isomorphic to {\rm GQ(1,3)}. 
Alternatively, assume that $\Xi$ possess a
triple $x$, $y$, $z$, of the pairwise not collinear points such that 
$\G(x,y)=\G(x,z)$. Then $\Xi\cong W(3)$.\qed
\label{GQ13inGQ33}
\end{lem}
\begin{lem} There exists $v\in\G_2(u)$ such that $\G(u,v)$ is of type 4.
\label{type4exist}
\end{lem}
\begin{pf}
Suppose that this is not the case. Hence, by Lemma \ref{2or4}, for any
$v\in\G_2(u)$ the hyperoval $\G(u,v)$ is of type 2. We have
$$|\G_2(u)|=40\cdot 27/16,$$ 
which it not integer. This is a contradiction.
\end{pf}

By Lemma \ref{type4exist}, there exists $v\in\G_2(u)$ such that 
$\G(u,v)$ is of type 4.
By Lemma \ref{lsW3}, 4, for each $w\in\G(v)\sm\G(u)$ the subgraph
$\G(u,w)$ is a hyperoval of type 4. Note that $|\G(v)\sm\G(u)|=32$
coincides with the number of the hyperovals intersecting $\O$ in a
$K_2$, and all such hyperovals are of type 4. The graph $\S$ defined on the
type 4 hyperovals by the rule that two vertices are incident if the
intersection of the corresponding hyperovals equals $K_2$ is isomorphic
to the complement of the collinearity graph of GQ(4,2), in particular it
is connected. Hence for any $x$ belonging to the connected component $\Xi$ of
$\G_2(u)$ containing $v$ one has $\G(u,x)$ is of type 4. 

Clearly $\Xi$ is a cover of $\S$. Hence {\em each} type 4 hyperoval of $\G(u)$
is a $\mu$-graph of $\G$. Since
$\mu(\S)=24$, which is greater than the size of any of hyperovals, 
$\Xi$ is a connected proper cover of $\S$. It implies that
for any $x\in\Xi$ there exists $y\in\Xi$ such that $\G(u,x)=\G(u,y)$.
Let $w\in\G(u,x)=\O$. Since $\O$ is a graph of valence 4, we have for
$\Theta=\G(w)$ that $\Theta(u,x)=\Theta(u,y)$. Note that $x$ is not
adjacent to $y$. Hence by Lemma
\ref{GQ13inGQ33}, $\Theta\cong W(3)$. 
Therefore there exists fourth point $t\in\Theta$ not collinear to $u$,
$x$, or $y$, such that $\Theta(u,x)=\Theta(u,t)$. 
By Lemma \ref{lsW3} 4),   $\G(u,t)=\G(u,x)$. Thus for {\em any} type 4
hyperoval $\O$ of $\D$ there exist {\em three distinct} points $x_1$,
$x_2$, $x_3$ of $\G_2(u)$ such that $\G(u,x_i)=\O$ for $i=1,2,3$.
Now, counting in two ways the edges between $\G(u)$
and $\G_2(u)$, one has that for any $x\in\G_2(u)$ the hyperoval $\G(u,x)$
is of type 4. In particular, $|\G_2(u)|=|\G|-40-1=135$.

Note that for any $z\in\D$ we have $\G(z)\cong W(3)$. Since
there exists $p\in\G_2(z)$ such that $\G(z,p)$ is of
type 4, the repetition of the above argument gives us that  for any
$q\in\G_2(z)$ the hyperoval $\G(z,q)$ is of type 4. Therefore $\G$ is a
strongly regular graph (SRG, for short)
with the same parameters as ${\cal U}_5$. 

To summarize, we state the following.
\begin{pr}
Let $\G$ be a triangular {\rm EGQ(3,3)} such that $\G(u)\cong W(3)$ for
some $u\in\G$. Then $\G(v)\cong W(3)$ for any $v\in\G$. Moreover,
$\G$ is an {\rm SRG} with the same parameters as ${\cal U}_5$,
and for any distance two pair of vertices $x$, $y$ there exist two more
vertices $w=w_{xy}$, $w'=w'_{xy}$ such that $\G(x,y)=\G(x,w)=\G(x,w')$.
\qed
\end{pr}

Using the above mentioned quadruples of points of $\G$, we may
define the structure of a partial linear space $\cL$ of line size 4 on $\G$.
The idea to consider $\cL$ is due to Hans Cuypers \cite{Cuyps}.

To complete the current case it suffices to show that $\cL$ is a (finite)
locally GQ with affine planes. It means the following. Consider the set
of minimal-by-inclusion subspaces generated by the pairs of intersected
lines of $\cL$ (such subspaces are usually called {\em planes}). 
We require all the planes which are linear spaces to be affine planes.
Now the incidence system of lines and affine planes through a given point
should be isomorphic to a (finite) GQ. 

\begin{pr} 
Let $\G$ be a triangular {\rm EGQ(3,3)} such that $\G(u)\cong W(3)$ for
some $u$. Then $\cL=\cL(\G)$ is a locally {\rm GQ(4,2)} with affine
planes, and therefore $\G\cong{\cal U}_5$.
\end{pr}

\begin{pf}
Let $w\in\G$, $l=wu$, $m=wv$ are two distinct lines of $\cL$
on $u$ such that 
$w$ and $v$ are nonadjacent, $p(l),p(m)$ the corresponding
isotropic points (cf. Proposition \ref{lsW3}) of $T(\G(w))$. 


Let $x\in l$, $y\in m$ such that $x$ and $y$ are not equal to $w$. Since
$x$ and $y$ lie at distance 2 in $\G$, there is a line through
them. Where are the other points, say $z$, on this line? Clearly $z\not=w$.
We claim that $z\not\in\G(w)$, as well. Suppose false. Then
by Proposition \ref{lsW3} $z$  is 
adjacent to some vertex $t$ of $\G(w,x)$. Since $\G(t)$ is a subspace of
$\cL$ and $xz$ a line of this subspace, we have $y\in\G(t)$, a
contradiction to $\G(w,x,y)=\emptyset$. So $z\in\G_2(u)$, moreover
$\G(w,x,z)=\G(w,y,z)=\emptyset$. Hence $p(z)$ belongs to the totally
isotropic line of $T(w)$ generated by $p(x)$ and $p(y)$.
We have shown that all the points in the subspace of $\cL$ generated by
$l$ and $m$ and 
distinct from $w$ belong to the set $\Pi=\{t\in\G_2(w)\mid p(t)\perp
p(x),\,p(t)\perp p(y)\}$. 
Clearly, $\Pi$ generates a linear space on 16 points with line size 4, that is
the affine plane of order 4. On the other hand, $l$
and $m$ generate a subspace of $\langle\Pi\rangle$. But two lines in
this plane generates it. So $\langle l,m\rangle\cong AG(2,4)$.

It is easy to check that a pair of lines through $w$
corresponding to the pair of nonorthogonal isotropic points of $T(\G(w))$
does not generate a linear subspace of $\cL$. Thus the only planes of $\cL$ are
affine ones. Moreover, the planes on $w$ are in the one-to-one
correspondence to the totally isotropic lines of $T(\G(w))$.
So $\cL$ is a locally GQ(4,2)  with affine planes.

The rest of the proof is straightforward and consists of implementation
of the corresponding classification result \cite{CuyEGQ}, which says 
that the locally (finite) GQ with affine planes 
are in one-to-one correspondence
with the standard quotients of the corresponding affine polar spaces,
see \cite{CoSh,CuPa}. It is easy to
deduce from \cite{CoSh,CuyEGQ} that our
object is the desired ${\cal U}_5$.
\end{pf}

\paragraph{Remark.} Alternatively, it suffices to
show that $\cL$ is a generalized Fischer space, and then apply their
classification. A generalized Fischer space is a connected partial linear space
such that each its plane is either affine or dual affine. 
Also, we know that there are only finitely many lines (points) on a given point
(resp. line) and the cocollinearity graph of $\cL$ (that is, $\G$)
is connected. Such spaces are called finite and irreducible.
Let us prove that $\cL$ is a generalized Fischer space. 
It was shown above that
if $l=wu$, $m=wv$ are two distinct lines on $w$ such that 
$u$ and $v$ are not adjacent in $\G$, then $\langle l,m\rangle\cong AG(2,4)$.
Thus it suffices to prove that if  
$u$ and $v$ are adjacent, then $\langle l,m\rangle\cong AG^*(2,3)$.
Indeed, in this case we may pick $t\in\G(w,u,v)$ and consider $\G(t)$ as
a subspace of $\cL$. Since $\langle l,m\rangle\cap\G(t)$ is   isomorphic
to $AG^*(2,3)$, so is $\langle l,m\rangle$. 

Now it easily follows from the classification of finite, irreducible
Fischer spaces given in \cite{Cuy:genfisp,CuySh} that $\G\cong{\cal U}_5$. 

\subsection{The remaining case}

Here we assume that $\G(u)\cong Q_4(3)$ for any point $u$ of $\G$.
It follows from Lemma \ref{GQ13inGQ33} 
that $\G$ has distinct $\mu$-graphs, that is, if
$\G(x,y)=\G(x,z)$, where $y$ and $z$ both at distance two from $x$, then
$y=z$.

We claim that hyperovals of type \ref{hs10} must appear as $\mu$-graphs of EGQ.
Assume to the contrary that all the $\mu$-graphs are arising from
hyperovals of type \ref{hs12}. By standard counting argument, there
must be 90 such $\mu$-graphs.
It follows from Proposition \ref{hoQ43} that there are exactly 24 $\mu$-graphs
intersecting the given one in $2K_2$, and all with such intersection
must be taken from the set $X$ of type \ref{hs12} hyperovals. 
So we may define an $O_5(3)$-invariant
graph on $X$ of valence $24$, two vertices are adjacent if the
corresponding hyperovals are intersecting in $2K_2$. A union of
connected components of this graph must be of size 90.
But the maximal possible size of the blocks of imprimitivity of $O_5(3)$ on
$X$ is 10. It implies a contradiction, since the valence of $X$ is
greater than 10. We are done.

Finally, we prove that, to the contrary,   hyperovals of type \ref{hs10} 
cannot appear as $\mu$-graphs of EGQ. We need one more statement
rectifying the embedding of a type \ref{hs10} hyperoval $\Phi$ in $\D$.
\begin{lem}
Let $\Phi\subset\D$, $X$ be the set of vertices of $\D$ outside $\Phi$ such
that each vertex $x\in X$ satisfies $\D(x)\cap\Phi\cong K_2$.
$\Phi$ is isomorphic to the two-fold antipodal cover of $K_5$ with the
antipodal equivalence relation $\phi=\phi_\Phi$.
$X$ is a connected graph of valence $6$, $x,y\in X$ are adjacent if
$\D(y)\cap\Phi\cap\phi(\D(x)\cap\Phi)$ is of size one.\qed
\label{adj:hs10}
\end{lem}
Now assume that given $u\in\G$, $v\in\G_2(u)$, we have
$\G(u,v)$ of type \ref{hs10}. Let $W\subset\G(v)\setminus\G(u)$ such
that $\G(u,v,w)\cong K_2$ for any $w\in W$. The subgraph $W$ is
isomorphic to $X$ defined in Lemma \ref{adj:hs10}.
Then $\G(u,w)$ is of type \ref{hs10}, by Proposition \ref{hoQ43} and
Lemma \ref{hoQ43:poss}
The set $Y$ of type \ref{hs10} hyperovals of $\G(u)$ intersecting $\G(u,v)$
in $K_2$ is of size 60 (Proposition \ref{hoQ43}), 
and the stabilizer of $\G(u,v)$ in $O_5(3)$ acts transitively
on $Y$. So in our attempt to select 20 of them
we could start from any element of $Y$. Let $\Phi_1$ be such a hyperoval. 
We try to form a graph isomorphic to $X$ defined in Lemma
\ref{adj:hs10}. There are exactly 6 elements $\Phi'$ of $Y$ such that 
$\Phi'\cap\G(u,v)\cap\phi_{\G(u,v)}(\Phi_1\cap\G(u,v))\cong K_1$ and
$\Phi_1\cap\Phi'\cong K_2$ or $2K_2$ (the latter is necessary condition for
the vertices in $W$ associated with $\Phi_1$ and $\Phi'$ to be adjacent, see
Proposition \ref{hoQ43}, for the former see Lemma \ref{adj:hs10}).
Since $X$ is of valence 6, 
we are forced to pick up all the 6 possible elements of $Y$. 
Proceeding in this manner (i.e.
considering $\Phi'$ instead of $\Phi_1$, etc.), we
however do not end up with 20 elements of $Y$, but with all 60 of them.
Therefore it is impossible to assign each vertex in $W$ a hyperoval from
$Y$, a contradiction.

The consideration of the case $\D\cong Q_4(3)$ is complete.
Hence the proof of the Theorem \ref{EGQ33} is complete.

\chapter{Geometric characterization of graphs from the Suzuki chain}
\label{chap:suz}
\newcommand{\HH}{{\bf H}}
\newcommand{\Case}{\paragraph}
\paragraph{Abstract.}
\footnote{This research was completed when this author held a position at the 
Institute for System Analysis, Moscow.}
Let $\Sigma_0,\dots,\Sigma_6$ be the graphs from Suzuki chain
\cite{Suz}. We classify connected locally $\Sigma_i$ graphs
$\Theta_{i+1}$ for $i=3,4,5$. If $i=3,4$ then $\Theta_{i+1}$ is
isomorphic to $\Sigma_{i+1}$, whereas $\Theta_6$ is isomorphic
either to $\Sigma_6$ or to its 3-fold antipodal cover
$3\Sigma_6$.
\section{Introduction and results}
\label{suz:sect1}
Recently  there  is  an extensive study of the geometries of the
classical groups as point-line  systems with fixed local structure.
This approach has led to many  beautiful  results.  We  refer  the
reader  to  a  recent paper of Cohen and Shult~\cite{CoSh} for a brief
survey about that.

This paper is an attempt to look 
at the Suzuki chain sporadic simple groups from the aforementioned 
point of view. 

It  is  necessary  to mention a lot of activity in this direction
under some additional assumptions. The first idea  is  to  assume
some  group  action  on  the  geometry  (usually, it is a diagram
geometry). Such an approach lead to study of  some  presentations
for  sporadic  groups, and sometimes to proving the faithfulness of
these presentations. 
For instance, see Ivanov~\cite{Ivnv:J4},  Shpectorov~\cite{ssh:Co2},
Soicher~\cite{Soi:Suz}. 
Another direction is via assuming some  global  property  of  the
geometry,  e.g.  assuming  that  the  related  graph  is strongly
regular    or    distance    regular,    see    the    book    by
Brouwer~et~al.~\cite{BCN}.  This  book  contains many results and many
references    to    results    of    this    kind.    See    also
Cameron~et~al.~\cite{CHP}. It is worth to mention that here we  do  not 
assume either any group action or global property. 

Note that graphs with constant neighbourhood often arise as collinearity
graphs of diagram geometries, see Buekenhout~\cite{Bue:geo:spo}.
Neumaier~\cite{Neu} observed that the graphs from Suzuki chain produce
nice flag-transitive diagram geometries for the related groups.

Throughout the paper we consider undirected graphs without loops
and multiple edges.
Given a graph $\Gamma$, let us denote  the  set  of  vertices  by
$\V=\V\Gamma$,  the  set  of edges by $\E=\E\Gamma$. Let $X\subseteq 
\V\Gamma$.
We denote by $\<X\>=\<X\>_{\Gamma}$ the subgraph $\Xi$
of $\Gamma$ {\em induced} by $X\/$
(i.e. $\V\Xi=X,\/\E\Xi=\{ (u,v)\in \E\Gamma|u,v\in X\}$).
Given two graphs $\Gamma$ and $\Delta$, the graph $\Gamma\cup\Delta$
(resp. the graph $\Gamma\cap\Delta$) is
the graph with the vertex set $\V\Gamma\cup \V\Delta$ (resp.
$\V\Gamma\cap \V\Delta$) and the edge set $\E\Gamma\cup \E\Delta$
(resp. $\E\Gamma\cap \E\Delta$). Given $v\in \V\Gamma$, we denote
$\Gamma_i(v)=\<\{x\in \V\Gamma|x$ at distance $i$ from $v\}\>$,
and $\Gamma_1(v)=\Gamma(v)$.
Furthermore, $\Gamma(X)=\bigcap_{x\in X} \Gamma(x)$.
To simplify the notation we use
$\Gamma(v_1,\dots,v_k)$ instead of $\Gamma(\{v_1,\dots,v_k\})$ and
$u\in\Gamma_i(\dots)$ instead of $u\in \V\Gamma_i(\dots)$.
As usual, $v=v(\Gamma)=|\V\Gamma|$,  $k=k(\Gamma)=v(\Gamma(x))$,
$\lambda=\lambda(x,y)=\lambda(\Gamma)=v(\Gamma(x,y))$, where
$x\in\V\Gamma$, $y\in\Gamma(x)$. Let $y\in\Gamma_2(x)$. We denote
$\mu=\mu(x,y)=\mu(\Gamma)=v(\Gamma(x,y))$. Of course, we use
$k$, $\lambda$, $\mu$ if it makes sense, i.e. if those numbers
are independent on the
particular choice of the corresponding vertices.
If $\Delta$ is a (proper) subgraph of
$\Gamma$ we denote this fact as $\Delta\subseteq\Gamma$ (resp.
$\Delta\subset\Gamma$).

We denote  the complete $n$-vertex graph by $K_n$,  the complete
multipartite graph with the $m$ parts of equal size $n$
by $K_{n\times m}$, the circuit of length $n$ by
$C_n$, the empty graph by $\emptyset$.

$\Aut(\Gamma)$ denotes the automorphism group of $\Gamma$.
Our group-theoretic notation is as in \cite{Atl}.

A {\em strongly regular} graph $\Gamma$ is a connected regular
graph of valency $k$, such that
$\lambda=\lambda(\Gamma)=\lambda(u,v)$ and
$\mu=\mu(\Gamma)=\mu(u,w)>0$ are independent on the choice of
$u\in \V\Gamma$,
$v\in\Gamma(u)$ and $w\in V\Gamma\setminus(\{u\}\cup\Gamma(u))$.
Note that $\mu>0$ implies the diameter of $\Gamma$ is less than
or equal to two.

Let $\Gamma$, $\Delta$ be two graphs. We say that $\Gamma$ is
{\em locally} $\Delta$ if $\Gamma(v)\cong\Delta$ for any $v\in \V(\Gamma)$.

Let $\Gamma,\ \overline\Gamma$ be two graphs. We say that $\Gamma$ is a
{\em cover} of $\overline\Gamma$ if there exists a mapping $\varphi$
from $\V\Gamma$ to $\V\overline\Gamma$ which maps edges to edges and for
any $v \in \V\Gamma$ the restriction of $\varphi$ to $\Gamma(v)$ is an
isomorphism onto $\overline\Gamma(\varphi(v))$. Note that, since
the latter assumption in the definition of cover is usually
omitted, our definition of cover is a bit nonobvious.

Suppose we have a chain of graphs $\Sigma_1,\dots,\Sigma_n$, such that
$\Sigma_i$ is locally $\Sigma_{i-1},\  i=2,\dots,n$. Then for any complete
$k$-vertex subgraph $\Upsilon$ of $\Sigma_m$, $1<k<m\leq n$
we have $\Sigma_m(\V\Upsilon)\cong\Sigma_{m-k}$.

{\em Suzuki chain} \cite{Suz} consists of the following graphs 
$\Sigma_1,\dots,\Sigma_6$.

$\Sigma_1$ is the empty four-vertex graph.

$\Sigma_2$ is the incidence graph of the (unique) $2-(7,4,2)$ design.

$\Sigma_3,\dots,\Sigma_6$ are strongly regular with the parameters
$(v,k,\lambda,\mu)$ equal to
$$(36, 14, 4, 6),\ (100, 36, 14, 12),\ (416, 100, 36, 20),\ (1782, 
416, 100, 96)$$ respectively.
The automorphism groups of $\Sigma_i$ are $S_4$, ${\it PGL}_2 (7)$, $G_2 (2)$,
$J_2\!:\!2$, $\Aut(G_2(4))$, $\Aut(\Suz)$ respectively.

Note that $\Sigma_i,\  i=2,\dots ,6$ is locally $\Sigma_{i-1}$.
Besides $\Sigma_6$, there exists another connected locally $\Sigma_5$ graph
$3\Sigma_6$ -- the 3-fold antipodal cover of $\Sigma_6$.
$3\Sigma_6$ is
a distance transitive graph (see \cite{BCN}) with the intersection
array $$\{ 416,315,64,1;1,32,315,416\}$$. It has been
constructed by Soicher~\cite{Soi:Suz} (see also \cite{Soi:drg}).

\begin{thm}
 Let $\Theta_{i+1}$ be a connected locally $\Sigma_i$
graph, $i=3,4,5$. Then $\Theta_4\cong\Sigma_4$, $\Theta_5\cong\Sigma_5$,
$\Theta_6\cong\Sigma_6$ or $3\Sigma_6$.
\label{suz:MTH}
\end{thm}

Note that under the additional assumption that $\Aut(\Theta_{i+1})$ is
transitive on the ordered $k$-cliques of $\Theta_{i+1}$ for any $k\leq i$,
this statement is proved in \cite{Soi:Suz} using a coset enumeration. 

\medskip
{\sl Remark.} Results from \cite{Soi:Suz} and our Theorem~\ref{suz:MTH}
immediately lead to
a computer-free proof of the faithfulness of the well-known presentations for
groups from Suzuki chain, given in \cite{Atl}. Namely if $G_n$ is given by
the presentation
$$\bigl\langle_{g_n}
\!\!\!\!\!\circ \!{\phantom{00000}\over \phantom{0}}
{\vphantom{\bigl\langle}}_{g_{n-1}}\!\!\!\!\!\!\!\!\!\!\!\!\!\circ \ \cdots\
{\vphantom{\bigl\langle}}_{g_3}
\!\!\!\!\!\!\!\circ \!\!{\phantom{00000}\over \phantom{0}}
{\vphantom{\bigl\langle}}_{g_2}
\!\!\!\!\!\!\!\circ \!\!{{\phantom{00} 8 \phantom{00}}\over \phantom{0}}
{\vphantom{\bigl\langle}}_{g_1}
\!\!\!\!\!\!\!\circ \!\!{\phantom{00000}\over \phantom{0}}
{\vphantom{\bigl\langle}}_{g_0}
\!\!\!\!\!\!\!\circ \!\!{\phantom{00000}\over \phantom{0}}
{\vphantom{\bigl\langle}}_a
\!\!\!\!\!\circ\ |\
(g_1 g_2)^4=a,\ (g_0 g_1 g_2 g_3)^8=1 \bigr\rangle, $$
then $G_4 \cong\Aut(J_2),\ G_5 \cong\Aut(G_2(4)),\ G_6 \cong 3\Aut(\Suz)$.

\section{Proof of Theorem~\protect\ref{suz:MTH}}
\label{suz:sect2}
Trying to simplify the reading of Sect.~\ref{suz:sect2},  
we  give  a  sketch  of proof here. 

Let $\Gamma=\Theta_i$. Pick $u\in \V\Gamma$ and look at 
$\Gamma(u,v)$ for each $v\in\Gamma_2(u)$. 
First, we establish the characterization of $\Gamma(u,v)$ in terms 
of their local structure (Lemma~\ref{suz:L1}~(1)). Then we classify  the 
connected components of the possible candidates  to $\Gamma(u,v)$ 
(Proposition~\ref{suz:P3}). Next, we list all  the  possible groupings  of 
$\Gamma(u,v)$ from the connected components (Lemma~\ref{suz:L5}). 
This list 
and  Lemma~\ref{suz:L1}~(2) allows us to determine the edges between 
$\Gamma(u)$ and $\Gamma_2(u)$. Using Proposition~\ref{suz:P2} we express the 
adjacency in $\Gamma_2(u)$ in terms of intersections of 
$\Gamma(u,v)$-subgraphs (Sect.~\ref{suz:sect23}
($i=4,5$), Sect.~\ref{suz:sect24} ($i=6$), also 
Lemmas~\ref{suz:L6} and~\ref{suz:L7}
are used). In cases $i=4,5$ we are done by Lemma~\ref{suz:L4}, 
which states that $\Gamma(u,v)$ are "hyperplanes" of $\Gamma(v)$. 
So holds in a subcase (Case $c(\Omega)=3$) of the case $i=6$, as well. 
In the remaining subcase (Case $c(\Omega)=1$, Case $c(\Omega)=2$ 
immediately leads to a contradiction) of the case $i=6$ we show 
that $\Gamma$ is a (unique) covering of $\Sigma_6$, hence 
$\Gamma\cong 3\Sigma_6$. 

\subsection{Suzuki's construction }
\label{suz:sect21}

We include the necessary information about the construction of 
$\Sigma_i,\/$ given by Suzuki~\cite{Suz}, $i=3,\dots,6$ (see also \cite{Atl}).

Let $\Delta=\Sigma_{i-1}$, $\Gamma=\Sigma_i$ and $H=\Aut(\Delta)$ 
for $i=3,\dots,6$. 
In fact, $\Gamma$ can be defined in terms of $\Delta$ and $H$. 

Let $\infty$ be an extra symbol. Let $S$ be the 
conjugacy class of $2$-subgroups $Z=\<z\>$ of $H$, where $z$ is a 
2A-involution of $H$, $i=3,4,5$. In the remaining case $i=6$ $S$ is the 
conjugacy class of $2^2$-subgroups $Z=\<x,y\>$ of $H$,
such that $[H:N_H(Z)]=1365$ and $x,\/y$ are 2A-involutions of $H$.
Then $\V\Gamma=\{\infty\} \cup \V\Delta \cup S.\/$

$\E\Gamma$ is exactly as follows. The vertex $\infty$ is adjacent to each 
vertex 
from $\V\Delta$, two vertices from $\V\Delta$ are adjacent if they are 
adjacent 
as vertices of $\Delta$, a vertex $x\in S\/$ is adjacent to a vertex
$v\in \V\Delta$ if a nontrivial element of $x$, considering as a subgroup
of $H$, fixes $v$, two vertices $x,y\in S\/$ are adjacent if $x,\/y$,
considering as subgroups of $H$, do not commute, but there exists $z\in S\/$
commuting with both of them.

\subsection{Neighbourhood of two vertices at distance two}
\label{suz:sect22}
We start with a simple general fact. Let $\Delta$ be a connected
graph satisfying the following property \medskip

\begin{itemize}
\item[$(*)$] For any $u\in \V\Delta$ and
$v\in \V\Delta\setminus (\V\Delta(u)\cup\{ u \})$\\
1)~the subgraph $\Delta(u,v)$ is isomorphic to some $M_{\Delta}$,
whose isomorphism type is independent on the particular choice of $u,\/ v 
$.\\ 
2)~for $w\in \V\Delta\setminus (\V\Delta(u)\cup\{ u \}),\/w\not=v,\/$
$\Delta(u,v)\not=\Delta(u,w)$.
\end{itemize}

\begin{lem}
Let $\Gamma$ be locally $\Delta$ graph, where $\Delta$
satisfies $(*)$. For any $u\in \V\Gamma$ and $v\in\Gamma_2(u)$
the following holds 
\begin{itemize}
\item[$(1)$] $\Gamma(u,v)$ is locally $M_{\Delta}$, 
\item[$(2)$] for $w\in\Gamma_2(u),\ w\not=v\/$, we have 
        $\Gamma(u,v)\not=\Gamma(u,w)$.
\end{itemize}
\label{suz:L1}
\end{lem}

\begin{pf}
By 1) of (*), we have $\Gamma(u,v,y)\cong  M_{\Delta}$  for 
any $y\in \Gamma(u,v)$. The first claim is proved.
By way of contradiction, let $w\in\Gamma_2(u),\/w\not=v\/$,\break
$\Gamma(u,v)=\Gamma(u,w)$. Pick $y\in\Gamma(u,v)$.
Looking at the common neighbourhood of $u,\/v$ and $u,\/w$
in $\Gamma(y)\cong\Delta$, we obtain a contradiction. \end{pf}

Note that $\Delta=\Sigma_i$ ($i=3,4,5$) satisfies (*) above. 
Indeed, the stabilizer of $u\in \V\Delta$ in $\Aut(\Delta)$ acts 
transitively and primitively on 
$\Delta_2(u)$. The transitivity  implies  1)  of  (*),  while  the 
primitivity 
2) of (*). Thus  Lemma~\ref{suz:L1}  holds  for  locally  $\Sigma_i$  graphs 
$\Gamma$   for   any   $i=3,4,5$.   In   particular,    for    any 
$u\in\V\Gamma$ and $v\in\Gamma_2(u)$  the  subgraph  $\Gamma(u,v)$ 
is a locally $M_{\Sigma_i}$ subgraph  of  $\Gamma(u)\cong\Sigma_i$ 
(i=3,4,5). We turn to the classification of the 
locally $M_{\Sigma_i}$ subgraphs of $\Sigma_i$ ($i=3,4,5$). 

First, we need a technical statement. Note that in fact we will prove
slightly more than we will really need. For a graph $\Gamma$ and $g\in 
\Aut(\Gamma)$ we denote the 
subgraph induced by the vertices fixed by $g$ as $Fix(g)=Fix_\Gamma(g)$.

\begin{pr}\label{suz:P2}
Let $g$ be a 2A-involution in $\Aut(\Sigma_i)$, $i=3,4,5$. Then
the number of 2A-involutions $h$ corresponding to particular isomorphism
type of the subgraph $Fix(g)\cap Fix(h)$ is presented in
Table~\ref{suz:tab1}. 
\end{pr}
\begin{table}
\begin{center}
\begin{tabular}{|c|c|c|c|c|} \hline
$i$ &\multicolumn{4}{|c|}{$Fix(g)\cap Fix(h)$}\\ \hline
3 &6 &24 &32 &--\\ \cline{2-5}
&$K_{4\times 2}$ &$K_{2\times 2}$ &$K_3$ &--\\ \hline
4 &10 &80 &160 &64 \\ \cline{2-5}
&$K_{4\times 3}$ &$K_{2\times 3}$ &$K_4$ &$\emptyset$ \\ \hline
5 &20 &320 &1024 &--\\ \cline{2-5}
&$K_{4\times 4}$ &$K_{2\times 4}$ &$K_5$ &--\\ \hline
\end{tabular}
\end{center}
\caption{ Number of 2A-involutions $h\in Aut(\Sigma_i)$
vs. the isomorphism type of $Fix(g)\cap Fix(h)$.  }
\label{suz:tab1}
\end{table}
\begin{pf}
We follow the notation of Sect.~\ref{suz:sect21}
Denote $\Gamma=\Sigma_i$, $\Omega=Fix(g)$, $G=\Aut(\Gamma)$.
Pick  $\infty\in  \V\Omega$,  denote  $\Delta=\Gamma(\infty)$.  It 
follows from the Suzuki's construction given in Sect.~\ref{suz:sect21} 
that $\Omega$ is locally $M_\Gamma$. It should be mentioned that
$\Omega_2(\infty)=\{x \in g^{\Aut(\Delta)}|(xg)^2=1\}$. Let us
denote by $\Xi_h$ the subgraph $\Omega\cap Fix(h)$ for $h\in g^G$. 

\Case{$i=3$} There exist two distinct subgraphs $\Xi$ of $\Omega$ 
isomorphic to $K_{4\times 2}$ and containing $\infty$. It is easy 
to check that 
$\Xi$ determines exactly two $h\in g^G$ such that
$\Xi_h=\Xi$. Hence there exist exactly
$2\cdot 2v(\Omega)/v(K_{4\times 2})=6$ involutions $h\in g^G$ such that
$\Xi_h\cong K_{4\times 2}$.

It is well-known that the action of $G$ on $g^G$ coincides with the
action on the set of points of a generalized  hexagon \HH\ of order 
(2,2) and that 
the corresponding subdegrees are $1,6,24,32$ according to the distance
in the collinearity graph of \HH. Hence if $h'\in g^G$ corresponds to
the subdegree 24 then there exists
$h\in g^G$ such that $\Xi_h\cong Fix(h)\cap Fix(h')\cong K_{4\times 2}$.
Hence $v(\Xi_{h'})\geq 4$. On the other hand it is easy to find 
elements 
$h',h''\in g^G$ with $\Xi_{h'}\cong K_{2\times 2}$ and $\Xi_{h''}\cong K_3$.
The proof for $i=3$ is complete.

\Case{$i=4$} By induction, there exist 6 elements $h_j\in g^G$ such that
$\Xi_{h_j}(\infty)\cong K_{4\times 2}$. Pick $y=h_1$ and $u\in
\V\Xi_y(\infty)$. Since $v(\Xi_y(u))\geq 1+k(K_{4\times 2})$, we have
$\Xi_y(u)\cong K_{4\times 2}$. Hence a connected component $\Upsilon$
of $\Xi_y$ is locally $K_{4\times 2}$. Thus $\Upsilon\cong
K_{4\times3}$ (cf. \cite[Proposition~1.1.5]{BCN}).
It can be easily shown that for any $w \in \V\Omega$ the subgraph
$\Omega(w)\cap\Upsilon$ is nonempty. Hence $\Xi_y=\Upsilon$. Thus there exist
exactly $6v(\Omega)/v(K_{4\times 3})=10$ involutions $h\in g^G$ such that
$\Xi_h\cong K_{4\times 3}$.

There exist 24 elements $h'_j\in g^G$ such that
$\Xi_{h'_j}(\infty)\cong K_{2\times 2}$. Pick $y=h'_1$ and $u\in 
\V\Xi_y(\infty)$.
Since $\Xi_y(u)$ contains a 2-path, we have $\Xi_y(u)\cong K_{2\times 2}$ or
$K_{4\times 2}$. The latter case is impossible. Indeed, as we have already
proved, if the latter case holds, we have $\Xi_y\cong K_{3\times 
4}$. Thus $\Xi_y(\infty)\cong K_{2\times 4}$, contradiction.
Hence a connected component $\Upsilon$ of $\Xi_y$
is locally $K_{2\times 2}$. Thus $\Upsilon\cong K_{2\times 3}$.
Therefore a connected component of $\Xi_y$ distinct from 
$\Upsilon$ (if such a component exists at all) is a
singleton $s\in \V\Xi_y$. Looking at the neighbourhood of $s$, we obtain a
contradiction, since $\Xi_y(s)\cong\emptyset$. Hence $\Xi_y=\Upsilon$.
Thus there exist exactly $24v(\Omega)/v(K_{2\times 3})=80$ 
elements $h'\in g^G$ such that $\Xi_{h'}\cong K_{2\times 3}$. 

There exist 32 elements $h''_j\in g^G$ such that
$\Xi_{h''_j}(\infty)\cong K_3$. Pick $y=h''_1$ and $u\in \V\Xi_y(\infty)$.
Similarly to the previous case, we obtain that a connected component
$\Upsilon$ of $\Xi_y$ is locally $K_3$. Thus $\Upsilon\cong K_4$.
Again, looking at $\Xi_y$ as at a subgraph of $\Omega$, we obtain
$\Xi_y=\Upsilon$.
Thus there exist exactly $32v(\Omega)/v(K_4)=160$ elements $h''\in 
g^G$ such that $\Xi_{h''}\cong K_4$.

Clearly for the other $h\in g^G$ we have $\Xi_h\cong\emptyset$.
The proof in the case $i=4$ is complete.

\Case{$i=5$} The argument runs parallel to the argument for the case $i=4$.
\end{pf}

\begin{pr}
Any connected locally $M_{\Sigma_i}$ subgraph
of $\Sigma_i$ coincides with the subgraph $Fix(g)$ for a
2A-involution $g\in \Aut(\Sigma_i)$ $(i=3,4,5)$.
\label{suz:P3}
\end{pr}

\begin{pf}
Denote $\Gamma=\Sigma_i$. Let $\Omega$ be a connected locally 
$M_\Gamma$ subgraph
of $\Gamma$, $i=3,4,5$. We will frequently use the construction of $\Gamma$
given in Sect.~\ref{suz:sect21}, as well as the corresponding notation.
Without loss of generality $\infty\in \V\Omega$.

\Case{$i=3$} Recall that $\Delta$ is the incidence graph of the
unique $2-(7,4,2)$ design. Since the latter design is complementary to
$\Pi=PG(2,2)\/$, the equivalent description of $\Delta$ (and $\Gamma$)
can be given in terms of $\Pi$. Namely, the vertices of $\Delta$ are
the points and the lines of $\Pi$, a point $p$ is adjacent to a line
$l$ if $p$ does not belong to $l$.

For each 2A-involution $v$ of $H$ there exists a unique flag 
(= incident point-line pair) $(p,l)$ such that $v$ leaves all points on
$l$ and all lines through $p$ fixed. Thus it is easy to check that
the subgraph $M_\Gamma=\Gamma(\infty,v)$ is isomorphic to
$K_1\cup K_1\cup C_4$.

It is easy to check that {\em any} (generated) subgraph of
$\Delta$ of the latter shape is $\Gamma(\infty,v)$ for some
$v\in \Gamma_2(\infty)$. Indeed, $\Delta$ has the remarkable property
that any 2-path in $\Delta$ lies in a unique $C_4$ (i.e. $\Delta$ 
is a {\em rectagraph} (see \cite{BCN})). Moreover,
$\Aut(\Delta)$   acts  transitively on the $C_4$-subgraphs of 
$\Delta$. Finally, it is
straightforward to check that there exist exactly two vertices of $\Delta$
which are nonadjacent to any vertex of a subgraph $\Xi$ of 
$\Delta$ isomorphic to $C_4$. 

Thus, without loss of generality, we may assume that
$\Gamma(\infty,v)\subset \Omega$ for some $v\in \Gamma_2(\infty)$. Let
$\Xi$ be the $C_4$-subgraph of $\Omega(\infty)$, $u\in \V\Xi$. 
Since $\Gamma(u)\cong\Delta$ is a rectagraph, a quick 
look at $\Xi(u)$ gives us $v\in \V\Omega$. Next, $v$ is adjacent to each
vertex in $\Omega(\infty)$. Hence any $w\in \Gamma_2(\infty)\cap\Gamma(v)$
does not belong to $\V\Omega$.

Let $x\in \Omega(\infty)\setminus \V\Xi$.
The unique $C_4$-subgraph of $\Omega(x)$ lies in
$\Gamma_2(\infty)$. On the other hand $\Gamma(x,v)\subset\Gamma_2(\infty)$.
Hence for $Q=\Gamma_2(\infty)\cap\Gamma_2(v)\cap\Gamma(x)$ we have
$v(Q)=4$. Thus $Q\subset \Omega$.

Consider $v$ as a 2A-involution of $H$. It is easy to
see, using \cite{Atl}, that $v$ is a 2A-involution of $\Aut(\Gamma)$. Therefore
its fixed vertices are
$\{\infty\}\cup\{v\}\cup \V\Gamma(\infty,v)\cup \V Q\subseteq \V\Omega$. Since
the subgraph $Fix(v)$ induced by the vertices fixed by $v$ is connected and 
locally
$M_\Gamma$, we have equality in the latter inclusion. This completes
the proof in the case $i=3$.

\Case{$i=4$} By induction, any subgraph of $\Sigma_3$ isomorphic to
$\Gamma(\infty,v)$ for some $v\in \Gamma_2(\infty)$ is in fact $Fix(g)$ for
some 2A-involution $g\in H\/$. Hence, without loss of generality, we may 
assume
that $\Gamma(\infty,v)\subset \Omega\/$.

Let $u\in \Omega(\infty),\ v\in \Omega_2(\infty)\cap\Omega(u)\/$.
Remembering the structure of $\Omega(u)$, we have
$v(\Omega(\infty,u,v))\geq 4$. Hence, by Proposition~\ref{suz:P2}
we have
$v(\Omega(\infty,v))\geq 8$. So $v(\Omega_2(\infty))\leq 7$. Let $Q$ be
the set of {\em all} possible candidates to $\V\Omega_2(\infty)$. Since
$\lambda(\Omega)=6$, there are precisely $60$ edges of $\E\Omega$ going
from $\V\Omega(\infty)$ to $\V\Omega_2(\infty)$. On the other hand precisely
$60$ edges are going from $Q$ to $\V\Omega(\infty)$. Hence 
$Q=\V\Omega_2(\infty)$.
In particular, a 2A-involution $w\in H$ belongs to $Q$. Next, $w$ is
a 2A-involution of $\Aut(\Gamma)$. Therefore $\Omega=Fix(w)$.

\Case{$i=5$} The argument runs parallel to the argument in the
previous case. 
\end{pf}

We have already known the connected components of locally $M_\Gamma$
subgraphs $\Omega$. Let us study the relationship between $\Omega$
and the vertices outside $\V\Omega$.

\begin{lem}
Let $\Omega$ be connected locally $M_\Gamma$,
$v\in \V\Gamma\setminus \V\Omega$, $\Gamma(v)\cap \Omega\not=\emptyset$,
where $\Gamma=\Sigma_i$ $(i=3,4,5)$.
Then $\Gamma(v)\cap\Omega\cong K_{2\times (i-1)}$.
If $i=3,4$ then for any $v\in \V\Gamma\setminus \V\Omega$ the subgraph
$\Gamma(v)\cap \Omega\not=\emptyset$.
In the remaining case $i=5$ for any $v\in \V\Gamma\setminus \V\Omega$
either the subgraph
$\Gamma(v)\cap \Omega$ is nonempty or $v$ lies in a uniquely determined
connected locally $M_\Gamma$ subgraph.
\label{suz:L4}
\end{lem}

\begin{pf}
If $i=3,4$ then the centralizer of a 2A-involution in $\Aut(\Gamma)$
has two orbits on $\V\Gamma$ (look at the scalar product of permutation
characters, given in \cite{Atl}). Hence the isomorphism type of
$\Xi=\Xi_v=\Gamma(v)\cap \Omega$ is independent on the choice of
$v\in \V\Gamma\setminus \V\Omega$. Next,
$$v(\Xi)(v(\Gamma)-v(\Omega))=v(\Omega)(k(\Gamma)-k(\Omega)).$$
Thus $v(\Xi)=4$ or $6$ according to $i=3$ or $4$.

If $i=3$ it is easy to
check that $k(\Xi)=2$. Indeed, pick a
vertex $u\in \V\Xi$. Look at $\Delta=\Gamma(u)$. Recall that 
$\V\Omega(u)$ is the set of fixed points of a 2A-involution $g\in 
\Aut(\Delta)$.
$C_{\Aut(\Delta)}(g)$ has one orbit on $\V\Delta\setminus \V\Omega(u)$. 
Performing
the same calculation as above, we have $k(\Xi)=v(\Omega(u)\cap\Delta(v))=2$.
The proof in the case $i=3$ is complete.

By induction if $i=4$ then $\Xi$ is locally $C_4=K_{2\times 2}$. Such
a graph $\Xi$ is unique and isomorphic to $K_{2\times 3}$.

Consider $\Sigma_5$ as the subgraph $\Delta=\Gamma(\infty)$
of $\Gamma=\Sigma_6$, as in Sect.~\ref{suz:sect21}. For a vertex $u\in
\Gamma_2(\infty)$ the group $H={\rm N}_{\Aut(\Delta)}(u)$ (recall that
$u$ is a $2^2$-subgroup of $\Aut(\Delta)$) has two orbits on $\Delta$.
By Sect.~\ref{suz:sect21} and Proposition~\ref{suz:P3}, $\Gamma(\infty,u)$
is the disjoint union of three copies of $\Omega$.
This enables us to calculate, for $v\in \V\Gamma(\infty)\setminus
\V\Gamma(u)$, the number $n$ of edges coming from $v$ to
$\Gamma(\infty,u)$. Namely $n=24$.
There exists $g\in H(v)$ of order 3 permuting cyclically the connected
components of $\Gamma(\infty,u)$.
Hence $v(\Xi)=n/3=8$. On the other hand $\Xi$ is locally $K_{2\times3}$.
Such a graph is isomorphic to $K_{2\times 4}$.

Finally, the last statement of the lemma follows directly from the
aforementioned consideration of $\Sigma_5$ as a subgraph of $\Sigma_6$. \end{pf}

Lemma~\ref{suz:L4} immediately implies

\begin{lem}
Let $\Omega$ be locally $M_{\Sigma_i}$ $(i=3,4,5)$.
If $i=3,4$ then $\Omega$ is connected. If $i=5$ then there exists a unique
equivalence relation $\varphi$ with the class size three on the set of
connected locally $M_{\Sigma_i}$ subgraphs such that all connected
components of $\Omega$ lie in one class of $\varphi$. \qed
\label{suz:L5}
\end{lem}

\subsection{The final part of proof for i=3 and 4}
\label{suz:sect23}
Let $\Delta=\Sigma_i$, $\Gamma$ be a connected locally $\Delta$ graph.
Let us identify $\Delta$ with $\Gamma(\infty)$ for $\infty\in \V\Gamma$.
Since $\Delta$ satisfies (*), by the classification of the locally $M_\Delta$
subgraphs of $\Delta$ given in Sect.~\ref{suz:sect22}, we have that
$\mu(\Gamma)$ equals to the number of the fixed points of some
2A-involution $g\in H=\Aut(\Delta)$. Since 
$$v(\Gamma_2(\infty))=k(\Gamma)\cdot (k(\Gamma)-k(\Delta)-1)/ \mu(\Gamma),$$
$v(\Gamma_2(\infty))=|g^H|$. Since $\Delta$ satisfies (*), by
Lemma~\ref{suz:L1} for a locally $M_\Delta$ subgraph
$\Omega$ of $\Delta$ there exists a unique $v=v_\Omega \in \Gamma_2(\infty)$ 
such that
$\Omega=\Gamma(\infty,v)$ and for two distinct locally $M_\Delta$ subgraphs
$\Omega,\;\Omega'$ we have $v_\Omega\not=v_{\Omega'}$. Thus we have
proved that we do not have any choice determining the edges from
$\Gamma(\infty)$ to $\Gamma_2(\infty)$.

Let $v\in \Gamma_2(\infty)$, $u\in \Gamma(v)\setminus\Gamma(\infty)$. It
follows from Lemma~\ref{suz:L4} , applied to $\Gamma(v)\cong\Delta$, that 
$\Gamma(\infty,v,u)\cong K_{2\times (i-1)}$
(in particular, $\Gamma$ has diameter two).
Next, $v(\Gamma_2(\infty)\cap\Gamma(v))=24$ if $i=3$ ($80$ if $i=4$). On
the other hand by Proposition~\ref{suz:P2}
if $i=3$ there exist exactly $24$ ($80$ if $i=4$) locally
$M_\Delta$ subgraphs of $\Delta$ intersecting $\Gamma(\infty,v)$ in a
subgraph of the shape $K_{2\times (i-1)}$. Thus the adjacency in
$\Gamma_2(\infty)$ is unique determined, as well.

The proof of Theorem~\ref{suz:MTH} is complete for $i=3,4$.

\subsection{The final part of proof for i=5}
\label{suz:sect24}
Let $\Delta=\Sigma_5$, $\Gamma$ be a connected locally $\Delta$ graph.
First we prove further lemmas about locally $M_\Delta$ subgraphs of $\Delta$.
Denote $G=\Aut(\Delta)$. Let $\Omega$ be a connected locally $M_\Delta$
subgraph of $\Delta$. Define the graph $\Lambda$ with
$\V\Lambda=\Omega^G,\ \E\Lambda=\{(\Omega',\Omega'')|\Omega',\Omega''\in
\Omega^G,\ \Omega'\cap\Omega''\cong K_{2\times 4}\}$.

\begin{lem}
The graph $\Lambda$ is a connected graph of valency
$k(\Lambda)=320$, on which $G$ acts vertex- and edge-transitively.
For each $v\in\V\Lambda$ and $u\in\Lambda_2(v)$ we have $k(\Lambda(v,u))\leq 
32$.
\label{suz:L6}
\end{lem}

\begin{pf}
The vertex-transitivity is clear. $k(\Lambda)=320$ follows immediately
from Proposition~\ref{suz:P2}. 
Let $v\in\V\Lambda$. Consider the action of the stabilizer 
$F$ of $\Omega$
in $G(v)$ on the set of 80 locally $M_{\Delta(\infty)}$ subgraphs
$\Xi_j(v)$ of $\Delta(v)$ such that $\Omega(v)\cap \Xi_j(v)\cong
K_{2\times 3}$ (cf. the case $i=4$ of Proposition~\ref{suz:P2}). 
This action is transitive, since
$\Xi_j(v)$ ($j=1,\dots,80$) constitute a suborbit of $\Aut(J_2)$ in its primitive
action of degree $315$. By Proposition~\ref{suz:P3}
each $\Xi_j(v)$ lifts to a unique
locally $M_\Delta$ subgraph $\Xi_j$ of $\Delta$. By Proposition~\ref{suz:P2}
we have $\Xi_j \cap \Omega \cong K_{2 \times 4}$. Hence $F$ acts transitively
on the $\Xi_j$ ($j=1,\dots,80$). Since the stabilizer of $\Omega$ in $G$
acts transitively on $\V\Omega$, the edge-transitivity is proved.

The unique imprimitivity system of $G$ on $\V\Lambda$ is given by 
the equivalence relation $\varphi$ defined in Lemma~\ref{suz:L5}
(cf. \cite{Atl}). Hence $\Lambda$ is connected.

The last statement follows from the fact that $\Lambda$ appears as
the subgraph $\Gamma_2(\infty)$ of the graph $\Gamma\cong 3\Sigma_6$ 
and that $\mu(\Gamma)=32$.
\end{pf}

Let us denote by $\overline\Omega$ the $\varphi$-equivalence class of
$\Omega$. Define the graph $\overline\Lambda$ by 
$\V\overline\Lambda=\overline\Omega^G$ and 
$\E\overline\Lambda=\{(\overline\Omega',\overline\Omega'')|
\overline\Omega',\overline\Omega''\in \overline\Omega^G,\
\Omega'\cap\Omega''\cong K_{2 \times 4}\}$.

\begin{lem}
\label{suz:L7}
For each $(\overline\Omega',\overline\Omega'')$,
we have $\overline\Omega'\cap\overline\Omega''\cong
K_{2 \times 4}\cup K_{2 \times 4}\cup K_{2 \times 4}$.
The graph $\overline\Lambda$ is connected graph of valency 
$k(\overline\Lambda)=320$.
\end{lem}

\begin{pf}
Let $S$ be the set of $\Xi\in\Omega^G$ satisfying 
$\Xi\cap\Upsilon\cong\ K_{2 \times 4}$ for some $\Upsilon$ 
$\varphi$-equivalent to $\Omega$. Since $\Omega$ does not contain a 
subgraph $K_{2 \times 4}\cup K_{2 \times 4}$, there is a unique 
such subgraph $\Upsilon$ for each $\Xi$ of $S$.
Hence the stabilizer $F$ of $\Omega$ in $G$ acts on $S$ transitively.
The subdegrees of $G$ on $\overline\Omega^G$ are $1,20,320,1024$.
Since $|S|=960,\ S=\bigcup_{\Xi\in S}\overline\Xi$. Hence
we have proved the first two claims of the lemma. The last one follows
from the primitivity of $G$ on $\overline\Omega^G$. \end{pf}

Let us identify $\Delta$ with $\Gamma(\infty)$ for $\infty\in \V\Gamma$.
Choose a vertex $v\in \Gamma_2(\infty)$. Denote $\Omega=\Gamma(\infty,v)$.
According to the results of Sect.~\ref{suz:sect22}, 
we know the connected components of
$\Omega$ up to conjugacy. By Lemma~\ref{suz:L5}
we have several possibilities on the number of the
components $c(\Omega)$ of $\Omega$.

\Case{$c(\Omega)=3$} We have $\Omega=\overline\Theta$
for a connected locally $M_\Delta$ subgraph $\Theta$ of $\Delta$. It
follows from Lemmas~\ref{suz:L4} and~\ref{suz:L5} that the subgraph
$\Xi(v)=\<\V\Gamma(v)\setminus\V\Omega\>$ lies in $\Gamma_2(\infty)$.
Clearly $v(\Xi(v))=v(\Delta)-v(\Omega)=320$. Let
$u \in \Xi(v)$. Then by Lemma~\ref{suz:L4} $\Gamma(\infty,v,u)
\cong K_{2 \times 4}\cup K_{2 \times 4}\cup K_{2 \times 4}$,
where the connected components of this subgraph lie in distinct
connected components of $\Omega$. Hence
$\Gamma(\infty,u)=\overline\Theta '$ for a connected locally $M_\Delta$
subgraph $\Theta'$ of $\Delta$. By Lemma~\ref{suz:L7}
a connected component $\Xi$ of
$\Gamma_2(\infty)$ containing $v$ is isomorphic to
$\overline\Lambda$. Counting in two ways the number of edges
from $\Gamma(\infty)$ to $\Gamma_2(\infty)$, we obtain that
$\Gamma_2(\infty)=\Xi$. Thus $\Gamma\cong\Sigma_6$ and the proof
in this case was completed.

\Case{$c(\Omega)=1$} Consider $u \in \Gamma(v)\cap \Gamma_2(\infty)$.
By Lemma~\ref{suz:L4} $\Gamma(\infty,v,u) \cong K_{2 \times 4}$. Therefore
$\Gamma(\infty,u)$ is connected. Look at the connected component $\Upsilon$ of
$\Gamma_2(\infty)$ containing $v$. By Lemmas~\ref{suz:L6} and~\ref{suz:L4}, 
we have $k(\Upsilon)\geq 320$ and $\Upsilon$
contains $\Lambda$ as a (possibly non-generated) subgraph.
Calculating in two ways the number of edges coming from $\Gamma(\infty)$ to
$\Gamma_2(\infty)$, we obtain that $\Gamma_2(\infty)=\Upsilon$.

Let us prove that $\Upsilon\cong\Lambda$, i.e. if $(u,v)\in\E\Upsilon$
then $\Gamma(\infty,u,v)\cong K_{2\times 4}$. Suppose contrary for
some $(u,v)\in\E\Upsilon$. We know that $\Gamma(v,u)\cong\Sigma_4$.
Let us determine to which layers of $\Gamma$ with respect to $\infty$ the
vertices of $\Gamma(v,u)$ do belong.
By Lemma~\ref{suz:L4} $\Gamma(\infty,u,v)=\emptyset$. Next, $u$ lies in a uniquely
determined locally $M_\Delta$ subgraph $\Omega'$  of  $\Gamma(v)$.
Namely, $\Omega'\in\overline\Omega$ (cf. Lemma~\ref{suz:L5}).
Dually, $v$ lies in such a subgraph of $\Gamma(u)$.
Therefore the number of vertices of
$\Gamma(v,u)$ belonging to $\Gamma_3(\infty)$ cannot exceed $k(\Omega')=20$.
Thus, $\Gamma(v,u)$ has at least $80$ vertices at the distance $2$ from
$\infty$. Let $w$ be a vertex of this set. It is easy to check that
$\Gamma(\infty,v,w)\cong\Gamma(\infty,u,w)\cong K_{2\times 4}$.
One can consider $(v,w,u)$ as a 2-path in $\Lambda$. There are at least
$80$ such paths from $v$ to $u$, which contradicts the last part of
Lemma~\ref{suz:L6}. We are done.

Let $w \in \Gamma_3(\infty)\cap\Gamma(v)$. Denote $\Theta=\Gamma(w)$.

\begin{lem}
\label{suz:L8}
The subgraph $\Xi=\Gamma_2(\infty)\cap\Theta$ equals
$\Theta_2(t)$ for some $t\in \V\Theta$.
\end{lem}

\begin{pf}
From the above we know that $w \in \V {\rm A}$ for a connected
locally $M_\Delta$ subgraph ${\rm A}\subset\Gamma_3(\infty)\cap\Gamma(v)$.
Thus $k(\Xi)=v(\Theta(v))-k({\rm A})=80$. By Proposition~\ref{suz:P3}
for any $u\in 
\V\Xi$ the subgraph $\overline\Xi(u)=\Theta(u)\cap\Gamma_3(\infty)$ equals
$\Theta(u,t^u)$ for some $t^u \in \Theta_2(u)$. This faces us to the 
problem of classification of the subgraphs of $\Theta$ satisfying the 
local property above.

Denote $t=t^v$. It suffices to show that  
$\Pi^u=\overline\Xi(u)\cap\Theta(u,t)\setminus\overline\Xi(v)$ 
is nonempty for any $u\in\Xi(v)$. Indeed, since the pointwise 
stabilizer of $v$ and $t$ in $\Aut(\Theta)$ acts transitively on
$\Theta(t)\setminus\overline\Xi(v)$, the sets $\Pi^u$ $(u\in\Xi(v))$
cover it, i.e. $\Theta(t)\cap\V\Xi=\emptyset$. 

Pick $u\in\Xi(v)$. Denote $\Pi=\Pi^u$. By Lemma~\ref{suz:L4}, 
$\Phi=\Theta(v,t,u)= \Theta(v)\cap\overline\Xi(u)$ 
is isomorphic to $K_{2\times 3}$.
Pick $(x,y)\in\E\Phi$. Consider $\Psi=\Theta(u,x,y)\cong\Sigma_2$.
The subgraph $\Phi(x,y)$ of $\Psi$ consists of two nonadjacent vertices 
$a$ and $b$. We are interesting in ways of completion of $\{a,b\}$ to
the subgraph of fixed points of a 2A-involution of $\Aut(\Psi)$ (one
such way gives $\Theta(u,t,x,y)$). Since $avb$ is a 2-path in $\Psi$, the 
involutions fixing both of $a$ and $b$ commute. Therefore if $t'$ 
gives us some $\Theta(u,t',x,y)$-subgraph we are interesting in,
$\Theta(u,t,x,y)\cap\Theta(u,t',x,y)$ contains $\{a,b\}$ properly.
Since $\Pi\subseteq\Theta(u,t,x,y)\cap\Theta(u,t',x,y)\setminus\{a,b\}$,
we are done. \end{pf}

Next, let $w' \in \Theta(t)$. By Lemma~\ref{suz:L8}
$\Gamma(w')\cap\Gamma_2(\infty)=\Gamma(w')\cap\Gamma_2(t')$ for some
$t'\in\Gamma(w')$. Suppose $t\not=t'$. Since $t'\in \V\Theta$,
we then have $\Theta(w,t')\cap\Theta_2(t)\not=\emptyset$, a contradiction.
Thus $t=t'$. 

Thus $\Gamma(t)\subseteq\Gamma_3(\infty)$. Therefore $t\in\Gamma_4(\infty)$.
By standard calculations, we can verify that $\Gamma$ is distance-regular
with the same intersection array as $3\Sigma_6$.

Note that we have shown that any nontrivial connected cover of $\Sigma_6$ is
distance-regular with the same intersection array as $3\Sigma_6$. Therefore
{\em any} (in our sense, see Sect.~\ref{suz:sect1}) 
nontrivial connected cover of $\Sigma_6$ is isomorphic to
$3\Sigma_6$. Hence in order to complete the proof, it suffices to show that 
$\Gamma$ is a cover of $\Sigma_6$.

Since $\Gamma_3(\infty)\cong\Sigma_5\cup\Sigma_5$, to be at the
distance 4 in $\Gamma$ is an equivalence relation on $\V\Gamma$.
Let $t\in\Gamma_4(\infty)$. Then for any $x\in\Gamma(\infty)$ there
exists a unique $y\in\Gamma(t)$ such that $x\in\Gamma_4(y)$. Hence
our relation is well-defined on $\E\Gamma$. Thus we may define the
quotient graph $\overline\Gamma$, whose vertices (resp. edges) are
equivalence classes of vertices (resp. edges) of $\Gamma$. Clearly,
for any $v\in \V\Gamma$ the
restriction of our quotient mapping to $\Gamma(v)$ 
is an isomorphism. Hence $\overline\Gamma$ is locally $\Sigma_5$.
Therefore $\overline\Gamma\cong\Sigma_6$.
The proof in this case was completed.

\Case{$c(\Omega)=2$} It follows from the above that in this case
$v(\Gamma(\infty,u))=64$ for any $u \in \Gamma_2(\infty)$. Now the standard
calculation of the edges coming from $\Gamma(\infty)$ to $\Gamma_2(\infty)$
shows that $v(\Gamma_2(\infty))$ is not integer, a contradiction.

The proof of Theorem~\ref{suz:MTH} is complete.

\medskip
{\bf Acknowledgment.} Several hints about subgraphs of graphs under
consideration were found using the computer package COCO \cite{COCO}.
The author thanks the designer of COCO Igor~Faradzev for help in the
use of it. The author thanks Sergey~Shpectorov who read the first
version of the paper and suggested some improvements to it.

\chapter{Geometric characterization of the sporadic groups $Fi_{22}$,
$Fi_{23}$ and $Fi_{24}$}
\label{chap:fi2n}
\newcommand{\Proof}{\noindent{\bf Proof :}$\quad$}        
\newtheorem{Result}     [thm]   {Result}
\paragraph{Abstract.}
Let $\Si_1,\dots,\Si_4$ be the 3-transposition graphs for Fischer's
sporadic groups $Fi_{21},\dots,Fi_{24}$ \cite{Fi}.
We classify connected locally $\Si_i$ graphs
$\Th=\Th_{i+1}$ for $i=1,2,3$. In the minimal case $i=1$ we also assume
that for every nondegenerate 4-circuit $abcd$ the subgraph
$\Th(a,b,c,d)$ is isomorphic to the disjoint union of three copies of
$K_{3,3}$. If $i=1,2$ then $\Th$ is
isomorphic to $\Si_{i+1}$, whereas if $i=3$ then $\Th$ is isomorphic
either to $\Si_4$ or to its 3-fold antipodal cover $3\Si_4$.

\section{Introduction}

There  have been extensive studies of the geometries of the
classical groups as point-line  systems with fixed local structure.
This approach has led to many  beautiful  results.  We  refer  the
reader  to  a  recent paper of Cohen and Shult \cite{CoSh} for a brief
survey on that. This paper is an attempt to look
at the Fischer's sporadic groups from the aforementioned
point of view.

It  is worth mentioning broad  activity in this area.
Usually certain additional properties are assumed. First,
some  group  action  on  the  geometry  (usually, it is a diagram
geometry) is often assumed. This approach has led 
to study of  some  presentations
for  sporadic  groups, and sometimes to proving the faithfulness of
these presentations. For such studies of the groups under consideration,
see Buekenhout and Hubaut \cite{BH}, 
Meixner \cite{Mei}, Van~Bon and Weiss \cite{vBW:fi:char}.
Another direction is via assuming some  global  property  of  the
geometry,  e.g.  assuming  that  the  related  graph  is strongly
regular    or    distance    regular,    see    the    book    by
Brouwer~et~al. \cite{BCN}.  This  book  contains many results and many
references    to    results    of    this    kind.

As well, the present paper may be viewed as an extension of Hall and
Shult's work \cite{HaSh}, were the same objects arising from classical
groups generated by 3-transpositions have been characterized
axiomatically.

Note that here we  do  not
assume either any group action or global property. Only a few
characterizations of sporadic groups under such weak assumptions are
know to the author. Hughes \cite{Hug1}, \cite{Hug2}, has done this for
Higman-Sims group. Cuypers \cite{Cuy:co2} has characterized Conway's second
group. The author \cite{Pa:suz} has produced such a characterization
for the Suzuki chain of groups. 

It is worth mentioning the concept of {\it Fischer space}, see
e.g. Buekenhout \cite{Bue:fi}, Weiss \cite{We1}, \cite{We2}, 
Cuypers \cite{Cuy:genfisp}, which
is in some sence dual to the concept of 3-transposition graphs. 
Lines of the Fischer space associated with a 3-transposition group $G$ 
are {\it triples} of the 3-transposition graph,
i.e. three points are on a line iff they, regarded as involutions of $G$,
generate $S_3$. Axiomatically, a Fischer space is a partial linear space
with line size 3, whose planes are either affine are dual affine, and
satisfying also certain irreducability condition. Finite Fischer spaces
were classified (in another terms) in \cite{Fi}. It should be stressed
that the class of objects classified in our paper is somewhat larger
than just a special family of duals of Fischer spaces. Indeed,
$3\Sigma_4$ is not the dual of any Fischer space at all. Finally, only
on the final stage of the proof in case of diameter 2 
it becomes clear that our objects admit
3-transposition groups (i.e. are associated with Fischer spaces).

Graphs with constant neighbourhood often arise as collinearity
graphs of diagram geometries, see Buekenhout \cite{Bue:geo:spo}.
One may look at our result as a characterization of 
certain diagram geometries, see the discussion following the statement
of Theorem \ref{fi2n:mainth}.

Throughout the paper we consider undirected graphs without loops
and multiple edges.
Given a graph $\G$, let us denote  the  set  of  vertices  by
$\V=\V\G$,  the  set  of edges by $\E=\E\G$. Let $X\subseteq
\V\G$.
We denote by $\langle X\rangle =\langle X\rangle _{\G}$ the subgraph $\Xi$
of $\G$ {\it induced} by $X\/$
(i.e. $\V\Xi=X,\/\E\Xi=\{ (u,v)\in \E\G\mid u,v\in X\}$).
Given two graphs $\G$ and $\D$, the graph $\G\cup\D$
(resp. the graph $\G\cap\D$) is
the graph with vertex set $\V\G\cup \V\D$ (resp.
$\V\G\cap \V\D$) and edge set $\E\G\cup \E\D$
(resp. $\E\G\cap \E\D$). Given $v\in \V\G$, we denote
$\G_i(v)=\langle \{x\in \V\G\mid x $ at distance $i$ from $v\}\rangle $,
and $\G_1(v)=\G(v)$.
Furthermore, $\G(X)=\bigcap_{x\in X} \G(x)$.
To simplify the notation we use
$\G(v_1,\dots,v_k)$ instead of $\G(\{v_1,\dots,v_k\})$ and
$u\in\G_i(\dots)$ instead of $u\in \V\G_i(\dots)$.
As usual, $v=v(\G)=|\V\G|$,  $k=k(\G)=v(\G(x))$,
$\lambda=\lambda(x,y)=\lambda(\G)=v(\G(x,y))$, where
$x\in\V\G$, $y\in\G(x)$. Let $y\in\G_2(x)$. We denote
$\mu=\mu(x,y)=\mu(\G)=v(\G(x,y))$. Of course, we use
$k$, $\lambda$, $\mu$ if it makes sense, i.e. if those numbers
are independent of the
particular choice of the corresponding vertices.
If $\D$ is a (proper) subgraph of
$\G$ we denote this fact by $\D\subseteq\G$ (resp.
$\D\subset\G$).

We denote  the complete $n$-vertex graph by $K_n$,  the complete
bipartite graph with parts of size $n$ and $m$
by $K_{n,m}$, the circuit of length $n$ by
$C_n$, and the empty graph by $\emptyset$.

Let $\G$, $\D$ be two graphs. We say that $\G$ is
{\it locally} $\D$ if $\G(v)\cong\D$ for each $v\in \V(\G)$.
More generally, if $\D$ be a class of graphs then $\G$ is
locally $\D$ if for each $v\in\V\G$ the subgraph $\G(v)$ is isomorphic
to a member of $\D$.

Let $\G,\ \overline\G$ be two graphs. We say that $\G$ is a
{\it cover} of $\overline\G$ if there exists a mapping $\varphi$
from $\V\G$ to $\V\overline\G$ which maps edges to edges and for
any $v \in \V\G$ the restriction of $\varphi$ to $\G(v)$ is an
isomorphism onto $\overline\G(\varphi(v))$. Note that, since
the latter assumption in the definition of cover is usually
omitted, our definition of cover is a bit nonobvious.

Suppose we have a chain of graphs $\Si_1,\dots,\Si_n$, such that
$\Si_i$ is locally $\Si_{i-1},\  i=2,\dots,n$. Then for any complete
$k$-vertex subgraph $\Upsilon$ of $\Si_m$, $1<k<m\leq n$
we have $\Si_m(\V\Upsilon)\cong\Si_{m-k}$.

We will freely use basic facts concerning {\sl polar spaces}
(see e.g. \cite{BuSh}, \cite{Ti:la}), {\sl generalized quadrangles} (GQ, for short)
(see e.g. \cite{BCN}) and their collinearity graphs. An
{\sl ($i$ times) extended} polar space is a graph which is a 
locally ($i-1$ times)
extended polar space (a 0 times extended polar space is a polar space).
Note that this definition
slightly differs from the common one concerning extended polar spaces,
given, say, in \cite{CHP}. In the
terminology of \cite{CHP} we consider {\sl triangular} extended polar spaces.
         
A {\sl 3-transposition group} is a group containing a conjugacy class $C$
of involutions (called 3-transpositions) such that for any $a,b\in C$, one
has either $(ab)^2=1$ or $(ab)^3=1$. The {\sl 3-transposition graph} associated with a 3-transposition group is the
graph with vertex set $C$, two vertices are adjacent if the corresponding
involutions commute. For a graph $\G$,
$\Aut(\G)$ denotes the automorphism group of $\G$.
Our group-theoretic notation is as in \cite{Atl}.

Let $\Si_1,\dots,\Si_4$ be the 3-transposition graphs for Fischer's
sporadic groups $Fi_{21},\dots,Fi_{24}$ \cite{Fi}. Then $\Si_i$ is
an $(i-1)$-fold extension of the polar space $\G$, where $\G$ is the polar space
arising from a 6-dimensional GF(4)-vector space carrying a nondegenerate
hermitian form. Also
$\Si_4$ admits a unique (antipodal) 3-cover $3\Si_4$, which is a distance
transitive graph in the sense of \cite{BCN}.

\begin{thm}
Let $\Th=\Th_{i+1}$ be a connected locally $\Si_i$ graph ($i=1,2,3$).
In the minimal case $i=1$ we also assume that,
for any $C_4$-subgraph $\Xi$ of $\Th$, one has
$\Th(\V\Xi)\cong K_{3,3}\cup K_{3,3}\cup K_{3,3}$.
If $i=1,2$ then $\Th\cong\Si_{i+1}$. If $i=3$ then
$\Th\cong\Si_4$ or $3\Si_4$.
\label{fi2n:mainth}
\end{thm}

One may view the graphs $\Th_{i+1}$ as the collinearity graphs
of certain $c^i\cdot C_3$-geometries ${\cal G}(\Th_{i+1})$, 
namely rank $i+3$ geometries with diagram
$$\circ \!{\phantom{00000}\over \phantom{0}}
\!\!\circ\cdots\circ \!\!{\phantom{00} c \phantom{00}\over \phantom{0}}
\!\!\circ \!\!{{\phantom{00000}}\over \phantom{0}}
{\vphantom{\bigl\langle}}_{_4}
\!\!\!\circ \!\! = \!\! = \!\! = \!\! = \!\! =
\!{\vphantom{\bigl\langle}}_{_2} \!\!\!\circ$$
Here $1\leq i\leq 3$.
Conversely, the elements of the geometry may be viewed as $k$-cliques
of the graph with natural incidence, $k=1,2,\dots, i+1, i+6, i+21$.

Note that $(i+21)$-cliques are of special interest, since it is well-known
that the geometries
induced on them are $i$ times extended projective planes of order 4. Such
extensions have been studied widely. In particular it has been shown that
$i$ may not exceed 3. This implies that $\Si_4$ and $3\Si_4$ do not posess
further extensions. It is well-known that, for a given $i$, the $i$-extension
of a projective plane of order 4 is
unique (i.e. one obtains Steiner systems $S(3,6,22)$, $S(4,7,23)$ and
$S(5,8,24)$ respectively).

Meixner~\cite{Mei} has proved the following result, which also follows
from Van Bon and Weiss \cite{vBW:fi:char}.

\begin{Result}
Let ${\cal G}$ be a residually connected flag-transitive
$c^i\cdot C_3$-geometry, $1\leq i \leq 3$.
Then if the $C_3$-residues of {\cal G} are isomorphic to
${\cal G}(\Si_1)$
then ${\cal G}={\cal G}(\Si_{i+1})$ (or ${\cal G}(3\Si_4)$ if $i=3$).
\label{fi2n:Meires}
\end{Result}

In his proof Meixner first shows that the collinearity graphs of the
geometries under consideration are locally $\Si_i$, then  he classifies
locally $\Si_i$ graphs having a certain automorphism group. We currently
are almost able to repeat the second part of his work without any assumption on
group action (almost, because if $i=1$ we assume a bit more than local
structure, cf. the statement of Theorem \ref{fi2n:mainth}). Buekenhout  and Hubaut
\cite{BH} have proved a result similar to Result \ref{fi2n:Meires} for $i=1$
under somewhat stronger assumptions (see also Del Fra, Ghinelli,
Meixner, Pasini \cite{DFGMP}).

\section{Proof of Theorem}
Let $\Th=\Th_{i+1}$ be a connected locally $\Si_i$ graph $(i=1,2,3)$.
The crucial point in our proof is the recognition of the
subgraphs $\Th(u,v)$ of $\Th$ induced on the common neighbourhood 
of two vertices $u$, $v$ at
distance two. The following well-known general fact (see, e.g. \cite{Pa:suz}),
has been used.

Let $\D$ be a connected graph satisfying the following property.

\begin{description}
\item[{\rm (*)}] For any $u\in \V\D$ and
$v\in \V\D\setminus (\V\D(u)\cup\{ u \})$
the subgraph $\D(u,v)$ is isomorphic to some $M_{\D}$, whose
isomorphism type is independent of the particular choice of $u$ and $v$.
\end{description}

\begin{lem} 
Let $\Ga$ be a locally $\D$ graph, where $\D$
satisfies {\rm (*)}. For any $u\in \V\Ga$, $v\in\Ga_2(u)$
the graph $\Ga(u,v)$ is locally $M_{\D}$. \qed
\label{lsmugraphs}
\end{lem}

Note that $\D=\Si_i$ satisfies (*).
Indeed, the stabilizer
of $u\in \V\D$ in $\Aut(\D)$ acts transitively on
$\D_2(u)$. This implies (*).
Thus Lemma \ref{lsmugraphs} holds for locally $\D$ graphs $\Th$.
We classify the locally $M_\D$ subgraphs of $\D$.
It turns out that these subgraphs are (disjoint unions of) $i$ times
extended GQ(4,2).
For $i=1$ \cite{BBegq42} (resp. for $i>1$ \cite{Pa:3tr}) gives us the
isomorphism type of these extended GQ.

Information on $\mu$-graphs gives us some control over the second
neighbourhood of $\Th$ as well.
Indeed, adjacency in $\Th_2(u)$ may be interpreted in terms of
intersections of two $\mu$-graphs, once we know the relationship
between vertices of $\D$ and $M_\D$-subgraphs.

The aforementioned outline for studying $\mu$-graphs of $\Th$ is worked out in
Subsection 2.1. In the remaining subsections the question of reconstruction
of the second neighbourhood of $\Th$ is considered.  It turns out that
for each $v\in\Th_2(u)$ there exists a unique $w\in\Th_2(u)$ such that
$\Th(u,v)=\Th(u,w)=\Th(v,w)$ (i.e. $\Th$ is a {\sl triple} graph).
Except a certain subcase of the case $i=3$, we find that $\Th$ has diameter
two, and we complete the proof along the lines of \cite{Pa:3tr}, constructing
certain automorphisms of $\Th$ which leave $u$ and $\Th(u)$ fixed and
interchange $v$ and $w$ in triples. Finally, we consider
the remaining subcase of case $i=3$, and show that we get a (triple) cover
of $\Si_4$ as the only other possibility. There the author adopts many ideas
from his paper \cite{Pa:suz}.

It is worth noticing, though we will not make use of these results, that
the covers of $\Si_i$ were classified by Ronan \cite{Ro:geo:covers}.

\subsection{Preliminaries}
\label{Prel}
Here we present technical
statements concerning certain subgroups of Fischer's groups and
graphs associated with.
Set $\O_i$ to be 
the graph of $(+)$-points of the
$(i+5)$-dimensional vector space over {\rm GF(3)}
carrying a nondegenerate bilinear form with discriminant 1, points are
adjacent if and only if 
they are perpendicular with respect to this form, for
$i=1,2,3$.
Equivalently, $\O_i$ is the 3-transposition graph for the group
$O_{i+5}^\varepsilon (3)$, where 
the Witt index $\varepsilon$ of equals $-$, empty or $+$
according as $i=1,2,3$.

The following statement is well known.
\begin{lem}
Graphs $\O_1$, $\O_2$, $\O_3$ are strongly regular graphs 
with parameters $(v,k,\lambda,\mu)$ equal to 
$$(126,45,12,18),\ (351,126,45,45),\ (1080,351,126,108)$$ 
respectively.\qed
\label{Oparam}
\end{lem}

Consider the case $i=1$. Let $\G$ be the collinearity graph of
$U_6(2)$-polar space, $H\cong O_6^-(3)$ be
a subgroup of $G\cong Fi_{21}\cong U_6(2)$. There are three conjugacy
classes of such
subgroups in $U_6(2)$. Each such subgroup $H$ has a unique 
orbit $\O_H$ of length 126 on
$\V\G$ (note that the subgraph induced on this 126-orbit is isomorphic
to $\O_1$). Let $O_k$ $(k=1,2,3$ refers to the conjugacy classes of the
$O_6^-(3)$-subgroups of $G)$ be the orbits of the action of $G$ on the set of
these 126-orbits, regarded as subsets of $\V\G$. 
Clearly, this action coincides with the action of $G$
on the union of conjugacy classes of $O_6^-(3)$-subgroups.

\begin{lem} 
Let $\O=\O_H\in O_1$.
$H=G(\O)$ has orbits $O_{1j}$ ($j=1,2,3$) on $O_1$ of
lengths $1$, $567$, $840$ respectively.
{\begin{enumerate}
 \item Let $\O'\in O_{1j}$. Then
$|\O\cap\O'|=30$ or $6$ according as $j=2$ or $3$.
 \item Let $k=2$ or $3$. Then $H$ has orbits $O_{kj}$ ($j=1,2$) on $O_k$ of
lengths $112$ and $1296$ respectively. Let $\O'\in O_{kj}$. Then
$|\O\cap\O'|=45$ or $21$ according as $j=1$ or $2$.
\end{enumerate}}
\label{intEGQ(4,2)}
\end{lem}

\Proof We exploit the following facts about the intersections of
$O_6^-(3)$-subgroups of $G$. Here we prefer to work with the subgroups
of index two of $H$ and its conjugates in $\Aut(\G)$. Such a
subgroup in $G(\O)$ we will denote by $F(\O)$. \hfil\break
1) If $\O'\in O_1$ then $F(\O)\cap F(\O')\cong F(\O)$, $2^4:A_6$ or
$3_+^{1+4}.2S_3$.\hfil\break
2) If $\O'\in O_k$ ($k=2,3$) then $F(\O)\cap F(\O')\cong 3^4:A_6$ or $A_7$.
\hfil\break
This follows from the information  found in \cite{Atl} about concrete types
of subgroups normalized by outer automorphisms of $U_6(2)$, etc.

Now for any possible $X=F(\O)\cap F(\O')$ we look at the action of $X$ on
the points of $\G$ inside and outside $\O$. Using a computer, we calculate
the orbits of this action and find that there is only one possiblity for
$|\O\cap\O'|$. (For instance, for $X\cong A_7$ we have found that $X$ has orbits
of length 21 and 105 on $\O$ and has no orbits of length 21 outside $\O$.)
\qed
                           
Consider the case $i=2$. Let $\G$ be the 3-transposition graph for
$G\cong Fi_{22}$, $H\cong O_7(3)$ be a subgroup of $G$.
There are two classes of such subgroups in $G$. Each such subgroup
$H$ has a unique orbit $\O_H$ of length $351$ on $\V\G$.
Let $O_k$ ($k=1,2$ corresponds to the conjugacy classes of
$O_7(3)$-subgroups of $G$) be the orbits of the action of $G$ on the set of
these 351-orbits, regarded as subsets of $\V\G$.

\begin{lem}
Let $\O=\O_H \in O_1$.
$H$ has orbits $O_{1j}$ (resp. $O_{2j}$) ($j=1,2,3$) on $O_1$ (resp.
on $O_2$) of lengths $1$, $3159$, $10920$ (resp. $1080$, $364$, $12636$)
respectively.
{\begin{enumerate}
\item Let $\O'\in O_{1j}$. Then $|\O\cap\O'|=63$ or $7$ according as
$j=2$ or $3$.
\item Let $\O'\in O_{2j}$. Then $|\O\cap\O'|=0$, $108$
or $36$ according as $j=1$, $2$ or $3$.
\end{enumerate}}
\label{intO73}
\end{lem}

\Proof
Action of $H$ on $O_1$ (as well as the corresponding intersections)
can be treated modulo the construction of $Fi_{23}$ as a transitive
extension of $Fi_{22}$.
Action on $O_2$ and correponding intersections follow from information
on the intersections given in \cite{Atl} and by induction from Lemma
\ref{intEGQ(4,2)}. Empty intersection appears, in particular, because
there are only two types of intersection in Lemma \ref{intEGQ(4,2)}.
\qed

The following statement is well-known.
\begin{lem}
Graphs $\Si_1$, $\Si_2$, $\Si_3$ are strongly regular graphs 
with parameters $(v,k,\lambda,\mu)$ equal to 
$$(693,180,51,45),\ (3510,693,180,126),\ (31671,3510,693,351)$$ 
respectively.\qed
\label{Fiparam}
\end{lem}

\subsection{$\mu$-graphs}
Here $\G=\Si_i$ ($i=1,2,3$).
Note that in the case $i=1$ we are unable to classify all the locally $M_\G$-subgraphs
of $\G$. However, our additional assumtion on $\Th$ says that for any
$\mu$-graph $\O$ of $\Th$ and for any pair $(x,y)$ of nonadjacent vertices
of $\O$ the subgraph $\O(x,y)$ is isomorphic to disjoint union of three
copies of $K_{3,3}$. Thus we are able to restrict our attention on the
locally $M_\G$-subgraphs of $\G$ satisfying the same condition as $\O$ does.

The following statement gives us a control over the isomorphism types of
locally $M_\G$-subgraphs $\O$ of $\G$.

\begin{pr}
Let $\O$ be a locally $M_\G$ subgraph of $\G$ satisfying for $i=1$
the additonal condition that for any pair of $(x,y)$ of nonadjacent vertices
$\O(x,y)\cong K_{3,3}\cup K_{3,3}\cup K_{3,3}$. Hence $\O$ is isomorphic
to (a disjoint union of if $i=3$) the graph(s) $\O_i$ 
defined in Subsection \ref{Prel}.
\label{muOn3}
\end{pr}

\Proof
Since $M_\G$ for $i=1$ is isomorphic to GQ(4,2),
this case immediately follows from \cite{BBegq42}. Thus we know
the isomorphism type of $\O_1$. Since $\mu$-graphs
of $\Si_i$ ($i<4$) are connected, $\O_i$ is locally $\O_{i-1}$ for
$1<i<4$. Thus in the remaining cases the classification of the
graphs with this neighbourhood structure given in \cite{Pa:3tr} is applicable.
(Since \cite{BBegq42} has not been published, it is worth noticing that a proof for
the case $i=1$ can be easily deduced from \cite{Pa:3tr}, as well).
\qed

Next, we establish the connection between the certain subgroups of
$Fi_{2i}$ and the subgraphs $\O$, whose isomorphism types $\O_i$ were
determined in Proposition \ref{muOn3}.

\begin{pr}
Let $\O\cong\O_i$ be a subgraph of $\G$. Then there exists
$H\cong O_{i+5}^\varepsilon (3)$-subgroup of $\Aut(\G)$ such that
$\V\O$ constitutes an orbit of $H$ in its action on $\V\G$, where
the Witt index $\varepsilon$ of $H$ equals $-$, empty or $+$
according as $i=1,2,3$.
\label{embedmu}
\end{pr}

\Proof
The following proof is suggested by the referee. The original one
was much longer.

It is enough to check the claim 
that $\O$ is a triple subgraph of $\G$. Indeed, in
this case the subgroup generated by the 3-transpositions
corresponding to $\V\O$ is a subgroup of $\Aut(\G)$, and may be easily
identified with $O_{i+5}^\varepsilon (3)$, as required.

Now let us prove the claim for $i=1$. 
Assume that $\{u,v,w\}$ is a triple of $\O$, but not of $\G$.
Hence, by Lemma \ref{Oparam}, $|\O(u,v,w)|=18$. 
Clearly, $|\G(u,v,w)|\ge |\O(u,v,w)|$. 
By a well-known property of polar spaces,
the points of $\G(u,v,w)$ constitute a proper geometric
hyperplane of GQ(4,2) induced on $\G(u,v)$, which cannot have more then
13 points.  Hence $\{u,v,w\}$ is a triple of $\G$.
This completes the proof in for $i=1$.

Let $i\ge 2$. Assume $\{x,y,z\}$ is a triple of $\O$. Fix a vertex 
$u\in\O(x,y,z)$. By induction, $\{x,y,z\}$, being a triple of $\O(u)$,
is a triple of $\G(u)$. But triples of $\G(u)$ are also triples of $\G$,
and we find that $\{x,y,z\}$ is a triple of $\G$. 
\qed

Let us turn to the relationship
between the vertices of $\G$ inside and outside a connected locally
$M_\G$ subgraph $\O$ of $\G$.

\begin{lem}
Let $v$ be a point of $\G$ outside $\O$, $L=\G(v)\cap\O$.
Then, according to the value of $i$,
\begin{enumerate}
\item
$|L|=30$, and $L$ is isomorphic to the 2-clique extension $L_1$ of
the point graph of GQ(2,2).
\item
$|L|=63$, and
$L$ is isomorphic to  the collinearity graph $L_2$ of $Sp_6(2)$-polar space.
\item
If $L$ is not empty then $|L|=120$, and
$L$ is isomorphic to  the collinearity graph $L_3$ of
$O_8^+(2)$-affine polar space.
\end{enumerate}
\label{ptOrel}
\end{lem}

\Proof Here $H$ denotes an automorphism group of $\O$ provided by
Proposition \ref{embedmu}. Below we denote $\L=\G(v)\cap\O$.
                                           
{\bf Case i=1.}
Since we know that $H$ has two orbits on
the points of $\G$, the first part of the  lemma follows from a
simple counting of the $\G$-edges going from $\O$ outside.

Since each maximal clique of $\O$ is a hyperoval in a plane of $\G$,
each line through $v$ intersects $\O$
in 0 or 2 points. Thus $L$ is the 2-clique extension of a
subgraph $\overline{\L}$ of the point graph $\overline{X}$ of GQ(4,2).
It is clear that each line of $\overline{X}$ intersects $\overline{\L}$ in
0 or 3 points. Thus $\overline{\L}$ has {\it lines} of size 3. Moreover,
if $l$ is such a line then $l\cup\overline{\L}(l)=\V\overline{\L}$. It
means that $\overline{\L}$ is the point graph of a GQ, evidently of
GQ(2,2).

{\bf Case i=2.}
Since we know that $H$ has two orbits on
the points of $\G$, the first part of the  lemma follows from a
simple counting of the $\G$-edges going from $\O$ outside.

By the case 1 of the current lemma each connected component of
$L$ is isomorphic to a graph such that  for each 
vertex the neighbourhood graph is the same as for $Sp_6(2)$ polar space.
According to J.~Tits \cite{Ti:la}, such a graph is unique. It follows from the
first part of the lemma that $L$ is connected.

{\bf Case i=3.}
Here $H$ is not maximal in $G\cong Fi_{23}$. According to \cite{Atl},
it is contained in a maximal subgroup isomorphic to $H:S_3$, which has
exactly two orbits $W$, $W'$ on $\V\G$. 
It is straightforward to check
that one of those orbits, say $W$, has length
$3\cdot |\V\G|$. Therefore $\O\subset W$, and $W$ contains two other copies
$\O'$ and $\O''$ of $\O$. 

If $(x,y)$ is an edge between, say, $\O$
and $\O'$, then, by the case 2 of the current lemma, 
the subgraph  $\O'(x)$ has valency
63. On the other hand, by \cite{Atl}, the stabilizator of $x$ in $H$ acts
on $\O'$ transitively, so $\O'(x)=\O'$, contradiction. 
Hence there are no edges between the copies of $\O$ in $W$.

Let $Y=\G(v)\cap W$. By the same argument
as in previous cases, we obtain that $|Y|=360$. Now, assuming that
there exists an element $z$ of order 3 from $(H:S_3)\setminus H$ fixing $v$
and interchanging the copies of $\O$ inside $W$, we immediately get that
$|L|=360/3=120$. A look at \cite{Atl} gives us that each 3-element of $G$
cannot either fix $W$ pointwise or act fixed point free on $\V\G$. Hence such
$z$ indeed exists, and the first part is proved.

By the case 2 of the current lemma each connected component of
$L$ is isomorphic to a graph such that  for each $x\in\V L$
the subgraph $L(x)$ is isomorphic to the collinarity
graph of $Sp_6(2)$-polar space.
By \cite{HaSh}, or \cite{CoSh}, 
there are exactly two such graphs, and only one of them
has the required number of vertices. It follows from the
first part of the lemma that $L$ is connected.
\qed
                                             
Now we are able to complete the classification of locally $M_\G$-subgraphs
of $\G$. The following statement is immediate consequence of
Lemma \ref{ptOrel}.

\begin{lem}
Let $\O$ be a locally $M_\G$-subgraph of $\G$. If $i=1,2$ then
$\O$ is connected. If $i=3$ there is a unique equivalence relation $\phi$
with class size three on the set of connected locally $M_\G$ subgraphs,
such that all the connected components of $\O$ lie in one class of $\phi$.
\qed
\label{locmugraphs}
\end{lem}

Finally, we need the following technical statement.

\begin{lem}
Let $\O$ be a locally $M_\G$-subgraph of $\G$, the vertices $v$ and $w$ of
$\G$ lie outside $\V\O$. Then $\G(v)\cap\O=\G(w)\cap\O$ implies $v=w$.
The subgraph of $\G$ induced on $\V\G\setminus\V\O$ is connected.
\label{Bij}
\end{lem}

\Proof
Follows from the primitivity of the action of the stabilizer of $\O$ in
$\Aut(\G)$ on the points outside $\O$.
\qed

\subsection{Final part of the proof}

Let $\Xi$ be the graph defined on the set of connected locally $M_\G$
subgraphs of $\G$, two vertices $\O$, $\O'$ are adjacent if
$\O\cap\O'\cong L_i$ ($L_i$ are defined in Lemma \ref{ptOrel}).
Let $\Th$ be a connected locally
$\G$ graph, $u\in\V\Th$, $v\in\Th_2(u)$.
           
{\bf Cases i=1 and 2.}
We start with a technical statement.

\begin{lem}
$\Xi$ has $4-i$ connected components, which are the orbits of $Fi_{2i}$
on the set of locally $M_\G$-subgraphs of $\G$.
$k(\Xi)=567$ or $3159$ according as $i=1$ or $2$.
\label{grphsonO12}
\end{lem}

\Proof
Immediately follows from Lemmas \ref{intEGQ(4,2)}, \ref{intO73}.
\qed

Now are able to reach our goal. By Lemma \ref{ptOrel},
there are no edges coming from $v$ to $\Th_3(u)$, that is, the diameter of
$\Th$ equals two. Lemma \ref{Bij} provides a bijection between
$\Th(v)\cap\Th_2(u)$ and $\Xi(\Th(u,v))$. 
Therefore $\Th$ is strongly regular with the same
parameters as $\Si_{i+1}$-graph. It follows from the fact that $\G$ is a
triple graph and the counting of $\mu$-graphs 
that $\Th$ is a triple graph. Therefore $\Th_2(u)$ is a
double cover of a connected component of $\Xi$.
Since this is a double cover, the covering map acts on
$\Th_2(u)$ fixed point freely. Thus for each $u\in\V\Th$ there exists an
involution $g_u$ fixing $\Th(u)$ pointwise and acting fixed point freely on
$\Th_2(u)$. The rest of the proof is straightforward along the lines of
\cite{Pa:3tr} and consists of identifying the group generated
by $g_w$ ($w\in\Th(u)$) with $\langle g_u \rangle .Fi_{2i}$ and applying
\cite{Fi} to identify the group generated by all the $g_u$ with
$Fi_{2,i+1}$.
\medskip

{\bf Case i=3.}
Let us recognize the situation with intersections of locally
$M_\G$-subgraphs in this case.

Let $H\cong O_8^+(3):S_3$ be a subgroup of $F=Fi_{23}$.
There is a unique class of such subgroups in $F$. Each such subgroup
has an orbit $\overline\O$ 
of length $3\cdot 1080$ on $\V\G$, whose elements are the vertices of
the subgraphs belonging to 
$\phi(\O)$ for a connected
locally $M_\G$ subgraph $\O$.
Let $O$ be the orbit of the action of $F$ on the set of
these $3\cdot 1080$-orbits regarded as subsets of $\V\G$.

\begin{lem}
Let $\Pi=\overline\O\in O_1$.
$H$ has the orbits $O_j$  on $O$
of lengths $1$, $28431$, $109200$ respectively ($j=1,2,3$).\hfil\break\indent
Let $\Pi'\in O_2$. Then $|\Pi\cap\Pi'|=360$.
The subgraph induced on $\Pi\cap\Pi'$
is isomorphic to the disjoint union of three copies of $L_3$.
Moreover, for any $\O\in\phi(\O)$ there exists unique
$\O'\in\phi(\O')=\Pi'$
such that $\O\cap\O'\cong\ L_3$.\hfil\break\indent
Let $\Pi'\in O_3$. Then $8\le |\Pi\cap\Pi'|<360$.
\label{intO83S3}
\end{lem}

\Proof
Action of $H$ on $O$ (as well as the corresponding intersections)
can be treated modulo the construction of $Fi_{24}$ as a transitive
extension of $Fi_{23}$.
\qed

Let $\overline\Xi$ be the graph defined
on $O$ such that two vertices $\phi(\O)$, $\phi(\O')$ are adjacent if
for some $\O_1\in\phi(\O)$ there exists $\O_1'\in\phi(\O')$ such that
$\O_1\cap\O_1'\cong L_3$. By Lemma \ref{intO83S3}, actually for {\it any}
$\O_1\in\phi(\O)$ such a $\O_1'$ exists. Clearly, $\overline\Xi$ is connected.
Thus $\overline\Xi$ is the quotient
of $\Xi$ defined by $\phi$, where $\Xi$ is defined in the beginning of the
current subsection.
\begin{lem}
$\Xi$ is connected, $k(\Xi)=28431$.
\label{grphsonconO}
\end{lem}
\Proof
Clearly $k(\Xi)=k(\overline\Xi)=28431$. If $\Xi$ were disconnected,
it would happen only if it has 3 connected components, but it is impossible,
since the simple group $Fi_{23}$ acts vertex transitively on it.
\qed

We need one more technical statement.
\begin{lem}
Let $\G\cong\Si_3$. Suppose that a subgraph $\D$ of $\G$ satisfies the following
property: for any $x\in\V\D$ one has $\G(x)\setminus\D(x)=\G(x,t^x)$ for
some $t^x\in\G_2(x)$. Then $\D=\G_2(t)$ for some $t\in\V\G$.
\label{last}
\end{lem}
\Proof
The major part of the following argument was suggested by the referee.

Fix $v\in\V\D$ and $t=t^v$. Let $t'$ be the third vertex in the triple
of $\G$ containing $v$ and $t$. The subgraphs $\G_2(t)$ and $\G_2(t')$
satisfy the conditions of the lemma. We show that $\D$ coincides with
one of $\G_2(t)$ and $\G_2(t')$.

Let $x\in\D(v)$. There are exactly two triples $\{x,t_{11},t_{12}\}$ and
$\{x,t_{21},t_{22}\}$ such that $\G(x,v,t,t_{\alpha,\beta})=\G(x,v,t)$,
where $\alpha,\beta=1,2$. Indeed, by Lemmas \ref{intO73} and \ref{ptOrel},  
$\G(x,v,t)$ is a maximal-by-size possible proper intersection of
$\G(x,a)$ and $\G(x,b)$ for $a,b\in\G_2(x)$.
Hence $\G(x,v,t,t_{\alpha,\beta})=\G(x,v,t)=\G(x,t,t_{\alpha,\beta})$,
where $\alpha,\beta=1,2$.  Therefore for $\alpha,\beta=1,2$ one has that 
$t_{\alpha,\beta}$ belongs to the Fischer subspace $\Pi$ generated
by $x$, $v$ and $t$, which is isomorphic to the dual affine plane of
order 2. Thus, the triples under question are the
two triples of $\Pi$ through $x$. 

Moreover, these triples $\{x,t_{11},t_{12}\}$, $\{x,t_{21},t_{22}\}$
meet the triple $\{v,t,t'\}$ in $\Pi$. So, without loss of generality we
have $t=t_{11}$, $t'=t_{21}$. 

Thus, without loss of generality we can take $t^x$ to be $t$.

Let $y\in\D(v,x)$. Then, repeating the arguments of the first part of
the proof with $y$ and $v$ or with $y$ and $x$, instead of $v$ and $x$,
we see that the triple $\{y,t^y,(t^y)'\}$ has to meet both the triples
on $v$ and $t$ and on $x$ and $t$, and thus contains $t$. Now, by
connectedness of $\D(v)$ and $\G_2(z)$ 
for any $z\in\V\G$, we find that all triples obtained in this way
contain $t$ and $\D=\G_2(t)$.
\qed

We start with identifying $\G_2(u)$. As a by-product, we settle the case
leading to $\Th\cong\Si_4$.
 
By Lemma \ref{locmugraphs},
the number of connected components of $\Th(u,v)$ equals
$m=m_v=m_{u,v}$ for $1\leq m\leq 3$. Let $w\in\Th(v)\cap\Th_2(u)$,
$\Th(u,v,w)\not=\emptyset$. Clearly $m_w=m$. By Lemma \ref{ptOrel},
$\Th(u,v,w)$ is isomorphic to the disjoint union of $m$ copies of
$L_3$. Set $\L=\L^v$ to be a connected (partial) subgraph of $\Th_2(u)$
containing $v$, two vertices $x$, $y$ are joined if $\Th(u,x,y)\cong\Th(u,v,w)$.
Lemma \ref{Bij} provides the natural bijection between $\L(v)$ and the
set of locally $M_\G$-subgraphs $\O$ of $\G=\Th(u)$ such that 
$\O\cap\Th(u,v)\cong\Th(u,v,w)$. 

Let $m\ge 2$. Then $\L$ is a cover of $\overline\Xi$ defined by
$y\mapsto\Th(u,y)$ for $y\in\V\L$. (This is straightforward for $m=3$.
For $m=2$ it follows from the fact that for any $\phi$-equivalence class
of connected locally $M_\G$-subgraphs of $\G$ exactly two of the three
members of the class appear as connected components of $\Th(u,y)$,
$y\in\V\L$.) Since
\[\mu(\overline\Xi)=5832>3\cdot 1080=
\sup_{x\in\V\Th,\ \mbox{\scriptsize $y\in\Th_2(x)$}}|\Th(x,y)|\]
(cf. \cite{BvL} on the value of $\mu(\overline\Xi)$),
$\L$ is either a connected proper cover, or all the distance two
vertices of $\L$ are joined by edges inside $\Th_2(u)$. Since the latter
implies $|\Th_2(u)\cap\Th(v)|>k(\Th)$, it is nonsense.
By consideration of the neighbourhood of a vertex from $\Th(u,v)$, the
foldness of this cover is two.

Let $m=3$. In this case, by Lemmas \ref{ptOrel} and \ref{locmugraphs},
$\Th(v)=\L(v)\cup\Th(u,v)$. By counting the edges between $\Th(u)$ and
$\Th_2(u)$, we have $\Th_2(u)=\L$.
The same arguments as in the already considered
cases $i=1$ or $2$ complete the proof of $\Th\cong\Sigma_4$.

Thus we may assume from now on that $m\leq 2$. Let $m=2$. Counting edges
between $\Th(u)$ and $\V\L$, we see that
\begin{equation} 
|\Th_2(u)|<2|\V\L|<2|\Th_2(u)|.
\label{sizeTh2}
\end{equation}
Let $v'\in\Th_2(u)\setminus\V\L$. Since $|\V \mbox{$\L^{v'}$} |<|\V\L|$, we have
$m_{v'}=1$. But $\L^{v'}$ is a (possibly, improper) cover of $\Xi$ defined
by $y\mapsto\Th(u,y)$ for $y\in\V \mbox{$\L^{v'}$}$. So $|\V \mbox{$\L^{v'}$} |\ge
|\V\Xi|>|\V\L|$, a contradiction to (\ref{sizeTh2}).

Hence $m=2$ is nonsense. In what follows we assume $m=m_{x,y}=m(\Th)=1$ does
not depend on the particular choice of $x\in\V\Th$ and $y\in\Th_2(x)$.
As we have already observed, $\L$ is a cover of $\Xi$.
Since
\[\max_{a\in\V\Th,\ \mbox{\scriptsize $b\in\Th_2(a)$}}|\Th(a,b)|\ge\mu(\overline\Xi)/3>\mu(\Xi),\]
$\L$ is either a connected proper cover, or all the distance two vertices  of
$\L$ are joined by the edges inside $\Th_2(u)$. By the same argument as
in the case $m\ge 2$, the latter is impossible.
By consideration of the neighbourhood of a vertex from $\Th(u,v)$, the
foldness of the cover is two.

Now we claim that $\Th$ is a
distance regular graph with the same interesection array as
$3\Sigma_4$. 

First, let us show that $\L=\Th_2(u)$. Suppose false, i.e. there
exists an edge $(x,y)$ in $\Th_2(u)$, which does not belong to $\L$. 
Clearly, $\Th(u,x)\not=\Th(u,y)$.
Hence
$v(\Th(u,x,y))=0$. Consider $\Th(x)$. The vertex $y$ belongs to a subgraph $\D$
from $\phi_x (\Th(x,u))$. Hence there are at most
$k(\D)=351$ vertices of $\Th(x,y)$ belonging to $\Th_3(u)$.
Moreover, for any vertex $w$ from $\Th(x,y)$ not in this 351-set $\D(y)$,
we have $\Th(u,x,w)\cong\Th(u,y,w)\cong L_3$, i.e. $(x,w)$ and
$(y,w)$ belong to $\E\L$. Since $\Th(x,y)\cong\Sigma_2$, we
have that there are at least 2808 such vertices $w$. Recall
$\L$ is a cover of $\Xi$. Hence there exists a pair of distance
two vertices of $\Xi$, which have at least 2808 common neighbours.
Therefore $\mu(\overline\Xi)=5832\geq 3\cdot 2808$, a contradiction.

Let $v\in\Th_2(u)$, $w\in\Th_3(u)\cap\Th(v)$. Denote $\Ga=\Th(w)$.
We claim that $\D=\Th_2(u)\cap\G$ satisfies the conditions of Lemma
\ref{last}. 
From the above consideration of $\Th_2(u)$, we know that for any
$x\in\Th_2(u)$ the subgraph $\Th_3(u)\cap\Th(x)$ is the disjoint union
of the two subgraphs from $\phi_x(\Th(x,u))$ distinct from $\Th(x,u)$.
Hence for any $x\in\V\D$ the subgraph $\G(x)\cap\Th_3(u)$ is isomorphic
to $\G(x,t^x)$, where $t^x\in\G_2(x)$. Thus Lemma \ref{last} is
applicable, and $\D=\G_2(t)$ for some $t\in\V\G$.

Let $w'\in\Ga(t)$. By the same arguments,
$\Th(w')\cap\Th_2(u)=\Th(w')\cap\Th_2(t')$ for some
$t'\in\Th(w')$. Suppose $t\not=t'$. Note that $t'\in\V\Ga$.
Clearly $t'\not\in\Ga_2(t)$. Since $\lambda(\Ga(t))=693<v(\Ga(w',t'))-1$,
we have $\Ga(w',t')\cap\Ga_2(t)\not=\emptyset$, 
contradiction to the definition of $t'$.

Thus $\Th(t)\subseteq\Th_3(u)$. Hence $t\in\Th_4(u)$. Now a standard 
counting of edges between $\Th_3(u)$ and $\Th_4(u)$ shows that $\Th$ has
the same intersection array as $3\Sigma_4$. 
It remains to show that $\Th$ is a cover of $\Sigma_4$.
This is enough to complete the proof. Indeed, it implies
that the universal cover of $\Sigma_4$, which has the same
neighbourhood as $\Sigma_4$, has the same intersection array as
$3\Sigma_4$. So this universal cover is isomorphic to 
$\Th$, and $\Th\cong\ 3\Sigma_4$.

Since $\Th_3(u)\cong\Sigma_3\cup\Sigma_3$, to be at distance 4 in $\Th$ is an 
equivalence relation on $\V\Th$. Then it is straightforward to check that 
for any $x\in\Th(u)$ there exists unique $y\in\Th(t)\cap \Th_4(x)$. 
Thus our equivalence relation is well-defined on $\E\Th$, and
we may define the quotient graph, whose vertices (resp. edges) are equivalence
classes of vertices (resp. edges) of $\Th$. Clearly, the restriction of our 
quotient map to the neighbourhood of any $z\in\V\Th$ is an isomorphism. 
Therefore our quotient graph is locally $\Sigma_3$ and of diameter 2, 
hence isomorphic to $\Sigma_4$.

The proof of the Theorem is completed.

\medskip
{\bf Acknowledgement.} The author thanks the referee for numerous
comments and improvements, particularly concerning proofs of Proposition
\ref{embedmu} and Lem\-ma \ref{last}.

\chapter{Extending polar spaces of rank at least 3}
\label{chap:epolsp}
\newcommand{\OO}{\Pi}
\renewcommand{\a}{\alpha}
\renewcommand{\b}{\beta}
\newcommand{\g}{\gamma}
\renewcommand{\d}{\delta}
\newcommand{\oOO}{\overline{\Pi}}
\paragraph{Abstract.}
It is shown that the extended polar space for the 
sporadic Fischer group $Fi_{22}$ is the only extended polar space which 
has more than two extended planes on a block
and is not isomorphic to a quotient of an affine polar space over GF(2).

New examples of EGQ(4,1) and EGQ(4,2) are presented, as well. 

\section{Introduction and the results}
This paper is a continuation of an earlier work \cite{Pa:fi} by the author, 
where, among other results, 
the $Fi_{22}$-extended polar space was characterized
by a stricter assumption involving certain 3- and 4-vertex configurations in
the point graph. We also prove certain conjectures made in \cite{BH} on the 
extendability of polar spaces.
Our main statement was known to be true under additional assumptions on the 
existence of a flag-transitive automorphism group acting on the geometry, see e.g. \cite{BH,DFGMP,Mei,vBW:fi:char}.

A connected incidence system $\G=\G(\cP,\cB)$ is an {\em extended polar space}
(respectively {\em extended (projective) plane}) if its point residues 
are finite, thick, nondegenerate polar spaces 
(respectively finite nondegenerate projective pla\-nes).

We say that an extended polar space $\G$ admits {\em extended planes} if there
exists a nonempty set $\S=\S(\G)$ of subsets of the set $\cP$ of points of $\G$
such that $\{x\}\cup\pi\in\S$ whenever $\pi$ is a plane in the residue
$\G_x$ of the point $x$ and for any $\Xi\in\S$ the incidence system of points 
and blocks of $\G$ on $\Xi$ is an extended projective plane.

\begin{thm}
Let $\G$ be an extended polar space admitting extended planes. Then exactly one of the following holds.
\begin{itemize}
\item[$(i)$] $\G$ is a standard quotient of an affine polar space over 
{\rm GF(2)}.
\item[$(ii)$] $\G$ is the $Fi_{22}$-extended polar space.
\item[$(iii)$] There are exactly two extended planes on each block, 
and $\G$ is an extended $Q^+_5(4)$-polar space.
\end{itemize}
\label{MTH}
\end{thm}

{\bf Remark.} The geometries in case $(i)$ are all known, see
\cite{BH,Bue:extpolarsp} (also \cite{HaSh,CoSh,CuPa} for various
generalizations). Namely, let $\G$ be an affine polar space over GF(2),
that is the geometry induced on the complement of a geometric hyperplane
in a nondegenerate polar space over GF(2). $\G$ admits a (proper)
standard quotient if and only if the diameter of the point graph on $\G$
is 3. In the latter case the point graph is a two-fold antipodal cover
of a complete graph. Here being at maximal distance is an equivalence 
relation with the classes of size 2,
and the quotient is defined on the equivalence classes of objects 
(that is, points, lines, etc.).

In case $(iii)$, however, no examples are known at all.

{\bf Remark.} It is easy to see that $\G$ has the structure of a geometry 
with the following diagram: (see \cite{Bue:diag} for the notion of diagram 
geometry)
\begin{center}
$$\node\stroke{\subset}\node\arc\node\cdots\node\arc\node\darc\node.$$
\end{center}
The types of the elements are, from left to right, as follows:
points, pairs of adjacent points, blocks, extended projective
2-spaces,\dots, extended projective $n-1$-spaces, where $n$ is the rank of
the polar space we are extending. Note that the elements corresponding
to the extended $j$-subspaces for $j>2$ should be recovered as
certain subsets of points, along the lines of the proof of
Lemma~\ref{res.GF2orGF4}, where extended 3-spaces are reconstructed.

The diagrams corresponding to the cases of the Theorem are as follows:
$$\hfill{\node\stroke{\subset}\node_2\arc\node_2\cdots
\node_2\arc\node_2 \darc\node_t,}\leqno(i)
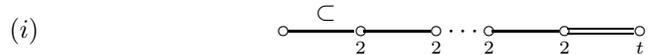
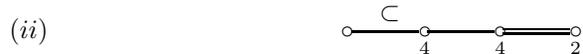
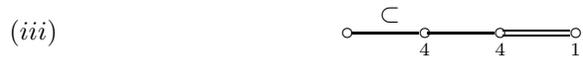
$$
$$\hfill{\node\stroke{\subset}\node_4\arc\node_4\darc\node_2.}\leqno(ii)$$
$$\hfill{\node\stroke{\subset}\node_4\arc\node_4\darc\node_1.}\leqno(iii)$$
%
Note that $t\in\{1,2,4\}$.

{\bf Remark.} Let $\G$ be an extended polar space with point residues of
rank at least 3. If each triple of pairwise adjacent points of $\G$ 
lies in a block (that is, $\G$ is {\em triangular}) then, by \cite{BH}, 
$\G$ admits extended planes.
However, since the proper standard quotients of affine polar spaces are not
triangular, the triangularity assumption is stronger than the assumption on 
existence of extended planes.
\medskip

Since in case $(iii)$ no examples are known, it seems quite natural to
ask the following.

{\bf Question.} Does there exist an extension $\G$ of the $Q_5^+(4)$-polar
space $\D$?

It appears that the technique employed in the paper cannot cope with
this question. Indeed, $\D$ admits hyperovals, cf.
Proposition~\ref{hyp:classes}, and obvious counting tricks do not rule
them out. By comparing the orders of the automorphism groups of the
hyperovals (see Section~\ref{sec5}) 
to the one of $\D$, one sees that the number of hyperovals
in $\D$ is huge. It makes attempts to reconstruct $\G$ using a
computer particularly difficult. In this respect,
it would be quite interesting to find a computer-free construction of
the hyperovals of $\D$. 
\medskip

{\bf Hyperovals and extensions of generalized quadrangles.}
A {\em hyperoval} $\O$ of a polar space $\Pi$ is a set of points of $\Pi$ such that
each line of $\Pi$ meets it in 0 or 2 points. 
(In \cite{BH} such objects are called {\em local subspaces}.)
They play an important role in extensions of $\Pi$.
As it was observed in \cite{BH}, a hyperoval $\O$ in a polar space $\Pi$ of rank
$r$ is a triangular extended polar space of residual rank $r-1$. The
structure of an extended polar space on $\O$ can be seen within the subgraph of
the point graph of $\Pi$ induced by $\O$.

In particular, as happens in the case we shall be mainly concerned with
(that is $r=3$, $\Pi$ is over GF(4)), $\O$ is an extended generalized
quadrangle of order $(4,t)$ (EGQ$(4,t)$, for short), 
where $t\in\{1,2,4,8,16\}$.  
EGQ$(s,t)$ are rather interesting in themselves, in particular a number of
finite simple groups, including sporadic ones, act on them
as flag-transitive automorphism groups, see e.g. \cite{CHP}.
As a by-product of our investigation of hyperovals we find new examples of
EGQ$(4,t)$, $t=1,2$.
However, these examples do not admit flag-transitive automorphism groups. For
$t=2$ our example is the only known example of an extension of a classical
generalized quadrangle of order $(s,t)$, $t>1$, which does not admit a
flag-transitive automorphism group.

These examples are discussed in greater detail in the last section of 
the paper.

Note  that  the  proof of Theorem \ref{MTH} depends upon  calculations
using the computer algebra system {\sf GAP} \cite{GAP} and its package
for computations with graphs and groups GRAPE \cite{GRAPE}.

\section{Preliminaries}
An {\em incidence system} is a pair $\G=\G(\cP,\cB)$, where $\cP$ is a set (of {\em
points}) and
$\cB$ is a set of subsets of $\cP$ (called {\em blocks}), 
each element of $\cB$ is of size greater
than one and the incidence between points and blocks is defined by inclusion.
Two points $p,q$ are said to be {\em adjacent} (notation $p\perp q$)
if they lie in a common block. The set of points adjacent to $p$ is denoted
by $p^\perp$. The {\em point graph} of $\G$ is a graph with vertex set $\cP$,
and adjacency the same as in $\G$.
The {\em residue} $\G_p$ of a point $p$ is the incidence system of points adjacent
to $p$ excluding $p$ itself, and blocks on $p$ with the point $p$ removed.
We say that $\G$ is {\em connected } if its point graph is connected.

From now on $\G$ is an extended polar space of residual rank at least 3
admitting extended planes.

\begin{lem}
The point residues of $\,\G$ are isomorphic.
\label{res.iso}
\end{lem}
\begin{pf}
Let $\Phi$ be an extended plane of $\G$ and $x\in\Phi$. Then $\G_x$ is a
polar space of rank at least $3$. Its isomorphism type is uniquely determined
by the isomorphism type of the polar space of blocks and extended planes on
$xy$, where $y\in x^\perp-\{x\}$, and incidence is by inclusion.
Hence $\G_x\cong\G_y$, and by connectivity of $\G$ the result follows.
\end{pf}

\begin{lem}
Let $\D=\G_x$, $x\in\cP$.
Then $\D$ is either
\begin{itemize}
\item[$(i)$] a polar space over {\rm GF(2)}, or
\item[$(ii)$] a rank $3$ polar space over {\rm GF(4)}.
\end{itemize}
\label{res.GF2orGF4}
\end{lem}
\begin{pf}
Note that $\D$ admits an embedding into a projective space over
GF$(q)$, cf. e.g. Tits~\cite{Ti:bldgs}. In particular the planes of
$\D$ are of a prime power order $q$. 
Since the only projective planes of such an order
which are extendable are of order 2 or 4 (cf. \cite{Hug65.1,Dem}),
$\D$ is defined over GF(2) or GF(4). It remains to check that in the latter
case the rank of $\D$ is 3.

Assume to the contrary that $\D$ is defined over GF(4) and has rank at least
4. Let $X=\{x\}\cup\tau$, where $\tau$ is the set of points of a
3-dimensional singular projective subspace of $\D$.
We claim that $\G(X,\cB_X)$ is an extended PG(3,4), where
$\cB_X=\{B\in\cB\mid B\subset X\}$.
It suffices to show that, given three points $a,b,c\in X$, there exists a
unique $B\in\cB_X$ containing them.
There exists a plane $\pi\subset\tau$ containing $a$,
$b$ and $c$.
Then $\Pi=\{x\}\cup\pi\in S$ is an extended plane, so there exists
$B\in\cB_\Pi$ containing $a$, $b$ and $c$.
The uniqueness of $B$ follows from the fact that $\G_a$ is a polar space and
$B-\{a\}$ is a line of it.
Thus $\G(X,\cB_X)$ is an extended PG(3,4).
On the other hand, according to \cite{Hug65.1} (see also \cite{Dem}), PG(3,4) is
{\em not} extendable. This contradiction implies $(ii)$.
\end{pf}

It was shown in \cite{BH,Bue:extpolarsp}
that in case $(i)$ of Lemma \ref{res.GF2orGF4}, $\G$ is
a standard (perhaps improper) quotient of an affine polar space over GF(2), and
in fact all such objects are known.
Thus we assume from now on that $\D=\G_x$ is a rank 3 polar space over GF(4).

\begin{lem}
Let $C$ be a block of $\G$, $r$ a point of $\G$, and suppose that
there is no extended plane on $r$ and $C$. Then $|C\cap r^\perp|=0,2$ or $4$.
\label{pt.vs.block}
\end{lem}
\begin{pf}
Assume that $C\cap r^\perp$ is not empty. Let $x\in C\cap r^\perp$. Then there
exists a unique point $x'\in C$ such that $\{x'\}=C\cap r^{\perp_{\G_x}}$.
So there exists a block $B^x$ containing $r$, $x$ and $x'$.
Clearly $C\cap B^x=\{x,x'\}$. Since obviously $B^{x'}=B^x$, there is an
equivalence relation on $C\cap r^\perp$ with classes $\{x,x'\},\dots$.
Thus $|C\cap r^\perp|$ is even. It remains to show that $C\not\subset
r^\perp$.

Let $\{p,p',q,q'\}\subseteq C\cap r^\perp$, and let $\Pi$ be an extended plane
containing $C$. Then $\Pi$ contains a unique block $C^q$ such that
$q,q'\in C^q$ and $C^q-\{q\}=\Pi\cap r^{\perp_{\G_q}}$. Also,
$\Pi$ contains a unique block 
$C^p$ such that $p,p'\in C^p$ and $C^p-\{p\}=\Pi\cap r^{\perp_{\G_p}}$.
Note that $|C^p\cap C^q|$ is 0 or 2.
Assume first that the latter holds, that is $\{y,y'\}=C^p\cap C^q$.
Then the three extended planes $\Pi$, $\langle C^p,r\rangle$ and 
$\langle C^q,r\rangle$ form a triangle in the generalized quadrangle 
of blocks and extended planes containing $yy'$, a contradiction.
So $C^p\cap C^q=\emptyset$.

Assume $\{x,x'\}=r^\perp\cap C-\{p,p',q,q'\}$. Then $\Pi$ contains a unique
block $C^x$ such that $x,x'\in C^x$ and $C^x-\{x\}=\Pi\cap r^{\perp_{\G_x}}$.
Repeating the argument above with $C^x$ in place of $C^p$, one has that
$C^x\cap C^q=\emptyset$. Similarly, $C^x\cap C^p=\emptyset$.
So the extended plane $\Pi$ admits three blocks with trivial pairwise
intersections, the situation well-known to be impossible.
This is a contradiction. 
\end{pf}

Note that if $\G$ is triangular then the case $|p^\perp\cap C|=4$ does
not occur.

\begin{pr}
There exists $p\in\cP$ such that $\D=\G_p$ contains a hyperoval $\O$.
Moreover, if $\G$ is not triangular then $\O\cap r^{\perp_\D}=\emptyset$ for
some $r\in p^\perp-\{p\}$.
\label{ho.exist} 
\end{pr}
\begin{pf}
If $\G$ is triangular then the statement is well-known, cf. \cite{BH}.
So we assume that $\G$ is not triangular.
There exist points $p$, $q$, $r$ such that $p\perp q\perp r\perp p$, but
there is no extended plane containing all of them.
Note that this implies $d_{\G_p}(q,r)=d_{\G_q}(p,r)=d_{\G_r}(p,q)=2$.
In what follows we call such a triple $pqr$ of points a {\em bad triangle}.
Let $C$ be a block on $pq$. By Lemma \ref{pt.vs.block}, $C\cap
r^\perp=\{p,p',q,q'\}$, where $p'$, $q'$ are such that there exist blocks
containing $\{p,p',r\}$, respectively $\{q,q',r\}$.
Now we are able to define a hyperoval $\O$ in $\G_p$.
For each block $B$ on $p$ such that $|B\cap r^\perp|=4$, set
$z,z'\in\O$, where $B\cap r^\perp=\{p,p_1,z,z'\}$ and
there exist blocks containing $\{r,p,p_1\}$ and $\{r,z,z'\}$.
Note that the triangles $prz$ and $prz'$ are bad.

The remaining possibilities for $B$ are  $|B\cap r^\perp|=2$ or $6$ 
and we set $B\cap\O=\emptyset$.
We shall check that our definition of $\O$ is correct.
If $|B\cap r^\perp|=2$ or $6$ then all the triangles $prw$, where $w\in
B\cap r^\perp$, are good. So there is no block $B'$ such that $|B'\cap
r^\perp|=4$ and $w\in B'$. We are done.
This proves the first part of the proposition. The second one follows from
the construction of $\O$.
\end{pf}

The following observation will be useful in determining the hyperovals of
$\D$.

\begin{lem}
Let $\Pi$ be a subspace of $\D$. Then $\Pi\cap\O$ is a hyperoval of $\Pi$
whenever $\O$ is a hyperoval of $\D$. \qed
\label{sspace.vs.ho}
\end{lem}

\section{Classification of the hyperovals}
Let $\D$ be a rank 3 nondegenerate polar space over GF(4). Then $\D$ is one
of the following polar spaces: $Q_5^+(4)$, $S_5(4)$, $Q_7^-(4)$, $H_5(4)$ or
$H_6(4)$. We say that two hyperovals of $\D$ are of the same {\em type } if
one can be mapped onto the other by an automorphism of $\D$.
\begin{pr} 
The hyperovals $\O$ in $\D$ are as follows.
\begin{itemize}
\item[$(i)$] $\D\cong Q_5^+(4):$ two types, one has $72$ points, the other
has $96$ points.
\item[$(ii)$] $\D\cong H_5(4):$ two types, one has $126$ points, the other
has $162$ points. Each point of $\D$ is collinear with a point of $\O$.
\item[$(iii)$] $\D\cong S_5(4)$, $Q_7^-(4)$ or $H_6(4)$. No hyperovals.
\end{itemize}
\label{hyp:classes}
\end{pr}
\begin{pf}
$(i)$ and $(ii)$ are results of a computer backtrack search.
Following Lemma \ref{sspace.vs.ho}, we first find all the hyperovals 
(up to type) of $x^\perp\subset\D$ containing $x\in\cP(\D)$. Then
we try to extend each of them to hyperovals of $\D$.

$(iii)$. For $\D\cong S_5(4)$ we use a computer and the fact that
$Q_5^+(4)\subset S_4(4)$, which allow us, by Lemma \ref{sspace.vs.ho}, to start
the search from a hyperoval of $Q_5^+(4)$.
It turns out that none of the hyperovals of $Q_5^+(4)$ are extendable to a
hyperoval in $\D$.

For $\D\cong Q_7^-(4)$, the Proposition follows immediately from Lemma
\ref{sspace.vs.ho} and the fact that $S_5(4)\subset Q_7^-(4)$ does not admit
hyperovals.

For $\D\cong H_6(4)$ we are able to give a proof which is
computer-free, apart from using part $(ii)$ of the Proposition.
Let $\O$ be a hyperoval of $\D$, and let $\cS$ be the set of the
$H_5(4)$-subspaces of $\D$. We claim that $\Phi\cap\O\not=\emptyset$ for any
$\Phi\in\cS$. 

Let $\Phi_0\in\cS$ intersect $\O$ nontrivially (such $\Phi_0$ clearly
exists, since there is an element of $\cS$ on any point of $\D$).
Each $\Phi_1\in\cS$ intersecting $\Phi_0$ in a hyperplane with a deep point
$p$ (that is, in the set $p^{\perp_{\Phi_1}}$, $p\in\cP(\Phi_1)$) satisfies
$\Phi_1\cap\O\not=\emptyset$, as well.
Indeed, by the second part of $(ii)$, 
$p^{\perp_{\Phi_0}}=p^{\perp_{\Phi_1}}\subset\Phi_1$
intersects $\O$ nontrivially, so $\Phi_1$ intersects $\O$ nontrivially.

Finally, note that
the graph defined on $\cS$ such that two vertices are adjacent
iff the corresponding subspaces intersect in a hyperplane with deep point is
connected. In particular, there is a path in this graph from $\Phi_0$ to
$\Phi$. Hence $\Phi\cap\O\not=\emptyset$.

So each of the 2752 elements $\Phi$
of $\cS$ corresponds to a subhyperoval of $\O$, that is,
to a nonempty set $\Phi\cap\O$.
Moreover, different elements of $\cS$ correspond to different subhyperovals.
There are 704 elements of $\cS$ on each point of $\D$. So there is the same
number of subhyperovals on each point of $\O$.
By $(ii)$, each of those subhyperovals has either 126 or 162 points.
Counting in two ways the number of incident point-subhyperoval pairs, we have
\begin{equation}
704|\O|=162a+126(2752-a),
\label{eq:pt-subho}
\end{equation}
where $a$ denotes the number of 162-point subhyperovals of $\O$.

Next, note that there are 176 elements of $\cS$ on any pair of noncollinear
points of $\D$. Hence each pair of noncollinear points of $\O$ is contained
in the same number of subhyperovals.
Counting in two ways the number of the pairs ``pair of noncollinear points
within a subhyperoval'', we have
\begin{equation}
176|\O|(|\O|-166)=162\cdot 116a+126\cdot 80(2752-a).
\label{eq:ptpt-subho}
\end{equation}

Since the system of equations (\ref{eq:pt-subho})--(\ref{eq:ptpt-subho}) does
not have any nonnegative integral solutions, $\O$ does not exist.
\end{pf}

\section{Extensions of $H_5(4)$}
Here we complete the proof of Theorem \ref{MTH}.
It follows from Proposition \ref{hyp:classes} and the following result.
\begin{pr}
Let $\G$ be an extension of $\D\cong H_5(4)$. Then $\G$ is the extended polar
space for $Fi_{22}$.
\label{ext:H54}
\end{pr}
\begin{pf}
It follows from Propositions \ref{hyp:classes} $(ii)$ and \ref{ho.exist} that
$\G$ is triangular. So $\G$ can be recovered from its point graph, cf.
\cite{BH}. For
simplicity, we denote the point graph of $\G$ by $\G$.

For a point $u\in\cP(\G)$, denote by $\G_2(u)$ the set of
points at distance 2 from $u$.
For each $x\in u^\perp-\{u\}$, $|x^\perp\cap\G_2(u)|=512$. 
For each $v\in\G_2(u)$, $|u^\perp\cap v^\perp|=126$ or $162$, cf. Proposition
\ref{hyp:classes} $(ii)$. Counting in two ways the number of edges between
$u^\perp$ and $\G_2(u)$, we see that there exists
$v\in\G_2(u)$ such that $|u^\perp\cap v^\perp|=126$. Indeed,
otherwise, for any $x\in u^\perp-\{u\}$, 
$$|\G_2(u)|=|u^\perp-\{u\}|\cdot |x^\perp\cap\G_2(u)|/162=693\cdot
512/162,$$ 
which is a non-integer.

Let $w\in v^\perp-u^\perp-\{v\}$. The subgraph induced by $\{u,v,w\}^\perp$
is isomorphic to the 2-clique extension $\Phi$ of the collinearity graph of
GQ(2,2), cf. \cite[Lemma~2.8~$(i)$]{Pa:fi}.
Thus the hyperovals $u^\perp\cap v^\perp$ and $u^\perp\cap w^\perp$ must
intersect in a subgraph isomorphic to $\Phi$.
A computer search shows that $|u^\perp\cap w^\perp|\not=162$.
By \cite[Lemma~2.10]{Pa:fi}, $\{u,v,w\}^\perp=\{u,v,w'\}^\perp$ implies
$w=w'$ for $w'\in v^\perp-u^\perp$. 

Let $\Xi$ be the graph defined on the 126-point hyperovals of $\D$, 
such that two vertices $\O$, $\O'$ are adjacent iff $\O\cap\O'\cong\Phi$.
By \cite[Lemma~2.11]{Pa:fi}, $\Xi$ has three connected components $\Xi'$,
$\Xi''$ and $\Xi'''$, each of size 1408 and valence 567. 
They are permuted by automorphisms of $\D$.
Since $|v^\perp-u^\perp-\{v\}|=567$, the connected component $\Lambda$ 
containing $v$ of the subgraph of $\G$ induced on $\G_2(u)$ is a cover of
$\Xi'$. By counting the edges between $u^\perp$ and $\G_2(u)$,
the index of this cover is at most 2.
Since $\mu(\Xi')=216$, a number bigger than the size of a hyperoval of
$\D$, we find that $\Xi'\not\cong\Lambda$. 
So $\Lambda$ is a 2-fold cover of $\Xi'$. 

Continuing
as in \cite[Sect.~2.3]{Pa:fi} (that is by observing that $\G$ is a
triple graph and recovering the 3-transposition group associated with
it), it follows
that $\G$ is the example related to $Fi_{22}$.

Alternatively, note that we have shown that $|\G(x,y)|=126$ for any two
vertices $x$, $y$ at distance 2. Since all the hyperovals of $H_5(4)$ of
size 126 are of the same type (cf. Proposition~\ref{hyp:classes}~$(ii)$), 
they must be isomorphic to the particular type of  the hyperovals
appearing in the $Fi_{22}$-example. Therefore the condition on
$C_4$-subgraphs of $\G$ in \cite[Theorem~1.1]{Pa:fi} holds, and so
$\G$ is indeed the $Fi_{22}$-example.
\end{pf}

\section{New EGQ(4,1) and EGQ(4,2)}
\label{sec5}
Here we describe new extended generalized quadrangles constructed as
hyperovals of a polar space $\D$. A computer-free proof of the existence of
the 162-point EGQ(4,2) is given. We slightly abuse notation by
sometimes identifying
the hyperoval with the respective EGQ and/or with the point graph of this
EGQ.
\paragraph{$\underline{\D=Q_5^+(4)}$.}
There are two types of hyperovals. They give two nonisomorphic  
EGQ(4,1). The first one, with 72 points, is isomorphic to the one in
\cite[Example 9.15 $(iii)$]{CHP}. Its automorphism group is of order 28800
and it acts flag-transitively.
The second one, with 96 points, is apparently new. Its automorphism group 
is of order 3200 and it has two orbits on points of length 16 and 80,
respectively. The distribution diagrams with respect to a point (see
\cite{BCN}) are different for points from different orbits. The diameter of
the point graph is 3.

\paragraph{$\underline{\D=H_5(4)}$.}
There are two types of hyperovals. They give two nonisomorphic  
EGQ(4,2). The first one, with 126 points, is well-known. See e.g.
\cite[Example 9.9 (b) $(ii)$]{CHP}, \cite{Pa:fi}.

A hyperoval $\Pi$ of the second type has 162 points.
The automorphism group of $\Pi$ is isomorphic to 
$(3^4:S_6).2$, and it is not flag-transitive (namely, there are two edge
orbits). We show its existence without a computer.
\begin{lem} {\rm (A.E.~Brouwer and H.~Cuypers)}
There exists a pair of $\,126$-point hyperovals 
$\,\Xi$ and $\,\Xi'$ intersecting
in a $45$-point Fischer subspace of $\,\Xi$.
The symmetric difference $\,\Pi$ of $\,\Xi$ and $\,\Xi'$ is a 
$162$-point hyperoval of $\D$.
\label{ho:162}
\end{lem}
\begin{pf}
The first part of the statement easily follows from 
\cite[Lemma 2.3 (2)]{Pa:fi}.
Let us turn to the second part.

Let $\Xi$ be a
126-point hyperoval of $\D$. It has Fischer subspaces $S$ of
size 45 (in the natural orthogonal description of $\Xi$, 
where the points are (+)-type points of the 6-dimensional GF(3)-space
equipped with a nondegenerate symmetric bilinear form with discriminant 1, 
adjacency coincides with perpendicularity, and
each $S$ consists of the points orthogonal to a given isotropic point). 
Pick one such $S$. There are
exactly three 126-point hyperovals intersecting in $S$, 
see \cite[Lemma~2.3]{Pa:fi}. 
Pick any two of them,
say $\Xi, \Xi'$. We claim that
the subgraph $\Pi$ induced on the symmetric difference of
the pointsets of $\Xi$ and $\Xi'$ is a hyperoval.

Let $x\in\Xi-S$, and let $l$ be a line of $\D$ on $x$. Let $\{x,y\}=l\cap\Xi$. 
Assume first that $y\in S$. Then $l\cap\Xi'=\{x',y\}$, so $l\cap\Pi=\{x,x'\}$,
as required. 

It remains to show that in the case $y\in\Xi-S$ we have $l\cap\Xi'=\emptyset$.
By \cite[Lemma~2.1~$(iv)$]{Pa:3tr},  the subgraph of the point graph of $\D$
induced on $x^\perp\cap S$ is isomorphic to the point graph $\Upsilon$ of GQ(2,2).
By \cite[Lemma~2.8~(1)]{Pa:fi},   the subgraph of the point graph of $\D$
induced by $x^\perp\cap\Xi'$ is isomorphic to the 2-clique extension of
$\Upsilon$. Thus if a line on $x$ intersects $\Xi'$ then it intersects $S$,
and we are in the case already considered.
\end{pf}

Note that the 126-point hyperovals and the sets $S$ described above
are Fischer subspaces of $\D$, and they are related to near subhexagons of
$H_5(4)$-dual polar space, see \cite{BCHW}. The proof of the lemma given
by Andries Brouwer used a technique developed in \cite{BCHW}.
The proof just given seems to be a streamlined version of Brouwer's proof.
\medskip

Next, we give two other descriptions of $\Pi$.
The first one is the author's description communicated to Andries Brouwer.
Let $\Lambda$ be the incidence graph of a $PG(6,6,2)$ defined on
\cite[p.373]{BCN} (Van Lint-Schrijver $PG$). 
Note that $\Lambda$ is not distance
transitive. Let $u$ be a vertex of $\Lambda$.
$\Lambda_3(u)$ is split into the two orbits $\Lambda_3^1(u)$
and $\Lambda_3^2(u)$ of lengths 60 and 15 respectively.   
$\Lambda_4(u)$ is split into the two orbits $\Lambda_4^1(u)$
and $\Lambda_4^2(u)$ of lengths 20 and 30 respectively.   
To obtain the point graph of $\Pi$, one has to choose as the set 
$\Pi(u)$ the union of $\Lambda_3^2(u)$ and $\Lambda_4^2(u)$, and then
apply $Aut(\Lambda)$ to get all the remaining edges.

Andries Brouwer (personal communication) also
gave the following elegant description of $\Pi$. 
Let $$f(x,y)=x_1y_1^2+\dots+x_6y_6^2$$ be the Hermitian form defining
$\D$. Then, without loss of generality,
the points of $\Pi$   are the points of $\D$ such that the product
of the coordinates lies in ${\rm GF(4)}\setminus {\rm GF(2)}$.

\paragraph{Remark.} Steven Linton (personal communication), computed,
using his vector enumeration program, which implements his algorithms
presented in \cite{Lin2}, that the fundamental group of the
162-point EGQ(4,2) is perfect. Hence the index of any its
proper covers is at least 60, the order of the smallest perfect group.
On the other hand, upper bounds on the number of points in EGQ$(s,t)$
given in \cite{CHP} show that such big  covers do not exist.
Thus the 162-point EGQ(4,2) does not have any proper covers.

Note that the 126-point EGQ(4,2) has the triple
cover, which is the only proper cover, see \cite{BBegq42,Mei,CHP}.

\begin{figure}
\setlength{\unitlength}{0.0125in}%
\begin{picture}(327,229)(59,570)
\thicklines
\put( 80,680){\circle{42}}
\put(160,740){\circle{42}}
\put(160,640){\circle{42}}
\put(265,765){\circle{42}}
\put(265,690){\circle{42}}
\put(365,690){\circle{42}}
\put(265,605){\circle{42}}
\put(100,690){\line( 5, 4){ 44.512}}
\put(100,675){\line( 2,-1){ 40}}
\put(160,720){\line( 0,-1){ 55}}
\put(180,745){\line( 3, 1){ 64.500}}
\put(180,730){\line( 2,-1){ 60}}
\put(175,725){\line( 3,-4){ 77.400}}
\put(180,650){\line( 2, 1){ 64}}
\put(180,635){\line( 5,-2){ 64.655}}
\put(265,745){\line( 0,-1){ 30}}
\put(265,670){\line( 0,-1){ 40}}
\put(285,760){\line( 6,-5){ 65.410}}
\put(285,690){\line( 1, 0){ 55}}
\put(280,620){\line( 4, 3){ 68.800}}
\multiput(280,750)(0.40000,-0.40000){26}{\makebox(0.4444,0.6667){\sevrm .}}
\put(290,740){\line( 1,-2){ 10}}
\put(300,720){\line( 1,-3){  5}}
\put(305,705){\line( 0,-1){ 25}}
\put(305,680){\line(-1,-4){  5}}
\multiput(300,660)(-0.25000,-0.50000){21}{\makebox(0.4444,0.6667){\sevrm .}}
\put(295,650){\line(-2,-3){ 10}}
\multiput(285,635)(-0.40000,-0.40000){26}{\makebox(0.4444,0.6667){\sevrm .}}
\put(155,735){\makebox(0,0)[lb]{\raisebox{0pt}[0pt][0pt]{\twlrm 30}}}
\put(155,635){\makebox(0,0)[lb]{\raisebox{0pt}[0pt][0pt]{\twlrm 15}}}
\put(100,705){\makebox(0,0)[lb]{\raisebox{0pt}[0pt][0pt]{\twlrm 30}}}
\put(130,725){\makebox(0,0)[lb]{\raisebox{0pt}[0pt][0pt]{\twlrm 1}}}
\put(100,655){\makebox(0,0)[lb]{\raisebox{0pt}[0pt][0pt]{\twlrm 15}}}
\put(130,645){\makebox(0,0)[lb]{\raisebox{0pt}[0pt][0pt]{\twlrm 1}}}
\put(150,705){\makebox(0,0)[lb]{\raisebox{0pt}[0pt][0pt]{\twlrm 3}}}
\put(150,665){\makebox(0,0)[lb]{\raisebox{0pt}[0pt][0pt]{\twlrm 6}}}
\put(155,605){\makebox(0,0)[lb]{\raisebox{0pt}[0pt][0pt]{\twlrm 6}}}
\put(160,765){\makebox(0,0)[lb]{\raisebox{0pt}[0pt][0pt]{\twlrm 9}}}
\put(260,760){\makebox(0,0)[lb]{\raisebox{0pt}[0pt][0pt]{\twlrm 20}}}
\put(260,685){\makebox(0,0)[lb]{\raisebox{0pt}[0pt][0pt]{\twlrm 60}}}
\put(260,600){\makebox(0,0)[lb]{\raisebox{0pt}[0pt][0pt]{\twlrm 30}}}
\put(360,685){\makebox(0,0)[lb]{\raisebox{0pt}[0pt][0pt]{\twlrm 6}}}
\put(345,715){\makebox(0,0)[lb]{\raisebox{0pt}[0pt][0pt]{\twlrm 10}}}
\put(180,750){\makebox(0,0)[lb]{\raisebox{0pt}[0pt][0pt]{\twlrm 8}}}
\put(185,730){\makebox(0,0)[lb]{\raisebox{0pt}[0pt][0pt]{\twlrm 12}}}
\put(180,660){\makebox(0,0)[lb]{\raisebox{0pt}[0pt][0pt]{\twlrm 24}}}
\put(185,635){\makebox(0,0)[lb]{\raisebox{0pt}[0pt][0pt]{\twlrm 8}}}
\put(235,705){\makebox(0,0)[lb]{\raisebox{0pt}[0pt][0pt]{\twlrm 6}}}
\put(235,680){\makebox(0,0)[lb]{\raisebox{0pt}[0pt][0pt]{\twlrm 6}}}
\put(230,770){\makebox(0,0)[lb]{\raisebox{0pt}[0pt][0pt]{\twlrm 12}}}
\put(235,600){\makebox(0,0)[lb]{\raisebox{0pt}[0pt][0pt]{\twlrm 4}}}
\put(230,625){\makebox(0,0)[lb]{\raisebox{0pt}[0pt][0pt]{\twlrm 12}}}
\put(250,630){\makebox(0,0)[lb]{\raisebox{0pt}[0pt][0pt]{\twlrm 10}}}
\put(255,655){\makebox(0,0)[lb]{\raisebox{0pt}[0pt][0pt]{\twlrm 5}}}
\put(280,710){\makebox(0,0)[lb]{\raisebox{0pt}[0pt][0pt]{\twlrm 21}}}
\put(275,635){\makebox(0,0)[lb]{\raisebox{0pt}[0pt][0pt]{\twlrm 6}}}
\put(260,570){\makebox(0,0)[lb]{\raisebox{0pt}[0pt][0pt]{\twlrm 12}}}
\put(290,610){\makebox(0,0)[lb]{\raisebox{0pt}[0pt][0pt]{\twlrm 1}}}
\put(260,790){\makebox(0,0)[lb]{\raisebox{0pt}[0pt][0pt]{\twlrm 9}}}
\put(295,755){\makebox(0,0)[lb]{\raisebox{0pt}[0pt][0pt]{\twlrm 3}}}
\put(280,735){\makebox(0,0)[lb]{\raisebox{0pt}[0pt][0pt]{\twlrm 9}}}
\put(250,730){\makebox(0,0)[lb]{\raisebox{0pt}[0pt][0pt]{\twlrm 12}}}
\put(255,715){\makebox(0,0)[lb]{\raisebox{0pt}[0pt][0pt]{\twlrm 4}}}
\put(165,705){\makebox(0,0)[lb]{\raisebox{0pt}[0pt][0pt]{\twlrm 12}}}
\put(330,695){\makebox(0,0)[lb]{\raisebox{0pt}[0pt][0pt]{\twlrm 30}}}
\put(290,675){\makebox(0,0)[lb]{\raisebox{0pt}[0pt][0pt]{\twlrm 3}}}
\put(335,670){\makebox(0,0)[lb]{\raisebox{0pt}[0pt][0pt]{\twlrm 5}}}
\put( 80,675){\makebox(0,0)[lb]{\raisebox{0pt}[0pt][0pt]{\twlrm 1}}}
\end{picture}
\caption{The distribution diagram of the 162-point EGQ(4,2).}
\end{figure}

\medskip
{\bf Acknowledgment.} The author thanks Andries Brouwer for valuable
comments and suggestions. 
\nocite{Ivnv:fi}

\chapter{Multiple extensions of generalized hexagons related to the simple
groups \McL\ and $\Co_3$}
\label{chap:hex}
\newcommand{\ZZ}{{\rm Z}\kern-3.8pt {\rm Z} \kern2pt}
\newcommand{\RR}{${\rm I\kern-1.6pt {\rm R}}$}
\newcommand{\NN}{{\rm I\kern-1.6pt {\rm N}}}
\renewcommand{\HH}{{\rm I\kern-1.6pt {\rm H}}}
\newcommand{\QQ}{${\rm Q}\kern-3.9pt {\rm \vrule height7.1pt
    width.3pt depth-2pt} \kern7.5pt$}
\newcommand{\CC}{${\rm \kern 3.5pt \vrule height 6.5pt
    width.3pt depth-1.3pt \kern-4pt C \kern.8pt}$}
\newcommand{\FF}{${\rm I\kern-1.6pt {\rm F}}$}
\newtheorem{lemma}[thm]{Lemma}
\newtheorem{prop}[thm]{Proposition}
\renewcommand{\pf}{\noindent {\em Proof.} \ }
\newcommand{\eop}{$_{\Box}$ \vspace{5mm} \relax}
\newtheorem{num}[thm]{ }
\newcommand{\split}{\colon\!}
\newcommand{\nsplit}{\cdot}

\begin{center}
Hans Cuypers,
Anna Kasikova, 
Dmitrii V. Pasechnik
\end{center}
\paragraph{Abstract.}
The groups $\McL$, $\Co_3$ and $2\times \Co_3$ are all contained
in the automorphism group of a $2$-, respectively, $3$-fold
extension of a generalized hexagon of order $(4,1)$.
We give a geometric characterization of these multiple extensions of
this generalized hexagon. 

\section{Introduction}
Several of the finite sporadic simple groups act as automorphism group
on extensions of generalized polygons. The large Mathieu groups act on
extensions of the projective plane of order $4$, the groups $\McL$,
$\Co_3$ and $\HS$ are automorphism group of extended generalized
quadrangles, the group $\Ru$ and $\He$ are automorphism groups of
extended octagons, and finally the Hall-Janko group $\HJ$ and the
Suzuki group $\Suz$, are known to act on extensions of generalized
hexagons of order $(2,2)$, respectively, $(4,4)$.  In general,
extensions of generalized hexagons and octagons have infinite covers,
see for example \cite{Pa:covers}.  Thus it seems natural to impose
extra conditions on extensions $\Gamma$ of generalized hexagons or
octagons to obtain characterizations of these geometries.  The two
extended generalized hexagons mentioned above satisfy the following
condition:

\begin{enumerate}
\item[$(*)$] $\{x_1,x_2,x_3\}$ is a clique  of the point graph not contained
in a circle, if and only if $x_2$ and $x_3$ are at distance 3 in the local
generalized hexagon $\Gamma_{x_1}$.
\end{enumerate}

It is just this condition that has been imposed on extensions of generalized 
hexagons in \cite[Theorem~1.1]{Cuy:suz}.
There the above mentioned extended generalized 
hexagons related to the sporadic groups
$\HJ$ and $\Suz$, as well as two other extensions of a generalized hexagon
of order $(2,1)$, respectively, $(4,1)$, related to
the groups $G_2(2)$, respectively, $PSU_4(3)$ were characterized.
In particular, it is shown \cite{Cuy:suz}, 
that these four extended generalized hexagons
are the only extensions of finite, regular, line thick, 
generalized hexagons satisfying $(*)$.

In this paper we are concerned with  (possibly
multiple) extensions of the
four extended generalized hexagons related to the groups $G_2(2)$, $\HJ$, $PSU_4(3)$ and $\Suz$
mentioned above. 
The interest in such geometries comes from the existence
of two examples related to the sporadic simple groups
$\McL$, respectively, $\Co_3$,  that are 1-, respectively,
2-fold extensions of the 
extended generalized hexagon 
on 162 points related to $PSU_4(3)$. Their point graphs are the complement
of the McLaughlin graph on 275 points, and with automorphism
group  $\McL\split2$, and a graph on 552 points
which is locally the complement of the McLaughlin graph and has as 
automorphism group the group $2\times \Co_3$.
This last graph is  one of the two  Taylor graphs
for the third Conway group, it is the one with intersection array
$\{275,112,1;1,112,275\}$, see \cite{BCN}.
These graphs contain a unique class of 7-, respectively, 8-cliques inducing
the multiple extended hexagon, 
which  will be called the {\em multiple extended hexagon related
to the group } $\McL$, respectively, $2\times \Co_3$.
These geometries $\Gamma$ are $r$-fold extensions of generalized hexagons
(with $r$ being $2$, respectively, $3$) satisfying
the following generalization of $(*)$:

\begin{enumerate}
\item[$(**)$] any $r+1$-clique is in a circle
and $\{x_1, \dots, x_{r+2}\}$ is a clique not in a circle,
if and only if $x_{r+1}$ and $x_{r+2}$ are at distance 3 in the 
generalized hexagon $\Gamma_{x_1,\dots, x_r}$.
\end{enumerate}

The second example related to $2\times \Co_3$ admits a quotient geometry on 276 points, which is a 
one point extension of the $\McL$ multiple extended hexagon, 
and has automorphism group $\Co_3$.
Obviously these multiple extended hexagons also 
carry the structure of a Buekenhout geometry
with diagram

\begin{picture}(500,100)
\put(10,70){\circle*{4}}
\put(50,70){\circle*{4}}
\put(90,70){\circle*{4}}
\put(130,70){\circle*{4}}
\put(10,70){\line(1,0){120}}
\put(90,72){\line(1,0){40}}
\put(90,68){\line(1,0){40}}
\put(28,72){$c$}
\put(68,72){$c$}
\put(150,70){for the group $\McL$, respectively,}

\put(10,20){\circle*{4}}
\put(50,20){\circle*{4}}
\put(90,20){\circle*{4}}
\put(130,20){\circle*{4}}
\put(170,20){\circle*{4}}
\put(10,20){\line(1,0){160}}
\put(130,22){\line(1,0){40}}
\put(130,18){\line(1,0){40}}
\put(28,22){$c$}
\put(68,22){$c$}
\put(108,22){$c$}
\put(190,20){for the groups $2\times \Co_3$ and $\Co_3$.}
\end{picture}

\noindent
The groups $\McL\split 2$, $\Co_3$ and $2\times \Co_3$, respectively,
act flag-transitively on these Buekenhout geometries.
As such these groups and geometries were characterized by Weiss in
\cite{We:MCL}, 
under some extra geometric condition closely related to $(*)$ and $(**)$.

In this paper we show that the  two multiple extensions
of a generalized hexagon related to the groups $\McL$ and $2\times \Co_3$
are characterized by the condition $(**)$.  We prove:

\begin{thm}
Let $\Gamma$ be an $r$-fold  extension of 
a finite line thick generalized hexagon satisfying $(**)$ with $r\geq 2$.
Then $\Gamma$ is isomorphic to the 
$\McL$ or $2\times \Co_3$  $2$-, respectively, $3$-fold 
extended generalized hexagon.
\end{thm}

Together with the results of \cite{Cuy:suz}, this supplies us with a complete
classification of multiple  extensions of finite line thick generalized
hexagons satisfying $(**)$.

In fact, the condition $(**)$ is equivalent to saying that
$\Gamma$ is a {\em  triangular} (multiple) extension of one of the four
extended hexagons satisfying $(*)$.
The $3$-fold extended hexagon on 276 points related to the group
$\Co_3$ is not triangular. As such it is characterized in the following
result.

\begin{thm}
  Let $\Gamma$ be a nontriangular $2$-or $3$-fold extension of a
  finite, regular, line thick generalized hexagon in which the residue
  of any point, respectively, pair of cocircular points satisfies
  $(*)$. Then $\Gamma$ is isomorphic to the $\Co_3$ $3$-fold extended
  generalized hexagon on $276$ points.
\end{thm}

We notice that the McLaughlin graph and the Taylor graph related to 
$2\times \Co_3$ with
intersection array $\{275,162,1;1,162,275\}$
are extensions of the unique generalized quadrangle of order $(3,9)$.
These geometries are characterized 
in \cite{Pa:mclco3}.

The point graph of the 3-fold extended generalized hexagon related to
$\Co_3$ is a regular 2-graph, as such it is characterized by 
Goethals and Seidel \cite{GoeSe}.

\section{Definitions and notation}

In this paper we use  notation and definitions of \cite{CHP,Cuy:suz}.
For convenience of the reader we recall some of these definitions
and fix the notation.

Let $\Gamma$ be an incidence structure $({\cal P},{\cal C})$
consisting of nonempty sets $\cal P$ of {\em points} and $\cal C$ of
{\em circles}, where a circle is a subset of $\cal P$ of size at least 2.
Then the {\em point graph} of $\Gamma$ is the graph with vertex set
$\cal P$ and as edges the pairs of distinct points that are {\em cocircular}
(i.e. are in some circle).
That two points $p$ and $q$ are cocircular is denoted by $p\perp q$.
For each subset $X$ of $\cal P$, the set $X^\perp$ consists of
all points $q$ cocircular with all the points of $X$.
We usually write $p^\perp$ instead of $\{p\}^\perp$ for a point $p$ in
$\cal P$.

Fix a point $p$ of $\Gamma$.
If all circles on $p$ 
contain at least 3 points, then  the residue
$\Gamma_p:=({\cal P}_p, {\cal C}_p)$
of $\Gamma$ at $p$, where ${\cal P}_p$ consists of all points of $\cal P$ distinct but cocircular
with $p$ in $\Gamma$, and ${\cal C}_p=\{C-\{p\} \ |\ C\in {\cal C}, p\in C\}$
is also an incidence structure.
For any graph $\cal G$ and point $p$ of the graph, the induced
subgraph on the neighbours of $p$ will be denoted by ${\cal G}_p$.

If $\cal F$, respectively, $\Delta$, is a family of incidence structures, 
respectively, just an incidence structure, then $\Gamma$ is called
a (1-fold) {\em extension of}  $\cal F$, 
respectively, $\Delta$, if and only if all circles of $\Gamma$ have at least
$3$ points, its point graph is connected 
and  $\Gamma_p$ is in $\cal F$, respectively,
is isomorphic to $\Delta$, for all points $p$ of $\Gamma$. 
Inductively, for each integer $r>1$, 
we define $\Gamma$ to be an {\em $r$-fold extension of} $\cal F$ 
or $\Delta$, if and only if all circles contain at least $r+2$ points,
its point graph is connected, and
$\Gamma_p$ is an $r-1$-fold extension of $\cal F$, respectively, $\Delta$
for all points $p$
of $\Gamma$. 
An {\em extended (generalized) hexagon} 
is an extension of a generalized hexagon.

In a generalized hexagon circles are usually referred to as {\em lines}.

An extension is called {\em triangular} if and only if any triangle
of the point graph is contained in a circle.

\section{The 2-fold extended hexagon related to $\McL$}

Let $\Gamma=(\cal P, \cal C)$ be a $2$-fold extension of a 
generalized hexagon satisfying
the hypothesis of Theorem 1.1.
For all points $p$ of $\Gamma$, the residue $\Gamma_p$ is an extended
generalized hexagon satisfying $(*)$. In particular, the results
of \cite{Cuy:suz} apply and we find that $\Gamma_p$ is isomorphic to one of the
four extended generalized hexagons related to $G_2(2)$, $\HJ$, $PSU_4(3)$
or $\Suz$. 

\begin{lemma}
For all points $p$ and $q$ of $\Gamma$ we have $\Gamma_p\simeq\Gamma_q$.
\end{lemma}

\pf Suppose $p$ and $q$ are adjacent points of the point graph of $\Gamma$.
Then $(\Gamma_p)_q\simeq (\Gamma_q)_p$. But, by the results of
\cite{Cuy:suz}, that implies that
$\Gamma_p\simeq \Gamma_q$. Since the point graph of $\Gamma$ is connected,
we have proved the lemma. \eop

Fix a point $p$ of $\Gamma$. By $\Delta$ we denote the residue $\Gamma_p$ 
at $p$.
The point graph of $\Gamma$ will be denoted by ${\cal G}$,
that of $\Delta$ by $\cal D$.
The above lemma implies that $\cal G$ is locally isomorphic to $\cal D$.
Suppose $(p,q,r)$ is a path of length $2$ in the point graph
$\cal G$ of $\Gamma$. Then let $H=H_q$ be the 
complement of $p^\perp \cap q^\perp
\cap r^\perp$ in $p^\perp\cap q^\perp -\{p,q\}$.
The set $H$ is a subset of the point set of the generalized hexagon
$\Delta_q=\Gamma_{p,q}$.
It plays an important role in \cite{Cuy:suz}, and will also be very useful in
the situation considered here.
The set $H$ is the complement of a $\mu$-graph in the point graph of $\Gamma_q$.
The following properties of $H$ are useful in the sequel of this section.

\begin{lemma}
Suppose $\Gamma$ is an extension of a generalized hexagon of order $(s,t)$.
Then we have the following:
\begin{enumerate}
\item[(i)] Each point of $\Delta_q$ not in $H$ is on a unique line
of the generalized hexagon $\Delta_q$ meeting $H$
in $s-1$ points.
\item[(ii)] $H$ contains $(s^2-1)(t^2+t+1)$ points.
\item[(iii)] Each line of $\Delta_q$ is in $H$, or disjoint from $H$,
or meets $H$ in $s-1$ points.
\end{enumerate}
\end{lemma}

\pf Inside $\Gamma_q$ which is one of the four extended generalized hexagons
the points $p$ and $r$ are at distance 2. Thus we can apply the results of 
\cite[Section~3]{Cuy:suz}. In particular Lemma 3.4. \eop

The possible embeddings of  sets $H$ in  $\Delta_q$ or 
${\cal D}_q$ having the properties $(i)-(iii)$ of Lemma 3.3 are
classified in 
\cite{Cuy:suz}.
To state that result we first have to fix some notation.
Suppose $\Delta$ is the $PSU_4(3)$ extended hexagon.
Then the generalized hexagon $\Delta_q$ of order $(4,1)$ will be identified with
the generalized hexagon on the flags of $PG(2,4)$.
A {\em hyperoval} of $PG(2,4)$ is a set of 6 points no three collinear,
its dual is the set of 6 lines missing all the 6 points.
This dual hyperoval is a hyperoval  of the dual plane, i.e.,
no three of the lines meet in a point.

\begin{lemma} Let $H$ be a set of points in $\Delta_q$ having the properties
$(i)-(iii)$ of Lemma $3.3$.
\begin{enumerate}
\item[(i)] If $\Delta_{q}$ has order $(2,1)$,  $(2,2)$ or $(4,4)$, 
then $H$ is the complement of a $\mu$-graph
$q^\perp \cap x^\perp$ for some point $x$ at distance two from $q$ 
in the point graph $\cal D$ of $\Delta$. 
\item[(ii)] If $\Delta_{q}$ has order $(4,1)$, then
$H$ consists of the flags of $PG(2,4)$ missing a hyperoval and its
dual hyperoval.
\end{enumerate}
\end{lemma}

\pf 
This is a straightforward consequence of the results of
\cite[Section~4]{Cuy:suz}. 
\eop

In the graph $\cal G$ the $\mu$-graph $M=p^\perp\cap r^\perp$ is a subgraph
of the graph $\cal D$, which is locally isomorphic to the complement of $H$ inside
the 
point graph of the local generalized hexagon of $\Delta$.

\begin{prop}
$\Delta$ is isomorphic to the $PSU_4(3)$ extended hexagon.
\end{prop}

\pf
Suppose $\Delta$ is not the $PSU_4(3)$ extended hexagon.
Consider the subgraph $M$  of $\cal D$. For
each point $q$ of $M$ the vertex set of the local graph $M_q$ is the complement of
a set $H$ in $\Delta_q$ satisfying $(i)-(iii)$ of Lemma 3.2.

Fix a point $q$ of $M$ and inside $M_q$ a point $x$.
By Lemma 3.3 there is a point $y\in \Delta_q$ at distance 2 from $x$
such that
$M_x$ is the induced subgraph on the set of common neighbours
of $x$ and $y$ inside $\cal D$. 

Inside the generalized hexagon $\Delta_q$, we find that the points $x$ and $y$
have mutual distance 2, so that there is a unique point $z$
collinear
with both $x$ and $y$.
By the choice of $y$, the points $x$ and $z$ are the unique points on
the line through $x$ and $z$ inside $M_q$.
Thus, by 3.2$(i)$, 
this line is the unique line on $z$ containing some point outside
$M$.
In particular, $y\in M_q$.

Repeating the above argument with the role of $x$ and $y$ interchanged,
we see that the line of $\Delta_q$ on $y$ and $z$ is the unique line
on $z$ containing a point outside $M$. A contradiction, and the 
proposition is proved. \eop

Thus from now on  we can assume that $\Gamma$ is locally the extended
generalized hexagon on 162 points related to the group $PSU_4(3)$.
The generalized hexagon of order $(4,1)$ is isomorphic
to the generalized hexagon on the 105 flags of the projective
plane $PG(2,4)$. We will identify $\Gamma_{p,q}=\Delta_q$
with this hexagon. By Lemma 3.3 we know that the set $H$
consists of all the flags missing an hyperoval and its dual in $PG(2,4)$.
The following  can be obtained easily,
especially using the information from
\cite{Atl}. 

The projective plane $PG(2,4)$ admits 168 hyperovals,
the group $P\Gamma L_3(4)$ being transitive on them.
Thus, there are 168 subgraphs in ${\cal D}_q$ isomorphic to $H$.
Of these subgraphs 56 are the complement of a $\mu$-graph of $\cal D$.
The group $PSL_3(4)$ has 3 orbits
on the hyperovals, all of length 56. Two hyperovals are in the same orbit 
if and only if they meet in an even number of points. 
The stabilizer of a hyperoval $O$ in $PSL_3(4)$ is isomorphic to $A_6$.
It has two more orbits on the $PSL_3(4)$-orbit of $O$,
one of length 45, consisting of all hyperovals meeting $O$ in $2$ points,
and one of length 10 consisting of the hyperovals disjoint from $O$.
This stabilizer has two orbits of length 36 and 20 on the two
other $PSL_3(4)$-orbits on the hyperovals, consisting of all hyperovals 
of that orbit meeting $O$ in, respectively, 1 or 3 points of $PG(2,4)$. 
The stabilizer of one of the orbits
in $P\Gamma L_3(4)$ is isomorphic to $P\Sigma L_3(4)$, and permutes the 2
other orbits.

We will use the following description of $\cal D$.
(See \cite{Cuy:suz}.) First we fix the point $q$. The neighbours of $q$ are the
105 flags of $PG(2,4)$, two flags adjacent if and only if
they are at distance 1 or 3 in the generalized hexagon on these flags.
Now fix one of the 3  orbits
of of $PSL_3(4)$ of length 56 on the hyperovals, say $\cal O$.
The vertices at distance 2 from $q$ are the elements of $\cal O$.
Two hyperovals are adjacent if and only if they meet in 2 points.
A flag and a hyperoval are adjacent if and only if the point
(respectively, line) of the flag
lies in the hyperoval (respectively, dual hyperoval).
Fix a hyperoval $O$ not in $\cal O$, and consider the subgraph
$\cal M$ of $\cal D$ consisting of $q$, all flags on $O$ or its dual
and all elements of
$\cal O$ that meet $O$ in 3 points of the plane.
Then $\cal M$ consists of $1+60+20=81$ points.
This graph is the $\mu$-graph appearing in the complement
of the McLaughlin graph and thus locally the complement of $H$ in the 
in the distance 1--or--3 graph of the generalized hexagon of order $(4,1)$.
It is a strongly regular graph with parameters $(v,k,\lambda,\mu)$
equal to $(81,60,45,42)$.
The complement of $\cal M$ in $\cal D$ is also isomorphic to $\cal M$.
There are 112 elements in the $PSU_4(3)$ orbit of $\cal M$,
the complement of $\cal M$ not being one of them.
Any subgraph of this orbit meets $\cal M$ in
81, 45 or 27 points, and hence the complement of $\cal M$
in 0, 36 or 54 points. (This can be checked within the McLaughlin graph.)

We will show that in fact all
subgraphs of $\cal D$ that are locally the complement of $H$ in
the distance 1--or--3 graph of the generalized hexagon of order $(4,1)$
are in the $PSU_4(3)$ orbit of $\cal M$ or its complement in $\cal D$.

\begin{lemma}
The graph $M$ is in the $PSU_4(3)$ orbit of $\cal M$ or of its complement
in $\cal D$.
\end{lemma}

\pf Fix the point $q$, and consider $M_q$, and its complement in ${\cal D}_q$
the set $H_q$.
The graph $M_q$ is  isomorphic to the graph whose vertex set consists
of the 60 flags of $PG(2,4)$ meeting an hyperoval $O_q$ or its dual. 
Two such flags
are adjacent
if and only if they are at distance 1 or 3 in the generalized hexagon
of order $(4,1)$ defined on the flags of $PG(2,4)$.
In particular, one can easily check that
inside $M_q$ two nonadjacent vertices have at least
32 common neighbours.
Thus inside $M$ the $\mu$-graphs consist of at least 33 vertices.
By the above, $M_q$ is either a $\mu$-graph in $\cal D$, or it consists
of the 60 flags of $PG(2,4)$ meeting one of the 112 hyperovals
not in $\cal O$.

The same arguments as used in the proof of Proposition 3.4 rule out
the case where $M_q$ 
is the $\mu$-graph for some point $r$ in $\cal D$
at distance 2 from $q$, i.e., $M_q=r^\perp\cap q^\perp$ in $\cal D$.

Thus assume that $M_q$ is not the $\mu$-graph of some point
$r$ at distance 2 from $q$, i.e. $O_q\not \in \cal O$. 
There are 36  hyperovals in $\cal O$
meeting the hyperoval $O_q$ 
in a unique point of $PG(2,4)$.
Thus each of these 36 points of $\cal D$ is adjacent to 30 points in $H_q$ and 30 points in
$M_q$.
Since $\mu$-graphs in $M$ contain at least 
32 points, none of the 36 points is in $M$.

The remaining 20 points at distance 2 from $q$ 
are hyperovals of $\cal O$ meeting the hyperoval
$O_q$ in 3 point of $PG(2,4)$. 
Hence these 20 points are adjacent to
18 points of $H_q$ and $42$ points an $M_q$.
So $\mu$-graphs of $M$ contain 42 points, and
as $M_q$ has valency 45, there are at least
$60.(60-45-1)/42=20$ points at distance two from $q$. 
Thus,  all 20 points in $\cal O$ meeting $O_q$ in 3 points
of $PG(2,4)$ are in $M$.
In particular, $M$ consists of $q$, $M_q$ and the 20 
points in $\cal O$ meeting $O_q$ in 3 points of the projective plane $PG(2,4)$. 

Since there are 112  hyperovals not
in $\cal O$, there are 
$162.112/81=224$  subgraphs in $\cal D$ that are locally isomorphic
to the complement of $H$ in the distance 1--or--3
graph of the generalized hexagon of order $(4,1)$.
The lemma follows now easily. \eop

\begin{prop}
If $\Delta$ is the extended generalized hexagon on $162$ points related
to $PSU_4(3)$, then  $\Gamma$ is isomorphic to the $2$-fold extended 
hexagon related to $\McL$.
\end{prop}

\pf By the above lemma we find that for any two points at distance
2 in $\cal G$ the number of common neighbours is 81.
Moreover, since each point of the complement of $\cal M$ in $\cal D$
is adjacent to some point in $\cal M$, the diameter of $\cal G$
is 2. Hence $\cal G$ is a strongly regular graph
with parameters $(v,k,\lambda,\mu)=(275,162,105,81)$.

A strongly regular graph with these parameters has been shown to
be isomorphic to the complement of the McLaughlin graph by Cameron, Goethals and Seidel \cite{CGS}.
Here however, we quickly obtain uniqueness of $\Gamma$ in the following way, without using \cite{CGS}.

There are 112 vertices at distance 2 from $q$.
We can identify each of these points with a
subgraph of $\cal D$ isomorphic to $\cal M$.
Without loss of generality we may identify one of these points with
$\cal M$.

The subgraph of $\cal G$ induced on the 112 points not adjacent to $q$
is locally isomorphic with the complement of $\cal M$ in $\cal D$,
and thus with $\cal M$. Since its valency is larger than $112/2$, this subgraph
is connected.

Two  adjacent vertices 
$x$ and $y$ in this graph have 60 common neighbours, and thus
45 common neighbours in ${\cal G}_q$.
But that implies that the two subgraphs $q^\perp \cap x^\perp$
and $q^\perp\cap y^\perp$ are in the same $PSU_4(3)$ orbit.
In particular, by connectivity of the subgraph, the points at distance 2 from $q$ are in one orbit
under the action of $PSU_4(3)$.
But then it is easy to see that the graph $\cal G$ is unique
up to isomorphism. Moreover, $\Gamma$ is unique: as there is only
one set of 5-cliques in $\cal D$  inducing an extended
generalized hexagon on $\cal D$, there is only one
set of 6-cliques making $\Gamma$ locally $\Delta$.
The $\McL$ 2-fold extended generalized hexagon does satisfy
the hypothesis, and we find it isomorphic to $\Gamma$. \eop

\section{The 3-fold extended hexagon related to $2\times \Co_3$}

The purpose of this section is to show that the $2\times \Co_3$
3-fold extended hexagon is the only extension of the $\McL$
2-fold extended hexagon satisfying the hypothesis of
Theorem 1.1, and that this extension can not be extended further.

Assume $\Delta$ is the 2-fold extension related to $\McL$, 
$\cal D$ its point graph, and $\Gamma$ an
extension of $\Delta$ with point graph $\cal G$
as in Theorem 1.1. 
Then $\cal G$ is locally $\cal D$.

Let $(p,q,r)$ be a path of length 2 in $\cal G$,
and set $M$ to be the subgraph induced by $p^\perp \cap r^\perp$ in $\cal D$.
Then the graph $M$ is locally isomorphic to the $\mu$-graph
of the complement of the McLaughlin graph, and thus to the graph $\cal M$ of the previous
section. By Lemma 3.5 we find that $M_q$ is either a $\mu$-graph
in $\Delta$, or the complement of a $\mu$-graph.
For each point $x$ of $M$ denote by $x'$ the unique point at distance 2 from $x$
in $\cal D$ such that $M_x$ is either $x^\perp \cap x'^\perp$ or its complement
in ${\cal D}_x$.

Let $x$ be a point in $M$, and $y\in M_x$.
Then $M_x\cap M_y$ consists of 60 points.
The points of ${\cal D}_y-M_y$ have 45 neighbours in $M_y$, so that we can conclude that
$x'$ is not adjacent to $y$, and $M_x$ consists of the points of ${\cal D}_x$ not adjacent to
$x'$. 
By the same argument we find that  $M_y$ consists of all the points in ${\cal D}_y$
not adjacent to $y'$.

The 81 points at distance 1 from $x'$, but 2 from $x$ 
have $36$ neighbours in $M_x$, and thus also some neighbour in $M_x\cap M_y$.
The 30 vertices at distance 2 from both $x$
and $x'$ have 54
neighbours in $M_x$, and therefore also neighbours in $M_x\cap M_y$.
Thus the point $y'$ has to be equal to $x'$. Since each point at distance
2 from both $x$ and $x'$ is adjacent to some $y\in M_x$ and thus in $M_y$, we find that $M$ contains and hence consists of the 112 points
not adjacent to $x'$.

Now we can count the number of points at distance
$2$ from $r$ in $\cal G$.  
There are exactly $275.112/112=275$ such points.

The point $q'$ is at distance 3 from $r$.
Any point of ${\cal G}_p$ different from $q'$ is either adjacent $q'$,
or adjacent to $r$. Thus all common
neighbours of $q'$ and $p$ are at distance 2 from $r$.
Since the graph ${\cal G}_{q'}$ is connected we find that all neighbours of $q'$ are
at distance 2 from $r$. Moreover, this implies that $q'$ is the unique point at distance 3 from $r$.

But now uniqueness of the graph $\cal G$ is obvious, and, since there is a unique way to
fix a set of 6-cliques in $\cal D$ making it into a 2-fold extended hexagon,
we also obtain uniqueness of $\Gamma$. We have proved:

\begin{prop}
If $\Gamma$ is locally the $\McL$ $2$-fold extended hexagon, then it is isomorphic
to the $2\times \Co_3$ $3$-fold extended hexagon.
\end{prop}

We finish with the following proposition that together with the
Propositions 3.4, 3.6 and 4.1 proves Theorem 1.1.

\begin{prop}
The $2\times \Co_3$ $3$-fold extended hexagon admits no triangular extension.
\end{prop}
 
\pf Suppose $\Gamma$ is a triangular extension of $\Delta$ which
is the $2\times \Co_3$ multiple extended hexagon, with point graph
$\cal G$, respectively, $\cal D$.
Then $\cal G$ is locally $\cal D$.

Fix a path $(p,q,r)$ of length 2 in $\cal G$, and set $M$
to be the induced subgraph of $\cal D$ with point set
$p^\perp \cap r^\perp$.
Then $M$ is locally isomorphic to the $\mu$-graph of $\cal D$.
We show that the assumption that such a graph $M$
exists leads to a contradiction.

Consider $M_q$. By the above this is a graph
on 112  points which is locally the graph $\cal M$ inside
${\cal D}_q$. Moreover, there is a point $y$ in
${\cal D}_q$ such that inside this graph $M_q$ consists
of all the  points not adjacent to $y$.
Let $z$ be the point in $\cal D$ at distance 3 from $y$.
Then $z$ is at distance 2 from $q$ and adjacent to all
points of $M_q$. Let $x$ be a point of $M_q$, and consider
$M_x$ inside $\Delta_x$. This graph is locally $\cal M$
inside ${\cal D}_x$. In particular, $M_x\cap M_q \simeq \cal M$.
But $M_x\cap M_q$ is  the $\mu$-graph
of the points $q$ and $z$ in ${\cal D}_x$, contradicting
the results of the third and fourth paragraph of this section.
\eop

\section{The nontriangular 3-fold  extended \protect\\
hexagon related to $\Co_3$}

In this final section we give a characterization of the 3-fold
extended hexagon on 276 points related to the sporadic group $\Co_3$.
But first we consider $2$-fold  extensions of generalized hexagons:

\begin{prop}
Let $\Gamma$ be an extension of one of the four extended hexa\-gons
related to $G_2(2)$, $PSU_4(3)$, $\HJ$, respectively, $\Suz$. 
Then $\Gamma$ is triangular, and thus isomorphic
to the $\McL$ $2$-fold extended hexagon.
\end{prop}

\pf Let $\Delta$ be one of the four extended hexagons satisfying $(*)$,
such that $\Gamma$ is an extension of $\Delta$.
The point graph of $\Delta$ is strongly regular
with parameters $(v,k,\lambda,\mu)$ say.

Let $x$, $p$ and $q$ be 3 pairwise adjacent points of $\Gamma$.
We want to show that they are in a common circle.

Assume to the contrary that $p$ and $q$ are not adjacent in $\Gamma_x$.
Let ${\cal B}=\{U\cap V-\{x\} \ |\ x\in U\in {\cal C}_p,\ x\in V \in {\cal C}_q,
|U\cap V| >2\}$, and set ${\cal N}= \bigcup_{U\in {\cal B}}\ U$.
Notice that $|{\cal N}|$ equals  $\mu$.
(See also \cite{Cuy:suz}.)

Fix a point $y\in {\cal N}$. Then there is no circle through
$p$, $q$ and $y$. Assume to the contrary that $W$ is such a circle.
There are circles $U\in {\cal C}_p$ and $V\in {\cal C}_q$ with $U\cap V$ 
containing
$x$ and $y$ together with some third point $z$.
Inside $\Gamma_y$ we see that $p$, $q$ and $x$ form a triangle not in a circle.
Thus the distance between the points $p$ and $q$ in the generalized
hexagon $\Gamma_{y,x}$ is 3. However, it is at most 2, since $z$ is collinear to
both $p$ and $q$ inside $\Gamma_{y,x}$. 
Thus indeed, there is no circle through $p$, $q$ and $y$.

This implies that any point of $\cal N$ is not collinear to $q$ inside
$\Gamma_p$. Thus $\cal N$ contains at most $k-\mu$ points. 
Hence, $\mu+1\leq k-\mu$, which is a contradiction in all the four cases.

We can conclude that $\Gamma$ is triangular, and satisfies
$(**)$. By Theorem 1.1, it is isomorphic to the $\McL$ 2-fold
extended hexagon. \eop

\begin{prop}
Suppose $\Gamma$ is an extension of the $\McL$ $2$-fold extended
hexa\-gon. Then $\Gamma$ is isomorphic to the $3$-fold extended
hexagon on $552$ points related to $2\times \Co_3$, or to
its quotient on $276$ points related to $\Co_3$.
\end{prop}

\pf Let $\Delta$ be the $\McL$ 2-fold extended hexagon.

If $\Gamma$ is a triangular extension of $\Delta$, then it satisfies
$(**)$, and, by Theorem 1.1, it is isomorphic to the $2\times \Co_3$
3-fold extended hexagon on $552$ points. 

Suppose $\Gamma$ is not triangular.
Let $p$, $q$ and $x$ be a triple of pairwise adjacent points of
$\Gamma$, such that $p$ and $q$ are not adjacent in $\Gamma_x$.
Let ${\cal B}=\{U\cap V-\{x\} |\ x\in U\in {\cal C}_p, x\in V\in {\cal C}_q,
|U\cap V|>3\}$, and set ${\cal N}=\bigcup _{U\in {\cal B}}\ U$.
Note that $\cal N$ contains 81 points.

Let $y\in \cal N$. Since $\Gamma_y$ is triangular, there is no circle on
$p$, $q$ and $y$. So $q$ and $y$ are not adjacent inside $\Gamma_p$.
The valency of the subgraph of the point graph of $\Gamma_p$ induced
on the vertices not adjacent to $q$ (inside $\Gamma_p$) is 81.
Moreover, this subgraph is connected. Hence all points of
$\Gamma_p$ not adjacent to $q$ inside $\Gamma_p$ are adjacent
to $q$ inside $\Gamma$. Since the complement of the point graph
of $\Delta$ is connected, we find that any two points of 
$\Gamma_p$ are adjacent inside $\Gamma$.
But that implies that all pairs of points of $\Gamma$ are adjacent,
and $\Gamma$ is a one point extension of $\Delta$.

Since the point graph of $\Delta$ is strongly regular with parameters
$k$ and $\mu$ such that $k=2\mu$, we find that $\Gamma$ carries the structure
of a regular 2-graph. (See 
\cite{BCN,GoeSe,Tay}.) 
Now we can finish the proof of the proposition
by referring to \cite{GoeSe}, 
or by considering the universal cover obtained by a standard
construction for 2-graphs, see \cite{BCN,GoeSe,Tay}. 
Its point graph is the Taylor graph of the regular 2-graph.
This cover is a triangular extension of
$\Delta$ and by Theorem 1.1 isomorphic to the $2\times \Co_3$ extension of $\Delta$.
This proves the  proposition. \eop

Theorem 1.2 follows from the two propositions of this section.



\chapter{Extended generalized octagons and the group \He\ }
\label{chap:EGO}
\newcommand{\oO}{\overline{\Omega}}
\newcommand{\oXi}{\overline{\Xi}}
\newcommand{\oU}{\overline{\Upsilon}}
\newcommand{\U}{\Upsilon}
\renewcommand{\OO}{\Pi}
\renewcommand{\a}{\alpha}
\renewcommand{\b}{\beta}
\renewcommand{\g}{\gamma}
\renewcommand{\d}{\delta}
\renewcommand{\oOO}{\overline{\Pi}}
\paragraph{Abstract.}
Let $\G$ be an extended generalized octagon such that the points of a triple
$\{u,v,w\}$ not on a block 
are pairwise adjacent if and only if the distance between $v$ and $w$ in the
local generalized octagon $\G_u$ equals 3 and there is a thick line through
any point of $\G_u$. Then $\G$ is one of the two examples
related to the groups $2\cdot L_3(4).2^2$ and $\He$. It is also shown that $\G$
does not admit further extensions.

\section{Introduction and the results}
A number of sporadic simple groups arise as automorphism groups of extensions
of classical geometries or buildings. The Mathieu groups act on (multiple)
extensions of PG(2,4) or AG(2,3), and the sporadic Fischer groups act on
extensions of the $U_6(2)$-polar space. The Suzuki chain groups act on
extensions of a certain subgeometry of the $O_8^+(2).3$-building. 
The latter groups also, along with $\McL$, $\Co_3$, $\HS$, $\He$
and $\Ru$, act on extensions of generalized polygons. The group $\Co_2$ acts on
an extension of the $U_6(2)$-dual polar space, the group $BM$ acts on an
extension of $^2E_6(2)$-building. These geometries were characterized under
the assumption that they admit a flag-transitive automorphism group, see e.g.
\cite{BH,DFGMP,Mei,Ivnv:BM,We:he,We:MCL,Pas:EGQIII}.

Purely combinatorial characterizations of the extensions of the
projective plane of order 4 and the affine plane of order 3, that is,
without any assumption on group actions, were obtained by Witt in the
1930's. Recently, similar characterizations of some other sporadic
geometries were obtained. The author gave such characterizations of
the geometries for the Suzuki chain groups related to the
$O_8^+(2).3$-building \cite{Pa:suz}, of those related to the sporadic
Fischer groups \cite{Pa:fi,Pa:epolsp}, and of the extensions of the
generalized quadrangle of order (3,9) related to the groups $\McL$ and
$\Co_3$ \cite{Pa:mclco3}.
Cuypers characterized the extended hexagons related to the Suzuki chain groups
and the geometry for $\Co_2$ \cite{Cuy:co2,Cuy:suz,Cuy:ls}.
Cuypers, Kasikova and the author characterized the multiple extensions of a
generalized hexagon related to $\McL$ and $\Co_3$ \cite{CKP}. 

Here we consider a class of extended generalized octagons. In general, they
have infinite universal covers, see \cite{Pa:covers}, thus additional
conditions are needed to characterize finite examples. 
The example $\G$ related to the group $\He$ satisfies the following condition.
\begin{itemize}
\item[$(*)$] The points of a triple $\{u,v,w\}$ of points not on a block of
$\G$ are pairwise adjacent if and only if $v$ and $w$ are at distance 3 in
the local generalized octagon $\G_u$.
\end{itemize}
There is one more example of an extended generalized octagon satisfying $(*)$,
it relates to the group $2\cdot L_3(4).2^2$.
It turns out that $(*)$, along with a natural assumption on point residues,
characterizes those two examples.
\begin{thm}
Let $\G$ be an extension of a nondegenerate generalized octagon satisfying
$(*)$. Assume that for any point $u$
the local generalized octagon $\G_u$ has at least one
thick line through each point. Then $\G$ is isomorphic either to the extended
generalized octagon on $2048$ points related to the group $\He$, or to the extended
generalized octagon on $112$ points related to $2\cdot L_3(4).2^2$.
\label{ThA}
\end{thm}
R.~Weiss \cite{We:he} characterized  these extended generalized octa\-gons as
flag-tra\-nsi\-tive geometries satisfying a property related to $(*)$ and having
diagram 
\begin{center}
$$\node\stroke{\subset}\node\farc\node.$$
\end{center}
The example related to $2\cdot L_3(4).2^2$ admits a 2-fold quotient
satisfying the property $(*)_{3,4}$ (cf. \cite{Cuy:suz})
obtained from $(*)$ by replacing the words
``distance 3" to the words ``distance 3 or 4". This geometry appears in the
list in \cite{We:he}, as well. Cuypers mentioned in \cite{Cuy:suz} that the
methods he uses there could be applied to the extended generalized octagons
satisfying $(*)_{3,4}$. 

Further, we settle the question about further extensions of the geometries
under consideration. 
\begin{thm}
Let $\G$ be as in Theorem \ref{ThA}. Then $\G$ does not admit any further extensions.
\label{ThB}
\end{thm}
Note that R.~Weiss shows in \cite{We:he} that there are no flag-transitive
geometries with the diagram
\begin{center}
$$\node\arc\node\stroke{\subset}\node\farc\node.$$
\end{center}
and with residues of the left-hand side type of elements as in the remark
following Theorem \ref{ThA}.

\section{Definitions and notation}
\label{sect:defs}
Let $\G$ be an {\em incidence system} $(\cP,\cB)$ of {\em points} and
{\em blocks}, where the latter are subsets of $\cP$ of size at least
2.  The {\em point graph} of $\G$ is the graph with vertex set $\cP$
such that two vertices $p$ and $q$ are {\em adjacent} (notation
$p\perp q$) if there is a block containing both of them.  The distance
between subsets $X$, $Y$ of points of $\G$ (notation $d(X,Y)$ or
$d_{\G}(X,Y)$) is the minimal distance in the point graph between a
point of $X$ and a point of $Y$.  We denote by $\G_n(X)$ the set of
points at distance $n$ from $X\subseteq\cP$, we also use $\G(X)$
instead of $\G_1(X)$.  If $X=\{x\}$ we often use $x$ instead of
$\{x\}$.  For any $X\subseteq\cP$, we denote $X^\perp=\{p\in\cP\mid
p\perp x {\rm\ for\ any\ } x\in X\}$.  The {\em subsystem} of $\G$
induced by the set $X$ is the incidence system $\G(X\cap\cP,\{B\cap
X\mid B\in\cB, |B\cap X|>1\})$.  We call $\G$ {\em connected} if its
point graph is connected.

An incidence system is called a {\em generalized $2d$-gon} if its
point graph has diameter $d$, for any $p\in\cP$ and $B\in\cB$ there
exists a unique point on $B$ closest to $p$, and for any $p,q\in\cP$
such that $d(p,q)=i<d$ there is a unique block on $p$ containing a
point at distance $i-1$ from $q$.  Since it follows that there is at
most one block on any pair of points, the blocks of $\G$ are usually
called {\em lines}, and the word {\em collinearity} is used instead of
the word {\em adjacency}.  Given $\G$, we define the dual system
$\G^*$, in fact also a $2d$-gon, whose points (respectively lines) are
lines (respectively points) of $\G$, incidence is by the inverse
inclusion. We say that $\G$ is {\em nondegenerate} if for each point
there exists a point at distance $d$ from it, and the same holds in
$\G^*$. A line of $\G$ is {\em thick} if it contains more than 2
points, otherwise it is {\em thin}. We say that $\G$ is {\em regular}
of order $(s,t)$ if each point is on exactly $t+1$ lines and each line
contains exactly $s+1$ points. The numbers $s$ and $t$ are called {\em
parameters} in this case.  The famous Feit-Higman theorem \cite{FeHi}
says that a line-thick finite regular $2d$-gon satisfies
$d\in\{1,2,3,4,6\}$. An incident point-line pair of $\G$ is called a
{\em flag}. Given a generalized $2d$-gon $\G$, one can always
construct a $4d$-gon $\G^F$, whose points are flags of $\G$ and whose
lines are the points and the lines of $\G$, incidence being the
natural one.  Note that it is customary to refer to a generalized
4-gon as a {\em generalized quadrangle} (GQ, for short, or GQ$(s,t)$).
Similarly, we refer to generalized 8-gons as {\em generalized
octagons} (respectively GO and GO$(s,t)$).

In what follows a well-known family of generalized quadrangles with
parameters ($s,s$), usually called $W(s)$ (see \cite{PayTh}, a
standard reference on GQ's), will play a significant role.  For a
prime power $s$, the points and the lines of $W(s)$ are the points and
the lines of the projective space PG($3,s$) which are totally
isotropic with respect to a nondegenerate alternating form $f$.
Without loss of generality, $f$ can be chosen as follows:
\begin{equation}
f(x,y)=x_1y_4-x_2y_3+x_3y_2-x_4y_1.
\label{eq:f}
\end{equation}
Note that $W(s)^*\cong W(s)$ if and only if $s=2^k$.

Let $\G$ be an incidence system.
Given $X\subset\cP(\G)$ satisfying $|B-X|>1$ for any $B\in\cB(\G)$, 
we refer to the incidence system 
$$\G_X=(\G(X), \{B-X\mid X\subset B\in\cB(\G)\})$$
as the {\em residue} of $X$, or as the 
{\em local} system in $X$. Let $\cD$ be a class of incidence systems. We say
that $\G$ is an {\em extension} of $\cD$ (or {\em extended} $\cD$) if, for any
$x\in\cP(\G)$, $\G_x$ is isomorphic to a member of $\cD$. If $\cD=\{\D\}$,
it is customary to refer to extensions of $\cD$ as extensions of $\D$
(or extended $\D$). 
If $\G$ is an extension of $\cD$ then 
the connected components of $\G$ are extensions of $\cD$,
also. Hence, unless otherwise stated, we assume our extensions to be connected.

In particular if $\cD$ is a class of generalized octagons, we say that $\G$
is an {\em extended generalized octagon}. 
Note that $\G$ can be considered as a
geometry with the following diagram (for this notion see e.g. \cite{Bue:geo:spo}).
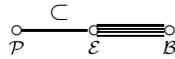
\begin{center}
$$\node_\cP\stroke{\subset}\node_\cE\farc\node_\cB$$
\end{center}
Here $\cE$ is the set of edges of the point graph of $\G$.
Note that the subgraph $\Xi$ of the point graph of $\G$ induced by $x^\perp-\{x\}$
is not necessarily isomorphic to the point graph $\Upsilon$ of $\G_x$. 
In particular, $\Xi$ and $\Upsilon$ are not isomorphic
if $\G$ is an extended generalized octagon satisfying $(*)$, 
namely two vertices of  $\Xi$ are adjacent if and only if the distance between them in $\Upsilon$ is 1 or 3.
 
\section{The type of point residues}
In what follows $\G=\G(\cP,\cB)$ denotes an extended GO satisfying the
conditions of Theorem \ref{ThA}. The main result of
this section is as follows. 
\begin{pr}
$\G$ is an extension of the generalized octagon $W(s)^F$ of flags of the
generalized quadrangle $W(s)$, where $s=2^k$.
\label{Pr1}
\end{pr}
In the remainder of this section we prove Proposition \ref{Pr1}. 
Then in Lemma \ref{LA1} we summarize a few facts arising mainly as by-products
of the proof. 
 
The following statement can be easily deduced from the main result of
Yanu\-sh\-ka \cite{Ya:order}.
\newpage
\begin{res}
{\rm (cf. \cite{Ya:order})} Let $\D$ be a nondegenerate GO such that for any point
$p\in\cP(\D)$ there exists a thick line on $p$. Then exactly one of the
following holds.
\begin{itemize}
\item[$(i)$]       $\D$ is regular of order $(s,t)$ for some $s>1$, $t\ge 1$.
\item[{\rm $(ii)$}]  $\D$ is a GQ$(s_1,s_2)^F$ for some $s_1<s_2$, $s_1>1$
and $s_2>1$.
\item[{\rm $(iii)$}] $\D$ is a $G^F$, where $G$ is the $(a+1)\times (b+1)$-grid
                        for some $a>1$ and $b>1$.\qed
\end{itemize}
\label{Res1}
\end{res}
Below we say that a GO is of {\em type } (x) if it belongs to the family (x)
of Result \ref{Res1}.

Let $u\in\cP$. Pick $x,z\in u^\perp-{u}$, $x\not=z$, $\{x,u,z\}\subseteq B\in\cB$.
If $\G_u$ is of type $(ii)$ (respectively, $(iii)$)
we make our choice to satisfy $|B|=s_1+2$ (respectively,
$|B|=3$ and the size of the other block on $\{u,x\}$ to be $a+2$). 
Choose $v\in x^\perp-u^\perp$ satisfying $v\perp_{\G_x}z$.
\begin{lem}
Let $\,\O$ be the subsystem of the residue $\,\G_u$ induced by $v^\perp$, and
$\Xi$ the
connected component of $\,\O$ containing $x$. Then 
\begin{itemize}
\item[$(a)$] $d_{\G_y}(u,v)=2$ for any $y\in\Xi$. 
\item[$(b)$] There is exactly one thin line of $\,\Xi$ on $y$, and the
remaining lines of $\,\Xi$ on $y$ coincide with those of $\,\G_u$, for any
$y\in\Xi$.
\end{itemize}
\label{L1}
\end{lem}
\begin{pf}
Note that $u\perp_{\G_x}z\perp_{\G_x}v$. In $\G_x$ we find that for any line $L$
on $u$ but not on $z$ each point $w\in L$ distinct from $u$ 
is at distance 3 from $v$, so $w\perp v$ by $(*)$.
The remaining line on $u$ has all its points, except $z$, at distance 2 from
$v$, so $z$ is the only point on it adjacent to $v$. 
Thus we have that $x^{\perp}\cap\Xi$ contains each
line $L$ of $\G_u$ on $x$ such that $z\not\in L$, and $xz$ (as a line of $\Xi$)
is thin. 
This proves part $(b)$ for $y=x$.

By the choice of $z$, there is a thick line of $\Xi$ on $x$, namely any line
of $\Xi$ on $x$ missing $z$ is thick.
On the other hand, let us 
choose $v'\in x^\perp$ satisfying $d_{\G_x}(u,v')=4$. Let $\Xi'$ be the
connected component containing $x$ of the subsystem of $\G_u$ induced by $v'^\perp$. 
Then  $x^{\perp}\cap\Xi'$ 
intersects each line of $\G_u$ on $x$ in exactly 2 points, that is, all the
lines of $\Xi'$ on $x$ are thin.  This gives us a criterion for checking
whether $d_{\G_y}(u,v)$ is $2$ or $4$.

Clearly $d_{\G_z}(u,v)=2$. If $y\in    
x^{\perp}\cap\Xi-\{x,z\}$ then the line $xy$ of $\G_u$ lies within
$\Xi$. By the choice of $z$, $xy$ is thick. 
Using the observation above in the case
$d_{\G_x}(u,v')=4$ we see that $d_{\G_y}(u,v)=2$. Hence  
$d_{\G_y}(u,v)=2$ for any $y\in\Xi$. Thus part $(a)$ is proved.

To complete the proof of $(b)$ it suffices to repeat
the first paragraph of the proof with $y$ in place of $x$ and with $z'\in
y^\perp$ satisfying $u\perp_{\G_y}z'\perp_{\G_y}v$ in place of $z$. 
\end{pf}
\begin{lem}
Let $\Xi$ be as defined in Lemma $\ref{L1}$, and set 
$X'=\Xi(x)\cup\Xi_3(x)$. Then
\begin{itemize}
\item[{\rm $(i)$}]      $|X'|=1+st+s^2t^2+st(1+st-s)^2$,
\item[{\rm $(ii)$}]     $|X'|=1+2s_2+s_2^2$,
\item[{\rm $(iii)$}]    $|X'|=1+a+b+ab$,
\end{itemize}
according as $\G_x$ is of type {\rm $(i)$, $(ii)$} or {\rm $(iii)$}.
\label{L3}
\end{lem}
\begin{pf}
By Lemma \ref{L1} $(b)$, we have that $|\Xi(x)|=1+st$, $1+s_2$,
or $1+a$ according as the type of $\G_u$ is $(i)$, $(ii)$, or $(iii)$.
Since the point graph of $\Xi$ behaves as a ``cactus'' in the first three
layers, it is straightforward to count $|\Xi_3(x)|$, and hence $|X'|$. 
\end{pf}

\begin{lem}
Let $\D=\G_x$, and set  
$$X=\bigcup_{ i=1,3\atop j=1,3}\D_i(u)\cap\D_j(v),$$
where $u$, $x$ and $v$ are as defined above.
Then 
\begin{itemize}
\item[{\rm $(i)$}]      $|X|=1+2st+s^2t^3+s^2t(t-1)$,
\item[{\rm $(ii)$}]     $|X|=1+s_1+s_2+s_1s_2$,
\item[{\rm $(iii)$}]    $|X|=2+2a$,
\end{itemize}
according as $\G_u$ is of type {\rm $(i)$, $(ii)$} or {\rm $(iii)$}.
Moreover, $X=x^\perp\cap\O-\{x\}$, where $\,\O$ is as defined in Lemma
$\ref{L1}$. 
\label{L2}
\end{lem}
\begin{pf} 
We begin by proving the first part of the lemma.
Observe that  the type
of $\D$ can be easily determined, using the blocks on $\{u,x\}$.
Namely, it is the same as the type of $\G_u$, with the exception that if $\G_u$ is
of type $(iii)$ then the corresponding parameters $b$ for $\G_u$ and $\D$ are
not necessary equal. 

For $\D$ of type $(i)$, we determine $|X|$ using standard calculations with
parameters of distance-regular graphs \cite{BCN}, as follows. The point graph of
$\D$ is distance regular with intersection array
$$\{s(t+1),st,st,st;1,1,1,t+1\}.$$
Given any two vertices $x$, $y$ satisfying $d_\D(x,y)=k$, the number of
vertices $w$ such that $d_\D(x,w)=i$ and $d_\D(y,w)=j$ is a constant
$p^k_{ij}$. Note that $|X|=p^2_{11}+2p^2_{13}+p^2_{33}$ and that
$p^2_{11}=1$, $p^2_{13}=st$. The constant $p^2_{33}=s^2t(t^2+t-1)$
is computed using the recurrence formulae in \cite[Lemma 4.1.7]{BCN}. 
The consideration of this case is complete.

For $\D$ of type $(ii)$, it is clear that
$$|\D(u)\cap\D(v)|+|\D(u)\cap\D_3(v)|+|\D_3(u)\cap\D(v)|=1+s_1+s_2.$$
Now let us consider $\Phi=\D_3(u)\cap\D_3(v)$. Let $vpqw$ be a path from $v$
to $\Phi$. Since there are exactly two lines on each point, it follows that
$p\in\D_3(u)$ and $q\in\D_4(u)$. By the same reason for any
$q\in\D_2(v)\cap\D_4(u)$ there is a unique point $w_q\in q^\perp\cap\Phi$. 
Moreover, $w_q=w_{q'}$ implies $q=q'$, for otherwise $\D$ possesses a
circuit on $q$, $q'$ and $w_q$ of length less than 8.
Hence $$|\Phi|=|\D_2(v)\cap\D_4(u)|=s_1s_2.$$ 
The consideration of this case is complete.

The remaining case of $\D$ of type $(iii)$
is dealt with by arguments similar to those used for type $(ii)$.

The last part of the lemma immediately follows from $(*)$.
\end{pf}

\begin{lem}
$\G_p$ is isomorphic to a GO$(s,1)$, for some constant $s>1$ independent from
particular choice of $p\in\cP$.
\label{EGO:L4}
\end{lem}
\begin{pf}
Note that $X'\subseteq X$, where $X$ and $X'$ are as defined in Lemmas \ref{L2} and
\ref{L3}.
It follows that $t=1$ if $\G_u$ is of type $(i)$. 
If $\G_u$ is of type $(ii)$ then $X'\subseteq X$ implies that
$s_1\geq s_2$, a contradiction. If $\G_u$
is of type $(iii)$ then $X'\subseteq X$ implies that $b=1$, which is
also a contradiction. 

Thus $\G_u$ is a GO$(s,1)$. Let $w\in u^\perp-\{u\}$.
In $\G_w$ we see that there are exactly two lines
on $u$, both of size $s+1$. Hence, by Result \ref{Res1}, $\G_w$ 
is a GO$(s,1)$, as well, and by the connectivity of $\G$ the result follows.
\end{pf}

Our next task is to investigate $\Xi$ more closely. 
\begin{lem}
$\O=\Xi$, that is, $\,\O$ is connected. Moreover,
$\O$ is a sub-GO of $\,\G_u$. It is of type {\rm $(iii)$} such that 
$a=b=s$. 
\label{L5}
\end{lem}
\begin{pf}
Using Lemmas \ref{L1} and
\ref{EGO:L4}, we observe that each point
of $\Xi$ has one thin line and one line of size $s+1$ through it. 
Combining Lemmas \ref{L3} and \ref{L2} with $t=1$, we have $X=X'$.
Hence $x^{\perp_\Xi}\cup\Xi_3(x)=x^\perp\cap\Xi=x^\perp\cap\O$.
It follows that every thick line of $\O$ intersects $x^\perp$, since each
line of $\G_u$ intersects $x^\perp$. Hence the
first part of the lemma holds, and the diameter of $\O$ is at most 4.

Since $\O$ has exactly $s$ thin lines $ww'$ satisfying
$d_\O(x,w)=d_\O(x,w')-1=2$, it has (exactly) $s$ thick lines $L$ at distance 3
from $x$, since for each $w$ with $d_\O(x,w)=2$ and such $L$
there exists $w'\in L$ such that $ww'$ is one of the $s$
thin lines just mentioned. Hence $|\O_4(x)|=s^2$. Also, there are $s^2$
points $p\in\O_3(x)$ such that the thick line on $p$ is at distance 2
from $x$. If $p$ is such a point then there exists a thin line $pp'$
satisfying $p'\in\O_4(x)$. It follows that for any $y\in\O_4(x)$ the
thin line on $y$ is at distance 3 from $x$. Hence $\O$ is a GO.
It is clear that it is of type $(iii)$ with $a=b=s$. 
\end{pf}

Let $\D=\G_u$, $x,z\in u^\perp-\{u\}$ and $x\perp_\D z$. 
Pick $L,M\in\cB(\D)-\{xz\}$ such that  $z\in L$ and $x\in M$.
We claim that there are $s$ points $p$ 
at distance 2 from $u$ in $\G_x$ satsifying
$L,M\subset u^\perp\cap p^\perp$ and  such that 
$u^\perp\cap p^\perp\not=u^\perp\cap p'^\perp$ 
for any pair $p,p'$ of those $s$ points. 

Pick $v\in x^\perp-u^\perp$ such that $v\perp_{\G_x}z\perp_{\G_x}u$.
We check that the set of $s$ points on the line $vz$ of $\G_x$
not containing $z$ satisfies the requirements
above. Let $p\in vz-\{v,z\}$. Clearly $L,M\subset u^\perp\cap p^\perp$. On
the other hand, we see within $\G_v$ that there is a point 
$q\in u^\perp\cap v^\perp$ such that $p\not\perp q$.  
Hence $u^\perp\cap v^\perp\not=u^\perp\cap p^\perp$, as required. 

Thus $\D$ admits at least $s$ distinct subGO's isomorphic to $\O$ which
contain $L$
and $M$, where $\O$ is the subsystem of $\D$ induced by $v^\perp$.
Let $\Theta$ be the GQ$(s,s)$ such that $\D=\Theta^F$, and $L,M$ are
(collinear) Points of $\Theta$ (we use capital `P' to distinguish Points of
$\Theta$ from points of $\D$). It is clear that $\O$ corresponds to a subGQ$(1,s)$
of $\Theta$, among the $s$ subGO's of $\D$ just constructed. So we have 
distinct subGQ$(1,s)$'s ${\cal S}_1,\dots,{\cal S}_s$ containing the Points $L$
and $M$ of
$\Theta$. Pick a Point $P\in L^{\perp_\Theta}-M^{\perp_\Theta}$. It is easy
to see that there is at most one subGQ$(1,s)$ through $P$ and $M$. Denoting
${\cal S}_i\cap LP=\{L,P_i\}$, we see that $P_i=P_j$ implies $i=j$, for
$i,j\in 1,\dots,s$.

We have established that for any $Q\in\cP(\Theta)-M^{\perp_\Theta}$ there is
a subGQ$(1,s)$ on $M$ and $Q$. In other words,
$|\{M,Q\}^{\perp_\Theta\perp_\Theta}|=s+1$, 
that is, $M$ is {\em regular}, cf. \cite{PayTh}.
Therefore each Point of $\Theta$ is regular, for we are free to choose
$M$ to be any line of $\D$. Hence $\Theta\cong W(s)$ by \cite[5.2.1]{PayTh}.

Let $N$ be a line of $\D$ intersecting $M-\{x\}$. 
The lines $N$ and $xz$ of $\D$ correspond to the collinear Points of
$\Theta^*$, the dual of $\Theta$. 
Repeating the argument above with $N$ and $xz$ in place of $L$
and $M$, we see that $\Theta^*\cong W(s)$ also. Hence $s=2^k$, cf.
\cite{PayTh}. The proof of Proposition \ref{Pr1} is complete.\qed

Several facts arising in the proof above will be required
later. We summarize them and some other facts in the following lemma.
\begin{lem}
Let $\,\Xi$ be a subGO of $\,\D=\G_u=\Theta^F$ isomorphic to $\,\O$. 
Then there exists
$q\in\G_2(u)$ such that $\,\Xi=u^\perp\cap q^\perp$. There are two classes
${\cal O,O}^*$ of such subGO, where ${\cal O}$ (respectively ${\cal O}^*$)
corresponds to the
subGQ$(1,s)$'s of $\,\Theta$ (respectively $\,\Theta^*$),  
$|{\cal O}|=|{\cal O}^*|=s^2(s^2+1)/2$,  and $\,\O\in{\cal O}$.
There are $(s-1)(s+1)^2$ elements $\,\Xi\in{\cal O}$ such that
$|\O\cap\Xi|=2s+2$ and $\,\O$ and $\,\Xi$ have two thick lines in common.
For all other $\,\Xi\in{\cal O}-\{\O\}$ we have $|\O\cap \Xi|=0$.\\
For $\,\Xi\in{\cal O}^*$, $\,\O\cap\Xi$ is either an ordinary octagon or
empty.
\label{LA1}
\end{lem}
\begin{pf}
Without loss of generality we can assume that the lines $L$, $M$ of $\D$
defined in the final part of the proof of Proposition \ref{Pr1}
lie in $\Xi$. We saw that there are exactly $s$ subGO's containing
$L$ and $M$, and all of them are of the form $q^\perp\cap u^\perp$ for some
$q\in\G_2(u)$. Thus $\Xi$ is among them, and the first statement is proved. 

The second statement just repeats an observation made at the end of the proof
of Proposition \ref{Pr1}, and the formula for $|{\cal O}|=|{\cal O}^*|$ is
well-known. 

Let $\Xi\in{\cal O}-\{\O\}$. If $\oO\cap\oXi$ contains a Line then there are
exactly 2 Points in $\oO\cap\oXi$. There are $(s+1)^2$ Lines in $\oO$ and there
are $s$ subGQ$(1,s)$ of $\Theta$
through each Line. Thus there are $(s-1)(s+1)^2$ elements of ${\cal O}$ 
intersecting $\O$ in two thick lines. If $\O\cap\Xi$ is nonempty then
$\oO\cap\oXi$ contains a flag. Hence $\O\cap\Xi$ contains a thin line.
Thus it also contains two thick lines, and we are in the already
considered situation. 

It remains to consider the case $\Xi\in{\cal O}^*$. 
Here $\oXi$ is a subGQ$(s,1)$ of $\Theta$.
Assume that $\oO$ and $\oXi$ have a common Point. It is easy to deduce
that they intersect in an ordinary quadrangle of $\Theta$, so $\O\cap\Xi$ is
an 8-gon. Otherwise, that is if $\oO$ and $\oXi$ do not have common Points,
$\O\cap\Xi=\emptyset$.
\end{pf} 

\section{More on point residues}
In this section we show that $\G$ is an extension of either $W(2)^F$ or
$W(4)^F$, and as a by-product we obtain more information about the set
$\G_2(u)$, where $u\in\cP$.
Let $v\in\G_2(u)$.
\begin{lem}
Let $x\in u^\perp\cap v^\perp$. Then, for any $y\in u^\perp\cap v^\perp$, 
$d_{\G_y}(u,v)=d_{\G_x}(u,v)=d\in\{2,4\}$.
\label{L7}
\end{lem}
\begin{pf}
If $d=2$ then the statement follows from Lemma \ref{L5}. The only other
possibility is $d=4$.
\end{pf}

\begin{lem}
Let $a,b\in u^\perp\cap v^\perp$ satisfy
$a\perp_{\G_u}c\perp_{\G_u}b$ for $c\in u^\perp-\{u\}$. Then $c\perp v$.
\label{EGO:L8}
\end{lem}
\begin{pf}
By Lemma \ref{L5}, if $d_{\G_a}(u,v)=2$ then the subsystem of 
$\G_u$ induced by $v^\perp$ is geodetically closed, so the statement holds in
this case. By Lemma \ref{L7} we may assume that 
$d_{\G_a}(u,v)=d_{\G_b}(u,v)=4$. Observe that the subsystem $\O$
of $\G_b$ induced by $a^\perp$ is a subGO, and that 
$uc$ is a thin line of $\O$. Since $d_{\G_b}(u,v)=4$, 
$v$ is at distance 3 from any line on $u$ of $\O$. Hence $d_\O(v,c)=3$,
so $v\perp c$.
\end{pf}

Let $d_{\G_x}(u,v)=4$ for some (and so, by Lemma \ref{L7}, for any) $x\in
u^\perp\cap v^\perp$, and let $\O$ be the subsystem of $\G_u$ induced by $v^\perp$.
It is easy to see that the connected components of $\O$ are ordinary
polygons. The following is an immediate consequence of Lemma \ref{EGO:L8}.
\begin{lem}
The connected components of $\O$ are ordinary polygons.
The distance between any two of them is at least 3.\qed
\label{L9}
\end{lem}
It turns out that the components are the polygons of the least possible
girth.
\begin{lem}
The connected components of $\O$ are ordinary $8$-gons. 
\label{L10}
\end{lem}
\begin{pf}
Let $axb$ be a 2-path contained in a component of $\O$.
Let $ubb_1dvca_1a$ be an
8-gon inside $\G_x$ containing $u$ and $v$ (points are listed in the natural
cyclic order). 
This 8-gon is unique, since $\G$ is extended GO$(s,1)$.
Since $d_{\G_x}(u,c)=d_{\G_x}(u,d)=3$, we have $c\perp u\perp
d$. Since $x,u\in a^\perp\cap c^\perp$ and $d_{\G_x}(a,c)=2$, by Lemma
\ref{L7} we have that $d_{\G_u}(a,c)=2$. So there exists $p=p_{ac}\in
u^\perp$ satisfying
$a\perp_{\G_u}p\perp_{\G_u}c$. By Lemma \ref{EGO:L8}, $p\in v^\perp$. 

Similarly we find $p_{bd},p_{cd}\in u^\perp\cap v^\perp$. So
$xap_{ac}cp_{cd}dp_{bd}b$ is an 8-gon within $\O$ containing $axb$.
By Lemma \ref{L9}, it is a full connected component of $\O$.
\end{pf}

Denote by $\Xi$ the connected component of $\O$ containing $axb$,
$\{c,d\}=\Xi_3(x)$.
Considering $\G_x$, we easily find that $X=\{u,v,x\}^\perp-\Xi$ satisfies
$|X|=4(s-1)^2$ and that $X\subseteq c^\perp\cup d^\perp$.
Moreover for any $y\in X$ there exists the unique $y'\in X-\{y\}$ such that
$y\perp_{\G_u}y'$. Denote by $\U$ the connected component of $\O$ containing
$y$. Let $y_1,y_2$ be the two points of $\U$ at distance 3 from $y$. As
$x\in\{u,v,y\}^\perp-\U$, we see that $x\in y_1^\perp\cup y_2^\perp$. So $\U$
has at least 4, and so exactly 4, of its points within $X$.
Summarizing, we have the following.
\begin{lem}
The set $X$ is the disjoint union of $(s-1)^2$ sets of size $4$
of the form $X\cap\U$, where $\U$ is a connected component of $\O$.\qed
\label{L11}
\end{lem}
%
Let $q\in\{x,v\}^\perp-u^\perp$ satisfy $d_{\G_x}(u,q)=2$.
We shall establish a correspondence between the subsystem
$\OO$ of $\G_u$ induced by $q^\perp$ and the subsystem $\OO'$ of
$\G_q$ induced by $u^\perp$.
\begin{lem}
Let $xz$ be a thin line of $\OO$. Then it is a thin line of $\OO'$, as well.
Let $\{z,a_1,\dots,a_s\}$ be the thick line on $z$ of $\OO$ and for
$i=1,\dots,s$ let $a_ia_i'$ be
the thin line on $a_i$ of $\OO$. 
Then the thick line of $\OO'$ on $x$ is $\{x,a_1',\dots,a_s'\}$.
\label{L12}
\end{lem}

\begin{pf}
The first statement follows immediately from the observation that for any
$x\in u^\perp\cap q^\perp$ the thin line $xz$ of $\OO$ is defined by
$q\perp_{\G_x}z\perp_{\G_x}u$. Let $L$ be the thick line of $\OO'$ on $x$,
and let $y\in L-\{x\}$. Note that $d_{\OO'}(z,y)=d_{\G_q}(z,y)=2$. Hence by Lemma
\ref{L7} we have that $d_{\G_u}(z,y)=2$. So either $y$ is one of the
$a_i'$ or $y\perp_{\G_u}x$. But the latter is clearly
impossible. Hence the lemma.
\end{pf}
 
\begin{lem}
The subsystem of $\OO$ induced by $v^\perp$ is the disjoint union of
$s$ copies of the $3$-path.
Moreover if $\a\b\g\d$ and $\a'\b'\g'\d'$ are any 
two of them  written in the natural order then, without loss in generality,
$\a,\g,\d\in\OO_3(\a')$ and $\a,\b,\d\in\OO_3(\d')$.
\label{L13}
\end{lem}
\begin{pf}
Let $\OO'$ be the subsystem of $\G_q$ induced by $u^\perp$.
Let $\Theta\cong W(s)$ such that $\G_q=\Theta^F$ and $\OO'$ corresponds
to a subGQ$(1,s)$ $\oOO'$ of $\Theta$. Denote by $ab$ the thin line of $\OO'$
at distance 2 from $v$. Let $A,B$ be the thick lines of $\OO'$ satisfying
$a\in A$, $b\in B$. Let $v$ correspond to the flag $(P,l)$
of $\Theta$. Then $P\in AB\in\cB(\oOO')$. 

For any $W\in\cP(\oOO')-\{A,B\}$ there exists a Line $w$ through $W$
intersecting $l$. There are $2s$ such Points $W$. 
As they correspond to thick lines of $\OO'$, they give us $2s$ thick
lines (and so $2s$ points)
of $\OO'$ at distance 3 from $v$. Also, there are $2s$ points
at distance 3 from $v$ lying on the thick lines $A,B$
of $\OO'$. Clearly, there are no more points at distance 3 from $v$ within
$\OO'$.
It shows that the subsystem of $\OO'$ induced by $v^\perp$ consists
of $4s$ points.

Next, for $w$ as above let $w\cap \oOO'=\{W,W'\}$.
This defines an equivalence relation with classes of size 2 on
$\cP(\oOO')-\{A,B\}$, and so gives $s$ thin lines of $\OO'$ at distance 3
from $v$.
The translation of this situation back into $\OO'$ is presented diagrammatically
in Figure \ref{Fig1}.
\begin{figure}[bhpt]
\leavevmode
\epsffile{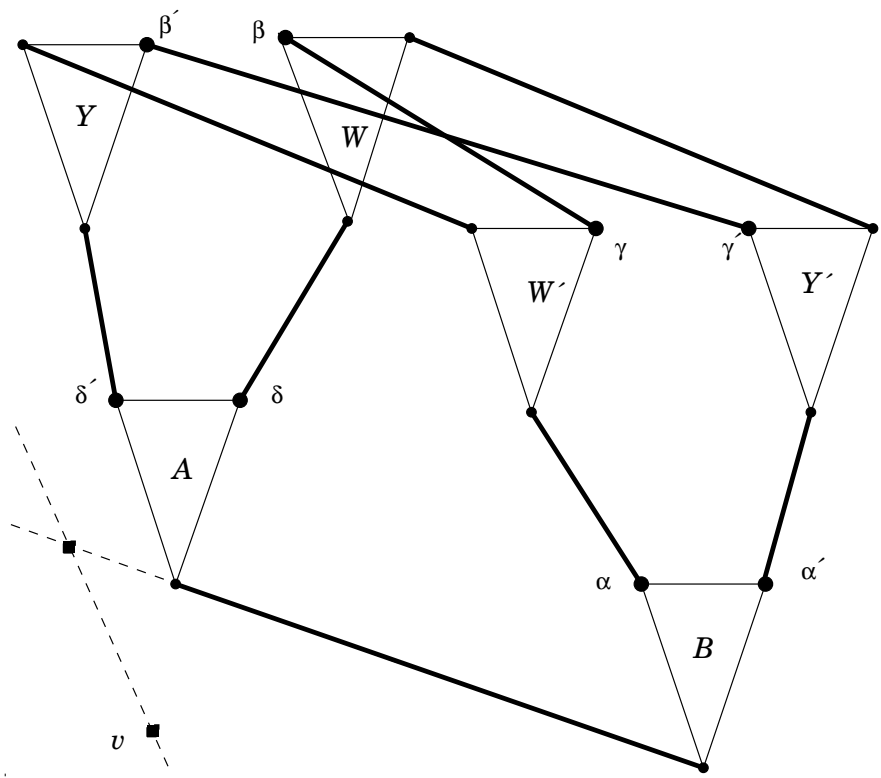}
\caption{The subsystem $\OO'$ of $\G_q$, see the proof of Lemma \protect\ref{L13}.
(Unless $s=2$, only a part of $\OO'$ is actually shown.)
Points inside $\OO'$ are shown as circles, whereas those outside are shown as
squares. Thick lines of $\OO'$ are shown as triangles, thin lines of $\OO'$
are shown as bold lines. Points of $\OO'$ at distance 3 from $v$ are shown by
bigger circles than the rest of them.
}
\label{Fig1}
\end{figure}
Now we can use Lemma \ref{L12} to look at the pointset we are interested in 
within $\OO$. By Lemma \ref{L12} we have that
$\a\perp_{\G_u}\beta\perp_{\G_u}\g\perp_{\G_u}\d$, and none of the points just listed
are collinear within $\G_u$ to any other point at distance 3 from $v$ within $\OO$.
The same holds for $\a'$, $\b'$, $\g'$ and $\d'$, and it is straightforward to
check the remaining claimed distances. 
\end{pf}

Now we return to the connected components $\Xi$ and $\U$ of $\O$, that we
considered in the text preceding Lemma \ref{L11}. We shall improve the result
of Lemma \ref{L11}. As above, denote $X=\{u,v,x\}^\perp-\Xi$. 
By Lemma \ref{L11}, the subsystem $X\cap\U$ has two connected components
with sets of points $\{y_1,y_1'\}$ and  $\{y_2,y_2'\}$ respectively. 
Let $\{y,y'\}=\Xi_3(x)$. As already noted, each point of $X$ is at
distance 3 within $\G_u$ from either $y$ or $y'$ (but not both of them).
\begin{lem}
Let $d_{\G_u}(y,y_1)=3$. Then $d_{\G_u}(y',y_1')=3$.\\
Also, $d_\U(\{y_1,y_1'\},\{y_2,y_2'\})=2$.
\label{L14}
\end{lem}
\begin{pf}
First, we claim that there exists $q\in\{x,y,v,y_1\}^\perp-u^\perp$ such that
$u\perp_{\G_y}w\perp_{\G_y}q$ for $w\in u^\perp\cap q^\perp$.
Note that $d_\D(y,y_1)=d_\D(v,y_1)=3$, where $\D=\G_x$. Let $L$ be the unique
line of $\D$ at distance 1 from $v$, $y$ and $y_1$, and let $z$ be the
point on $L$ satisfying $\{z\}=L\cap\D_3(u)$. Denote
$\{q\}=z^{\perp_\D}\cap\D_2(u)$. Clearly $d_\D(v,q)=d_\D(y,q)=3$ and
$d_\D(y_1,q)=1$ or $3$. We are done.

Now we apply Lemma \ref{L13}.
Since $x,y,y_1\in q^\perp$, the subsystem $Q$ of $\O$ induced by $q^\perp$
contains the 3-path joining $x$ and $y$ and a 3-path of $\U$ containing
$y_1$. In the notation of Lemma \ref{L13}, without loss in generality let us
take $\a'=x$, $\d'=y$.
Since $y_1$ is at distance 3 from both $x$ and $y$, either
$\a=y_1$ or $\d=y_1$. 

If $\a=y_1$ then, by Lemma \ref{L13}, $y_1'\not\in Q$ and, by the same
lemma,
$\{\b\}=\U\cap\{y,y_1\}^\perp\not=\{y_1'\}$. The line $y_1\b$ of $\G_u$ is
at distance 2 from $y$, so the line $y_1 y_1'$ is at distance 3 from $y$.
Hence $y\not\perp y_1'$. Therefore $y'\perp y_1'$.
 
If $\d=y_1$ then $y_1'=\gamma$ and so $y\not\perp y_1'$.
Hence again $y'\perp y_1'$.
This completes the proof of the first part of the lemma.

To prove the second part, observe that by Lemma \ref{L13} the system
$X\cap\U\cap Q$ has two components and that the distance between them equals 2.
\end{pf}

At this point we are able to analyze $\O$ yet more precisely.
We shall establish that $s=2$ or $4$ and that $\O$ is as in the known examples.
We translate the situation with the 8-gons of $\O$ into the corresponding one
with the 4-gons in $\Theta$, where $\D=\G_u=\Theta^F$.

It is obvious that an (nondegenerate) 8-gon of $\D$ corresponds to a
quadrangle $ABCD$ of $\Theta$, such that the points of the 8-gon are the
flags $(A,AB)$, $(A,AC),\dots,(C,AC)$ of $\Theta$. Denote by $\oO$ the set of
4-gons of $\Theta$ corresponding to the set of 8-gons of $\O$. 
We write the Points of 4-gons in the natural cyclic order. We call
the Lines $AB$, $AC$, $BD$ and $CD$ of $\Theta$ the {\em Sides} of $ABCD$.
The next statement follows immediately from Lemma \ref{L9}.
\begin{lem}
Let $P$, $P'$  (respectively $L$, $L'$)
be the sets of Points (respectively, of Sides) 
of two distinct $4$-gons of $\oO$ corresponding to the
$8$-gons $\U$, $\U'$ of $\O$.
Then $P\cap P'=L\cap L'=\emptyset$ and $Y\not\in s'$ for any $Y\in P$ and
$s'\in L'$.     \qed
\label{L15}
\end{lem}
\begin{lem}
With the notation of Lemma $\ref{L15}$, let $d_\D(\U,\U')=3$.
Then for any $Y\in P$ there exists a unique
$Y'\in P'$ such that $Y\perp_\Theta Y'$, and for any $s\in L$ there exists
a unique $s'\in L'$ intersecting $s$.
\label{L16}
\end{lem}
\begin{figure}[bhpt]
\leavevmode
\epsffile{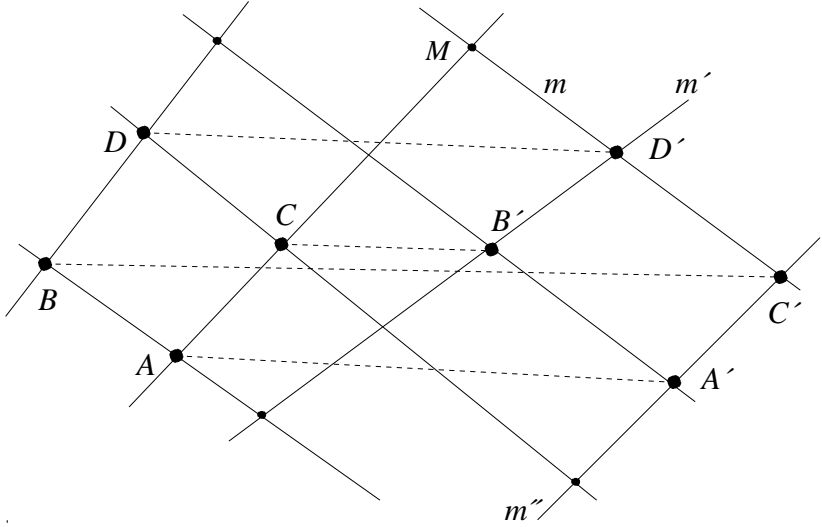}
\caption{
Points of $\Theta$ are shown as circles.
Points inside $\oO$ are shown as bigger circles than Points outside. 
Sides of $\oO$ are shown by continuous lines.
Other Lines of $\Theta$ are shown by discontinuous lines.
}
\label{Fig2}
\end{figure}
\begin{pf}
In this proof we shall write $\perp$ instead of $\perp_\Theta$.
Figure \ref{Fig2} illustrates the argument.
Let $P=ABDC$, $P'=A'B'C'D'$.
Without loss of generality $A\perp A'$.
So the lines $l_1=(A',A'B')(A',A'C')$ of $\U'$  and $(A,AB)(A,AC)$ of $\U$
are at distance 3 in $\D$. 
There is one more line $l_2$ of $\U'$ at distance 3 from $(A,AB)$ such that 
$l_2$ and $(A,AB)(B,AB)$
are at distance 3 in $\D$. By Lemma \ref{L14}, $d_{\U'}(l_1,l_2)=2$.
Without loss of generality, $l_2=(B',B'D')(D',B'D')$. It follows that the
Lines $AB$ and $B'D'$ of $\Theta$ intersect.
Since there must be two lines of $\U'$ at distance 3 from $(B,AB)$, and one
of them is $l_2$, we see that $B\perp C'$.
Continuing in the same vein we find that $AC$ intersects $C'D'$, $C\perp B'$,
and $BD$ intersects $A'B'$.
Finally, we see that $D\perp D'$ and $CD$ intersects $A'C'$.

Thus the existence of $Y'$ and $s'$ mentioned in the statement is
established.
The uniqueness easily follows from Lemma \ref{L11}.
\end{pf}

Our next goal is to reconstruct the whole of $\O$ starting from a single 8-gon
$\U\subset\O$.

Let $\oU$ and $\oU'$ be the 4-gons of $\Theta$ corresponding to $\U$ and
$\U'$ as in Lemma \ref{L16}. We shall show that the choice of the point
$M=AC\cap C'D'$ and the line $m=C'D'$, see Figure \ref{Fig2}, determine $\U'$
uniquely. That is, if there exists $\U'\subset\O$ having the Side $m$ then 
it is determined uniquely. 

Let $M$ and $m$ be chosen. 
Note that the points of $\oU'$ collinear to $B$ or to $D$ lie on $m$.
By Lemma \ref{L16}, $m\cap BD=\emptyset$. So the two points $D'$ and $C'$ of
$\oU'$ satisfying $D'\perp_\Theta D$, $C'\perp_\Theta B$ are distinct and
uniquely determined. 
Next we see that the Side $m'$ of $\oU'$ on $D'$ distinct from $m$ must
intersect $AB$. It determines $m'$. Similarly we reconstruct the remaining
Side on $C'$.
Finally we are forced to set $A'=A^{\perp_\Theta}\cap m''$,
and $B'=C^{\perp_\Theta}\cap m'$.
Thus $\oU'$, and so $\U'$, is determined.

There are $s-1$ possible choices of $M$ and $s-1$ possible choices of $m$
through $M$. These gives a total of at most $(s-1)^2$ 8-gons of $\O$ at
distance 3 from $x\in\U$.
Since by Lemma \ref{L11} this number must be exactly $(s-1)^2$, and hence
the 8-gons of
$\O$ at distance 3 from $x$ are determined once $\U$ is chosen.

Since $\Theta\cong W(s)$, $s$ even, the 8-gons of $\G_u$ all lie in one orbit of
$Aut(\G_u)$. Thus there is no loss of generality in choosing $\oU$ to have 
points $A=(1000)$, $B=(0010)$, $D=(0001)$, $C=(0100)$. (Here we present
$\Theta$ as it was explained in Section \ref{sect:defs} Collinearity is
determined by (\ref{eq:f}).)

Choose $M=(1100)$ and $m=MD'$, where $D'=(0111)$. We see that
$\oU'=\oU(m,M)$ has the following points, where the notation is chosen to
match those of Figure \ref{Fig2}:
$A'=(1110)$, $B'=(1101)$, $C'=(1011)$, $D'$ as above.

Let $M_1=(1a00)$, where $a\in GF(s)-GF(2)$, $m_1=M_1D_1$, where $D_1=(011a)$.
Let us find $\oU''=\oU(m_1,M_1)$ with points $A''$, $B''$, etc.
First, observe that $D\perp D_1$, so $D''=D_1$,
and $B''=B^\perp\cap m_1=(a^{-1},0,1,a)$.
Since $D'\perp D''$, it follows that $d_{\G_u}(\oU',\oU'')=3$.
Therefore there exists a point $Y$ of $\oU'$ satisfying $Y\perp B''$.

If $Y=C'$ then $f(C',B'')=a+a^{-1}=0$, so $a^2=1$, contradicting the
choice of $a$. (Here $f$ is the form in (\ref{eq:f}).)

If $Y=B'$ then $1+a+a^{-1}=0$, so $a^2+a=1$. The latter implies
$a^3+a^2=a$, so $a^3=1$. Hence in this case $s\le 4$.

Finally, $Y=A'$ implies $1+a=0$, again contradicting the choice of $a$, and
$Y=D'$ is impossible.
So $s\le 4$.

To summarize, we state the following.
\begin{pr}
$\G$ is either an extension of $W(2)^F$ or an extension of $W(4)^F$. \qed
\label{Pr17}
\end{pr}

\section{The remaining cases}
We shall reconstruct the point graph of $\G$ layer by layer.
Let $u\in\cP(\G)$.
The set $\G_2(u)$ is naturally divided into two parts $\G_2'(u)$ and
$\G_2''(u)$ according as $d_{\G_x}(u,v)=2$ or $4$, respectively, 
where $v\in\G_2(u)$ and
$x\in u^\perp\cap v^\perp$, cf. Lemma \ref{L7}.
We proceed by showing that $\G$ is isomorphic to one of the known examples.
Table \ref{T1} is self-explanatory. It shows certain information on the
subsystems $\O$ of $\D=\G_u$ induced by $v^\perp$, $v\in\G_2(u)$, for the 
examples related to the groups $2\cdot L_3(4).2^2$ (here $s=2$) and $\He$
(here $s=4$).
\begin{table}[bhpt]
\begin{center}
\newcommand{\Down}[2]{\raisebox{-#2ex}[1ex][.5ex]{#1}}
\begin{tabular}{|c|c|c|c|c|c|c|}
\hline
\Down{$s$}{1}
   &location             &\multicolumn{5}{|c|}{location of $x\in\D-\O$}\\\cline{3-7}
   & of $v$              &$d_\D(x,\O)$      &1    &2   &3   &4 \\       \hline
   &\Down{$\G_2'(u)$}{1} &$|x^\perp\cap\O|$ &6    &8   &    &  \\       \cline{3-7} 
\Down{2}{1}
   &                     &$\# x$            &9    &18  &0   &0   \\     \cline{2-7}
   &\Down{$\G_2''(u)$}{1}&$|x^\perp\cap\O|$ &6    &8   &8   &0 \\       \cline{3-7} 
   &                     &$\# x$            &16   &8   &4   &1   \\     \hline
   &\Down{$\G_2'(u)$}{1} &$|x^\perp\cap\O|$ &10   &16  &    &  \\       \cline{3-7} 
\Down{4}{1}
   &                     &$\# x$            &75   &300 &0   &0   \\     \cline{2-7}
   &\Down{$\G_2''(u)$}{1}&$|x^\perp\cap\O|$ &22   &16  &32  &  \\       \cline{3-7} 
   &                     &$\# x$            &240  &60  &45  &0   \\     \hline
\end{tabular}
\end{center}
\caption{Sets $x^\perp\cap\O$ in the examples, 
where $\O=u^\perp\cap v^\perp\subset\G_u$. }
\label{T1}
\end{table}
By Lemma \ref{LA1}, for any $v\in\G_2'(u)$ the subsystem  $\O$ of
$\D=\G_u$ induced by $v^\perp$ is a subGO,
and each subGO of $\D$ isomorphic to $\O$ appears in this form for a point of
$\G_2'(u)$.
We shall prove similar statement for $\G_2''(u)$.
\begin{lem}
Let $s=2$, $\Phi=\cP(\D)$ and $m=4$, respectively let
$s=4$, $\Phi$ be the set of subGO$(2,1)$'s of $\D$ forming an orbit of
$Aut(\D)$ of length $1360$ and $m=3$.\\
Then for any $v\in\G_2''(u)$ there exists unique $\phi\in\Phi$ such that
$$u^\perp\cap v^\perp=\{x\in\cP(\D)\mid d_\D(x,\phi)=m\}.$$
Moreover $\phi$ is at maximal distance in $\D$ from 
$u^\perp\cap v^\perp$. 
\label{L18}
\end{lem} 
\begin{pf}
The statement of the lemma merely describes the situation with $\G_2''(u)$
and $\Phi$ in the known examples.
All we need to show is that once $\D$ is given, the set $\G_2''(u)$ and the
edges joining $u^\perp$ and $\G_2''(u)$ are uniquely determined.

Let $v\in\G_2''(u)$. Pick  an 8-gon $\U$ of $\O=u^\perp\cap v^\perp\subset\D$.
We know that $Aut(\D)$ acts transitively on the 8-gons of $\D$, so there is
no loss in generality in choosing $\O=u^\perp\cap v^\perp$.
We saw in the final part of the proof of Proposition \ref{Pr17} that once
$\U$ is chosen, the set $\O^x$
of the 8-gons of $\O$ at distance 3 from a point $x$ of
$\U$ is determined. It is easy to check that $\O^x$ does not depend upon the
particular choice of the point $x$, so we denote $\O^\U=\O^x$. 
It is straightforward to check that 
$\Pi=\{\U\}\cup\O^\U=\{\Xi\}\cup\O^{\Xi}$ for
any 8-gon $\Xi$ of $\Pi$.
It follows that $\O$ is the disjoint union of several Components (to
distinguish from 8-gons of $\O$) of the form
$\{\Xi\}\cup\O^{\Xi}$, where $\Xi$ is an 8-gon of $\O$.
Each Component is one of the sets arising in the known examples, so it
remains to show that there is only one Component.

Clearly the Components are at distance 4 from each other.
But in both cases the size of the set of points of $\D$ at distance 4 
from a Component is less than the size of a Component. We are done.
\end{pf}

It immediately follows from Lemma \ref{L18} and the observation preceding it
that the partial subgraph of $\G$
consisting of the points $u^\perp\cup\G_2(u)$ and the edges at distance at
most 1 from $u$ is isomorphic to the same subgraph in the example.
Next task is to show that the edges within $\G_2(u)$ are uniquely determined,
as well.

\begin{table}[bhpt]
\tabcolsep=3pt
\newcommand{\Down}[2]{\raisebox{-#2ex}[1ex][.5ex]{#1}}
\begin{center}
\begin{tabular}{|c|c|c|c|c|c|c|c|c|c|c|c|}
\hline
\Down{$s$}{1}
   &location   &         &\multicolumn{9}{|c|}{location of $x$}\\\cline{4-12}
   & of $v$    &         &\multicolumn{3}{|c|}{$\G_2'(u)-\{v\}$}
                                 &\multicolumn{6}{|c|}{$\G_2''(u)-\{v\}$}\\\hline
   &\Down{$\G_2'(u)$}{1} &$|\{x,u,v\}^\perp|$ &8 &6 &0 
                &\multicolumn{2}{|c|}{$8^*$}
                 &\multicolumn{2}{|c|}{$8$}
                  &\multicolumn{2}{|c|}{$4$}\\ \cline{3-12} 
\Down{2}{1}  
   &                     &$\# x$            &9    &9  &1   
        &\multicolumn{2}{|c|}{9}
         &\multicolumn{2}{|c|}{18}
          &\multicolumn{2}{|c|}{18}\\ \cline{2-12} 
   &\Down{$\G_2''(u)$}{1} &$|\{x,u,v\}^\perp|$ &$8^*$ &8 &4 
        &8 &6
         &\multicolumn{2}{|c|}{5}
          &\multicolumn{2}{|c|}{4}\\ \cline{3-12} 
   &                     &$\# x$            &4    &8  &8   
        &4 &16
         &\multicolumn{2}{|c|}{16}
          &\multicolumn{2}{|c|}{8} \\\hline 
   &\Down{$\G_2'(u)$}{1} &$|\{x,u,v\}^\perp|$ &10 &8 &0 
                &\multicolumn{2}{|c|}{16}
                 &\multicolumn{2}{|c|}{8}
                  &\multicolumn{2}{|c|}{0}\\ \cline{3-12} 
\Down{4}{1}
   &                     &$\# x$            &75    &100  &$60+36$   
        &\multicolumn{2}{|c|}{300}
         &\multicolumn{2}{|c|}{$100+900$}
          &\multicolumn{2}{|c|}{60}\\ \cline{2-12} 
   &\Down{$\G_2''(u)$}{1} &$|\{x,u,v\}^\perp|$ 
                &16 &8        &0  &32 &25  &22  &12  &8   &0 \\ \cline{3-12} 
   &                     &$\# x$ 
                &60 &$20+180$ &12 &45 &144 &240 &720 &180 &30 \\\hline 
\end{tabular}
\end{center}
\caption{The intersections of sets $u^\perp\cap v^\perp$ in the examples.}
\label{T2}
\end{table} 
In Table \ref{T2} we present information about the sets $\{x,u,v\}^\perp$,
where $x,v\in\G_2(u)$. We denote there by $8^*$ the entries corresponding to
the subsystems of $\G_u$ not isomorphic to those mentioned in Lemma \ref{L13},
that is the disjoint union of 2 copies of the 3-path.

Comparing Tables \ref{T1} and \ref{T2}, we see that the edges within
$\G_2(u)$ are uniquely determined, as well. As an example, let as consider
the case $s=4$, $v\in\G_2'(u)$.
In Table \ref{T1} we see that the set $v^\perp-u^\perp-\{v\}$  consists of 75
points $x$ such that $|\{u,v,x\}^\perp|=10$ and 300 points $x$ such that  
$|\{u,v,x\}^\perp|=16$. In Table \ref{T2} we look at the corresponding
entries and see that in $\G_2(u)$ there are exactly 75  points satisfying
$|\{u,v,x\}^\perp|=10$ and exactly 300 points satisfying 
$|\{u,v,x\}^\perp|=16$. Hence the $v^\perp$ is determined.
 
Thus the consideration of the case $s=4$ as complete, since by Table \ref{T1}
the diameter of the point graph is 2.

Let us turn to the case $s=2$. We see that for any $v\in\G_2''(u)$ there is a
unique $z_v=v^\perp\cap\G_3(u)$. Each $w\in v^\perp\cap\G_2''(u)$ lies in a
block on $v$ and $z_v$. Hence $z_w=z_v$.
The subgraph of the point graph induced on $\G_2''(u)$ is isomorphic to those
induced on $\G(u)$, in particular it is connected. Hence $\G_3(u)=\{z_v\}$
and $z_v^\perp=\{z_v\}\cup\G_2''(u)$.

The proof of Theorem \ref{ThA} is complete.

\section{Proof of Theorem \protect\ref{ThB}}
Let $\G$ be an extension of $\cD$, where $\cD$ be the class of
extended generalized octagons satisfying the conditions of Theorem \ref{ThA}.

Pick $u\in\cP(\G)$, $x\in u^\perp-\{u\}$. Note that 
$\G_{ux}$ is a GO$(s,1)$, for some $s\in\{2,4\}$. By the connectivity argument, the
value of $s$ does not depend on the particular choice of $u$ and $x$.
There exists a point $v\in\cP(\G_x)$
such that for some $y,z\in\{u,v,x\}^\perp$ one has 
$u\perp_{\G_{xy}}z\perp_{\G_{xy}}v$. 
Define the incidence system $\O$ as follows. Let 
$$ \cB_0(\O) =  \{B\cap B'\mid v\not\in B\in\cB(\G_u),\ 
u\not\in B'\in\cB(\G_v),\ |B\cap B'|=3\}$$
and let $\cP(\O)=\bigcup_{B\in \cB_0(\O)} B$. 
Then let
$$\cB(\O)=\cB_0(\O)\cup \{B-\{u\}\mid v\not\in B\in\cB(\G_u),\ 
B-\{u\}\subset\cP(\O)\}.$$
Note that $\O_x$ is a subGO of $\G_{ux}$ of type $(iii)$, where $a=b=s$.
Since there exist blocks of $\G$ containing, respectively, $\{u,x,y,z\}$ and
$\{v,x,y,z\}$, one has that $\O_y$ is a subGO of $\G_{uy}$ and
$\O_y\cong\O_x$ for any $y\in\cP(\O_x)$ and an appropriate choice of $z$.

The idea of the following argument is adapted from \cite{CKP}.
By Lemma \ref{LA1}, there exists $r\in\cP(\G_{ux})$ such that $\O_y$ is the
subsystem of $\G_{uy}$ induced by $r^{\perp_{\G_u}}$.
In particular, $y\perp_{\G_{ux}}w\perp_{\G_{ux}}r$ for $w\in\cP(\G_{ux})$.

By the choice of $r$, we have $wy\cap\cB(\O_x)=\{w,y\}$, where
$wy\in\cB(\G_{ux})$. Hence the other line $wr$ of $\O_x$ on $w$ is thick.
In particular, $r\in\cP(\O_x)$.
Repeating the argument above with the roles of $y$ and $r$ interchanged, we
have that the line $wr$ of $\O_x$ is thin.
This is the contradiction.

The proof of Theorem \ref{ThB} is complete.

\bigskip

{\bf Acknowledgement.} The author used the computer algebra system
{\sf GAP}
\cite{GAP} and its shared package GRAPE \cite{GRAPE} to investigate
combinatorial properties of the extended generalized octagons
for the groups $\He$ and $2\cdot L_3(4)$ and compute contents of Tables  
\ref{T1} and \ref{T2}.\nocite{vBCC} \nocite{Ti:la}

\chapter{Review}
\label{chap:review}
Diagram geometries are a tool for investigating multi-dimensional
incidence geometries starting from certain information on their
2-dimensional residues. In this chapter we survey the development of
the theory of diagram geometries, with the emphasis on
characterizations of geometries for finite simple groups, extensions
of buildings, and graphs with a given neighbourhood. The main goal is
to place the research carried out in this thesis into the
context of our current understanding of diagram geometries.

\section{A brief history}
Although diagram geometries, as they are currently known, were introduced
only in the late seventies by Buekenhout~\cite{Bue:diag}, similar
structures have been investigated since the first axiomatic studies of
projective and affine spaces in the beginning of this
century, cf. Veblen and Young~\cite{VY}.  Later, in the thirties,
Witt~\cite{Wi1,Wi2} realized that the 5 sporadic simple groups
discovered by Mathieu~\cite{Math1,Math2,Math3} last century can be
adequately described as automorphism groups of (repeated) extensions
of the affine plane of order 3 or the projective plane of order 4.

During  the  fifties, Veldkamp~\cite{Ve} started  investigating  axiom
systems  for polar spaces. This eventually led Tits~\cite{Ti:bldgs} to
introduce  buildings of  spherical  type as a  special case of diagram
geometries.   At  the  same time,  an  enormous amount of activity was
directed   at    classifying   finite   simple   groups,    see   e.g.
Aschbacher~\cite{As:spo}.  Geometries  of subgroups of  finite  simple
groups played an  important role in this research.  

Perhaps  one  of  the  most notable examples  in this context  was the
investigation  of  3-transposition  groups by  Fischer~\cite{Fi}  (see
Chapter~\ref{chap:fi2n}  for  more  information).  Here  an  essential
combinatorial  tool was the graph of commuting  3-transpositions.   We
shall discuss this in greater detail in Section~\ref{rev:loc3tr}.

Recall  that a {\em  transitive extension} of  a permutation group $H$
acting on $\O'$ is a  transitive  permutation group $G$ acting on  the
set  $\O=\O'\cup\{\infty\}$  such  that the  stabilizer  in $G$ of the
point $\infty$ equals $H$.  Many sporadic simple groups, including the
Mathieu groups, $J_2$, $\HS$,  $\McL$, $\Suz$, $\He$, $\Ru$, Fischer's
sporadic   groups  and  $\BM$,  were  discovered   or  constructed  as
transitive  extensions  of classical  or other sporadic simple groups,
cf. \cite{As:spo}.

The latter discoveries  promoted  interest  in  graphs  with  constant
neighbourhood,  as  a  natural  combinatorial  generalization  of  the
orbital   graphs    of    transitive    extensions.    For   instance,
Shult~\cite{Sh:grext}  derived  graph-theoretic   insight   into   the
existence   of  doubly   transitive  permutation  representations   of
$Sp_{2n}(2)$, $\HS$ and $\Co_3$ along with the  related Taylor  graphs
\cite{BCN}  (note  that  Taylor  graphs  as  a concept were yet to  be
introduced  at  the  time  of  the  publication  of  \cite{Sh:grext}).
Buekenhout  and  Hubaut~\cite{BH}  gave a  local  characterization  of
``natural'' graphs for several of the finite simple  groups, including
$Fi_{22}$.

When  the  group  $H$  mentioned  above is  transitive on  $\O'$,  any
extension $G$ is at least doubly transitive. So  the  orbital graph of
$G$  is  complete  and  carries  no  information.   In  this  case  an
appropriate object to consider is an $H$-invariant design on $\O'$ and
its natural extension.  Finding extensions of  designs (see 1.1.3) was
yet another combinatorial generalization  of  the transitive extension
problem,   and  intense  research  has   been   conducted  there,  cf.
Dembowski~\cite{Dem}.   The  results  by Hughes~\cite{Hug65.1,Hug65.2}
are of  particular relevance to us,  since he classified extensions of
finite projective spaces, which is a special case  of the more general
problem of classifying the  extensions of polar  spaces considered  in
this   thesis,    see    Chapters~\ref{chap:o3n},    \ref{chap:EGQ33},
\ref{chap:fi2n},
\ref{chap:epolsp}.

These  investigations were similar in  their local  approach,  in that
each attempted to derive global properties of the underlying geometric
object  from   its  local   properties.    The  celebrated  paper   by
Buekenhout~\cite{Bue:diag} introduced a common  theoretical background
for these topics.  One of the  major hopes  was to  create  a unifying
theory which encompassed both classical geometries and sporadic simple
group geometries.

\section{Finite simple group geometries}
First we recall some definitions and notation.

A (Tits-Buekenhout) geometry of rank $r$ is defined inductively, as follows.
\begin{description}
\item[$(i)$] A {\em geometry of rank $1$\/} is a set with at least
two elements (namely, the null graph with at least 2 vertices).
\item[$(ii)$] A {\em geometry of rank $r\ge 2$} is a pair $(\G,I)$,
where $\G$ is a connected $r$-partite graph and $I$
is its $r$-partition such that for any vertex $x$ of $\G$ the subgraph
$\G(x)$ is a geometry of rank $r-1$ with respect to the partition
of the set of vertices of $\G(x)$ induced by $I$.
\end{description}

The classes  of $I$ are  called {\em types}. Note that $(\G,I)$ can be
recovered  from $\G$. So with a  slight abuse of  notation  we  denote
$\G=(\G,I)$.   The  vertices  of  $\G$  are  called   {\em  elements}.
Elements are {\em incident\/} if they are adjacent in $\G$.  A set $X$
of pairwise incident elements of $\G$  is a {\em  flag}.  Given a flag
$X$ with $|X|<n$, the  geometry $\G(X)$ is the {\em residue\/} of $X$.
The restriction of $I$ to $\G(X)$ is called the {\em cotype\/} of $X$,
whereas the size of this restriction is called  the {\em  corank\/} of
$X$.

Let $X$ be a flag of cotype $\{\omega\}$. $\G$ is said {\em to admit 
order} $w$ in $\omega$ if $w=|\G(X)|-1$ does not depend upon the
particular choice of $X$ of cotype $\{\omega\}$.

We illustrate these  definitions  by  considering the affine  geometry
$\Theta=(\Theta,I)$  of rank $3$ (that  is,  of  dimension  $3$)  over
GF$(q)$.  Here  $I=\{${\it  points}, {\it lines}, {\it planes}$\}$ and
the incidence is given by inclusion.   The  residue of  a point $p$ is
formed by the  lines  and  the planes  on  $p$  and  is isomorphic  to
PG$(2,q)$.   The residue of a plane  $\pi$ is formed by the points and
the lines on  $\pi$ and is AG$(2,q)$. The residue of  a line is formed
by the points  and  planes  on  it,  and  is a  generalized digon (see
below).  Note  also  that the  orders  of  $\G$  are $q-1$, $q$,  $q$,
respectively.

\medskip

The  rank  2 residues  of a  geometry  $\G$ already contain  a lot  of
information  about  $\G$.    This  suggests  the  idea  of  describing
isomorphism  types of rank 2  residues in a compact graphical way.   A
{\em diagram} $\cD$ of $\G$ is a (possibly directed) graph with vertex
set $I$, such that an edge $(u,v)$ is labeled by the isomorphism types
of residues of flags of cotype $\{u,v\}$.  Also, the vertices of $\cD$
are sometimes labeled by the names of types along with their orders.

Common  conventions  have  been  developed  for the  diagrams of  some
important rank 2 geometries.

A    diagram     for    a    {\em    generalized    $n$-gon\/}    (see
Chapter~\ref{chap:intro}) is  drawn as  $\node\stroke{{}^{(n)}}\node$.
The symbol $(n)$ is usually omitted for  $n=2,3,4,6,8$ and $0,1,2,3,4$
ordinary  parallel  bonds,  respectively,   are  drawn  instead.   For
instance, to denote a generalized quadrangle of order $(s,t)$ we  draw
$\node_s\darc\node_t$.  We define a {\em generalized digon\/}  to be a
rank 2 geometry whose incidence graph is complete bipartite.
\medskip

Similarly, a {\em circular space\/}, which is  the incidence system of
elements  and (unordered)  pairs of elements of  a set, is  denoted by
$\node\stroke{\subset}\node$.  An {\em affine  plane\/} is denoted  by
$\node\stroke{\it Aff}\node$.

\medskip

Thus,   in    our    example    the    diagram    of    $\Theta$    is
$$\node_{q-1}\stroke{\it  Aff}\node_q\arc\node_q.$$  Here  the  types,
from left to  right, are {\em  points, lines} and {\em planes}.  Their
orders are displayed below the corresponding  nodes.  We note that the
residue  of  a  line  is  a  generalized digon,  so, according  to the
convention, no bonds are drawn between {\em  points}  and {\em planes}
nodes.

\medskip

An {\em automorphism\/} of a geometry $\G$ is a type-preserving
automorphism of its incidence graph. The automorphisms of $\G$ form a
group denoted by $\Aut(\G)$. We say that $\G$ is {\em
flag-transitive\/} if $\Aut(\G)$ acts transitively on the sets of
flags of each type.

\bigskip

Most  finite  simple groups  are  known  to  act flag-transitively  on
Tits-Bueken\-hout  geometries  with  diagrams  obtained  from  spherical
Coxeter diagrams by replacing, if  necessary,  one or  more projective
plane    strokes    $\node\arc\node$   by   circular   space   strokes
$\node\stroke{\subset}\node$ or affine plane strokes $\node\stroke{\it
Aff}\node$. In  particular, finite groups of Lie type  of sufficiently
large  rank are automorphism  groups of the corresponding buildings of
spherical type (that is, no $\node\arc\node$ strokes are replaced).

Geometries with spherical Coxeter diagrams, and thus also simple
groups of Lie type, form perhaps the most famous class of geometries
and groups characterized by their diagrams, see Tits~\cite{Ti:la} and
Pasini~\cite{Pas:book}.  However, the finite sporadic simple groups
are left untouched by this classification. The purpose of this thesis
is to try to bridge this gap. 

Note that  the main  body  of work carried  out  in  this direction is
concerned with flag-transitive geometries. So in most results not only
diagram and, maybe, some extra  local geometric assumptions are  made,
but  also  that  the  geometries  in  question  admit  flag-transitive
automorphism    groups.    The   survey   article   by   Pasini    and
Yoshiara~\cite{PaYo}, (see also \cite{Pas:book}) mentions many results
in  this  direction.  More references will follow.  Unfortunately, the
results for flag-transitive geometries  fall  far short  of  those for
buildings   mentioned   above,  and,  more   generally,  for   various
subgeometries of buildings (for instance, affine polar spaces
\cite{CoSh}),  see  e.g. \cite{Pas:book}.  Indeed, the  latter do  not
require any assumptions on the group action.

When  the author  started working  on the  material included  in  this
thesis, only a  few  similar results  for  sporadic  group  geometries
requiring no assumptions on group actions were known.

Now we are going  to  survey classification results on several classes
of  diagram geometries (including  results from this thesis) which are
of similar flavour to those on buildings and related geometries.

\subsection{$c^k.A_n$-geometries}
Historically, the first class of geometries apart from buildings
which has been completely classified is the class of
$c^k.A_n$-geometries. It includes geometries for three of the sporadic
simple groups, namely $M_{22}$, $M_{23}$ and $M_{24}$.

A $c^k.A_n$-geometry is a geometry of rank $k+n$ with diagram
$$\node^1\arc\node^2\cdots\node^{k-1}\arc\node^k\stroke{\subset}
\node^{k+1}_s\arc\node^{k+2}_s\cdots\cdots\node^{k+n-1}_s\arc\node^{k+n}_s.$$
Here the labels above the nodes reference the types.

\smallskip

According to Hughes~\cite{Hug65.1}, any $c^k.A_n$-geometry is 
one of the following.
\begin{itemize}
\item[$(i)$] An affine geometry over GF(2); here $s=2$ and $k=1$.
\item[$(ii)$] One of the Steiner systems $S(3,6,22)$, $S(4,7,23)$ 
or $S(5,8,24)$ for the sporadic simple groups $M_{22}$, $M_{23}$, 
or $M_{24}$; here $s=4$, $n=2$ and $k=1$, $2$, or $3$, respectively.
\item[$(iii)$] An extension of a projective plane of order 10; here
$k=1$ and $s=10$. 
\end{itemize}
Meanwhile, $(iii)$ has been ruled out due to the non-existence of
projective planes of order 10, see Lam, Thiel and Swiercz~\cite{PG10}.

\subsection{$c^k.C_n$-geometries}
The class of $c^k.C_n$-geometries includes geometries for the sporadic
simple groups $\McL$, $\HS$, $\Co_3$, $\Suz$, $\Co_1$, $Fi_{22}$, $Fi_{23}$
and $Fi_{24}$, cf. \cite{PaYo}. Before considering results in this
area relevant to us, we refer the reader to \cite{PaYo} for a survey
on results classifying $c^k.C_n$-geometries under a flag-transitivity
assumption.

Chapter~\ref{chap:intro} presents a detailed summary
of the statements of the author's
results in this area covered in the thesis, namely those from
Chapters~\ref{chap:o3n}, \ref{chap:EGQ33}, \ref{chap:fi2n} and
\ref{chap:epolsp}.

A $c^k.C_n(s,t)$-geometry has a diagram
$$\node^1\arc\node^2\cdots\node^{k-1}\arc\node^k\stroke{\subset}
\node^{k+1}_s\arc\node^{k+2}_s\cdots\cdots\node^{k+n-2}_s\arc\node^{k+n-1}_s
\darc\node^{k+n}_t.$$
Here the labels above the nodes reference the types.  

Recall   that  we  consider  geometries  satisfying  the  Intersection
Property  (see  \cite{Bue:diag}).   Thus  their  residues  satisfy the
following  property  (LL).
\begin{enumerate}
\item[(LL)] Given  an  ordered  subset  of  types $S=\{t_1,\dots,t_m\}$ 
such that the subdiagram induced by  $S$ is connected, and  a  residue
$\Pi$ of cotype $S,$ there  is at most one element of  $\Pi$  of  type
$t_2$ incident to two elements of $\Pi$ of type $t_1.$
\end{enumerate}
In particular, it  follows that  the $C_n$-residues are polar  spaces.
From a combinatorial point  of view, this means  that the  cases $n=2$
and $n>2$ are rather  different.  Indeed,  in  the  former  case these
residues  are  generalized  quadrangles  (GQ),  which  form  a  rather
``wild'' variety of objects in  the sense  that  the number  of  their
isomorphism types grows  fast  as their  size  grows,  cf.  Payne  and
Thas~\cite{PayTh}.  This  indicates that  a satisfactory combinatorial
characterization of $c^k.C_2$-geometries for arbitrary  $C_2$-residues
is highly unlikely.

First,  we  discuss  the  case  $n=2$.  Even  in  order  to  obtain  a
classification under the  assumption of flag-transitivity,  it appears
necessary  to  restrict  the  $C_2$-residues  to  those  of  so-called
classical type, that is, to the  GQ's  admitting a classical group  as
flag-transitive automorphism  group, see \cite{PaYo}.   Combining  the
classification  of flag-transitive  automorphism  groups of  classical
GQ's by Seitz~\cite{Seit:73} with results on  doubly transitive groups
(noting that the line stabilizer in such a flag-transitive group is at
least  doubly  transitive on  the line),  one obtains a short list  of
possibilities for GQ's which might appear as $C_2$-residues.

This  suggests  that  any  attempt  to  classify  $c^k.C_2$-geometries
without  a  flag-transitivity   assumption  should  be  restricted  to
$C_2$-residues  contained in the  list just mentioned.  The first such
case is  $s=2$.   A complete classification of $c.C_2(2,t)$-geometries
follows      from      Buekenhout     and     Hubaut~\cite{BH}     and
Buekenhout~\cite{Bue:extpolarsp}.

For $k\ge 2$, most of the {\em known} infinite families of
$c^k.C_2(2,t)$-geometries are characterized in Hall and
Shult~\cite{HaSh}, and Pasechnik~\cite{Pa:3tr} (Chapter~\ref{chap:o3n}
contains \cite{Pa:3tr}). In these classifications it is essentially
assumed that all the $c^{k_0}.C_2$-residues are given for some fixed
$k_0$.  For a more detailed treatment see Chapter~\ref{chap:intro}.

\smallskip
The second case to consider is  $s=3$.   Note  that $t=1$, $3$, $5$ or
$9$ are the only possibilities, since these are the only values of $t$
for   which  a   GQ$(3,t)$   exists,   cf.~\cite{PayTh}.   A  complete
classification of $c.C_2(3,1)$-geometries can be found in Blokhuis and
Brouwer~\cite{BlBr44} and Fisher~\cite{Fi44}.  The examples satisfying
(T) (that is, the property that any subset of collinear points lies in
the   residue   of   some   element)   have    been    classified   by
Pasechnik~\cite{Pa:EGQ33} (Chapter~\ref{chap:EGQ33} contains
\cite{Pa:EGQ33}) for  $t=3$,  and by  Makhnev~\cite{Makhps}
for $t=5$.  For $t=9$, Pasechnik~\cite{Pa:mclco3}
classified all examples.

For $k>1$ the following results are known. For $t=3,$ Pasechnik
(unpublished) and Cuypers (personal communication) independently
obtained a characterization of the only known infinite family of
$c^k.C_2(3,3)$-geometries. Again the assumption is made that for
certain fixed $k_0$ all $c^{k_0}.C_2$-residues belong to this family.
For a description of this family see e.g. \cite{PaYo} or
Meixner~\cite{Mei}.  For $t=9$, Pasechnik~\cite{Pa:mclco3} classified
the $c^2.C_2(3,9)$-geometries. Also, in \cite{Pa:mclco3} it has been
shown that if a $c^k.C_2(3,9)$-geometry satisfies (T) then $k<3$.

Note  that the unique  $c.C_2(3,9)$-geometry  admits  $\Aut(\McL)$  as
flag-transi\-tive     automorphism      group,      and     the      two
$c^2.C_2(3,9)$-geometries   admit  $\Co_3$  or   $2\times  \Co_3$   as
flag-transitive automorphism groups.

\smallskip
The only other values of $s$ and $t,$ for which several results of the
same  flavour  are   known,  are   $s=4$   and  $t=2$.   Blokhuis  and
Brouwer~\cite{BBegq42}  classified  $c.C_2(4,2)$-geometries satisfying
(T), such  that  their point graphs  have $\mu$-graphs  (that is,  the
subgraphs  induced  on the common  neighbourhood of  two  vertices  at
distance 2)  with $K_{3,3}$ as the only possible connected components.
As \cite{BBegq42}  remains unpublished,  we note that this  result can
also  be derived from  Pasechnik~\cite{Pa:3tr} (Chapter~\ref{chap:o3n}
contains  \cite{Pa:3tr}).  An  example of a $c.C_2(4,2)$-geometry  not
satisfying    this   property   on    $\mu$-graphs    is    given   in
Pasechnik~\cite{Pa:epolsp} (see Chapter~\ref{chap:epolsp}).

An infinite family of $c^k.C_2(4,2)$-geometries has been
characterized in Pasechnik~\cite{Pa:3tr} 
(Chapter~\ref{chap:o3n} contains \cite{Pa:3tr}). 
This result plays an
important role in the classification of $c^k.C_3$-geometries, see
Chapter~\ref{chap:fi2n}.

\bigskip
In the  case $n>2$  the $C_n$-residues  are finite nondegenerate polar
spaces,  which  are comple\-tely classified. Moreover, using the results
on $c^k.A_{n-1}$-geomet\-ries (see above), we have  that either $s=2$ or
$s=4$ and  $n=3$.  For  $s=2$, a complete classification  is  given by
Buekenhout and Hubaut~\cite{BH} and Buekenhout~\cite{Bue:extpolarsp}.

Pasechnik~\cite{Pa:epolsp,Pa:fi}   (see  Chapters~\ref{chap:fi2n}  and
\ref{chap:epolsp}) gives a complete answer in the remaining case under
the  additional  assumption that $t>1$ and  that (T) holds  for $k>1$.
This   is   a   generalization   of   results   obtained   under   the
flag-transitivity assumption in \cite{BH,DFGMP,Mei,vBW:fi:char}.   The
examples  appearing  in the  classification  have  automorphism groups
$Fi_{22}.2$, $Fi_{23}$, $Fi_{24}$ and $3\cdot Fi_{24}$ (here $k=1$, 2,
3 and 3, respectively).

\subsection{Extensions of nondegenerate generalized polygons
satisfying $(**)_J$} 
Several sporadic simple groups arise as
(repeated) transitive extensions of automorphism groups of generalized
$n$-gons for $n>4$.  So they act flag-transitively on (repeated)
extensions of generalized $n$-gons for $n=6$ or 8, that is on a
geometry with diagram
$$\node^1\arc\node^2\cdots\node^{k-1}\arc\node^k\stroke{{}^\subset}
\node^{k+1}_s\stroke{{}^{(n)}}\node^{k+2}.$$
The universal cover of such a geometry may be infinite, e.g. see
Pasini~\cite{Pa:covers}.

However,  it  is  possible to  try  to  classify  all  extensions  not
satisfying some  strong  form of the  Intersection Property.   Roughly
speaking,  we ensure that their point graphs contain ``non-geometric''
cliques (in other words, (T) is  not satisfied).  Under the assumption
of   flag-transitivity,  such   a   program   was   carried   out   by
Weiss~\cite{We:MCL,We:he}.  A combinatorial generalization of  Weiss's
assumption is  the following property $(**)_J$ of the point graph $\G$
of   an  $r$-fold  extension  $\cG$   of  a  (finite,   nondegenerate)
generalized $n$-gon, where  $n=2d>4$ and  $J\subseteq\{2,\dots,n/2\}$.
We  call the elements of the  rightmost type  in the  diagram of $\cG$
{\em blocks}.

An  $r$-fold extension  $\cG$  of  a generalized  $n$-gon  is said  to
satisfy $(**)_J$ if the  following  holds for the  point graph $\G$ of
$\cG$.
\begin{enumerate}
\item[$(**)_J$] Any $(r+1)$-clique of $\G$ is in a block;
a clique $\{x_1, \dots, x_{r+2}\}$ is not in a block
if and only if $x_{r+1}$ and $x_{r+2}$ are at distance $d\in J$ in the 
point graph of the generalized $n$-gon $\cG_{x_1,\dots, x_r}$.
\end{enumerate}

For $J=\{3\}$ the classification of such extensions is completed in
Cuypers \cite{Cuy:suz}, Cuypers, Kasikova and Pasechnik~\cite{CKP}
(see Chapter~\ref{chap:hex}) and Pasechnik~\cite{Pa:he} (see
Chapter~\ref{chap:EGO}).  It seems quite remarkable that the resulting
list of examples is finite; in particular it contains sporadic group
geometries for $J_2$, $\Suz$, $\McL$, $\Co_3$ and $\He$. In other
words, a parallel combinatorial interpretation of the exceptional
group-theoretic situations leading to the existence of the sporadics
just mentioned is found.

Cuypers also considered a wider class of objects, namely extensions of
near polygons satisfying $(**)_J$, see \cite{Cuy:ls}. In particular
a geometry for the sporadic group $\Co_2$ appears in a classification
of this nature, of extended near hexagons satisfying $(**)_{\{2\}}$.

\subsection{$c.A_n.c^*$-geometries}
A $c.A_n.c^*$-geometry is a geometry of rank $n+2$ with diagram
$$\node\stroke{\subset}\node_s\arc\node_s\dots\node_s\arc\node_s
\stroke{\supset}\node.$$
It follows from the classification of $c.A_n$-geometries (see above)
that either $s=2$, or $s=4$ and $n=2$.

The case $n=2$ in the classification of $c.A_n.c^*$-geometries has
been settled by Hughes~\cite{Hug1,Hug2}. The list obtained consists of
two subgeometries of PG(4,2) (here $s=2$) and a geometry with
automorphism group $\Aut(\HS)$ (here $s=4$).

If $n>2$ then, according to Van Nypelseer~\cite{VanNy}, the only
examples are obtained from PG$(n+2,2)$.

\section{Graphs with prescribed neighbourhood}
Let $\cD$ be a class of graphs. A graph $\G$ is said to be {\em
locally} $\cD$ (or {\em locally} $\D$ if $\cD=\{\D\}$) if for each
vertex $v$ of $\G$ the graph $\G(v)$ induced on the set of neighbours of $v$
is isomorphic to a member of $\cD$.
The literature on the topic of local characterization of graphs, that
is finding all graphs which are locally $\cD$ for a given $\cD$, is
rather extensive, see e.g.
\cite{BFdFSh,BlBr44,Coh:loc,BH,HaSh,MS:alt,Pa:3tr,Pa:suz,Wee1,Ha:locpet,BlBr:locK33pet}
and also references in \cite{BlBr:locK33pet}.

Here we give a survey of results which are most relevant to the theme
of this thesis. Namely, the areas we shall touch upon are:
\begin{enumerate}
\item towers of graphs related to finite sporadic simple groups; and
\item towers of locally 3-transposition graphs, including locally
cotriangular graphs, and their generalizations.
\end{enumerate}
A {\em tower} (or a {\em chain}) of graphs $\Si_1$, $\Si_2$,\dots,
$\Si_n$ is a set of graphs such that $\Si_{i+1}$ is locally $\Si_i$
for $i=1,\dots,n-1$.

\subsection{Towers of graphs related to finite sporadic simple groups}
Several towers of graphs related to sporadic simple groups have been
characterized locally. These include graphs from the Suzuki chain
$\Si_1$,\dots,$\Si_5$, namely the graphs with automorphism groups
$PGL_2(7)$, $G_2(2)$, $\Aut(J_2)$, $G_2(4).2$ and $\Aut(\Suz)$ on 14, 36,
100, 416 and 1782 vertices, respectively.  Pasechnik~\cite{Pa:suz}
(see Chapter~\ref{chap:suz}) has shown that for $2\le n\le 4$ a
locally $\Si_n$-graph $\G$ is isomorphic to $\Si_{n+1}$ (or, if $n=4$,
to the 3-fold antipodal cover of $\Si_5$).  Brouwer, Fon-der-Flaas and
Shpectorov~\cite{BFdFSh} have classified all locally $\Si_1$-graphs.
\medskip

Another family is the Fischer chain $\Phi_1$,\dots,$\Phi_4$, namely
the graphs with automorphism groups $\Aut(U_6(2^2))$, $\Aut(Fi_{22})$,
$Fi_{23}$ and $Fi_{24}$ on 693, 3510, 31671 and 306936 vertices,
respectively.  Pasechnik~\cite{Pa:fi,Pa:epolsp} (see
Chapters~\ref{chap:fi2n} and \ref{chap:epolsp}) has shown that $n\le 3$ and 
any locally $\Phi_n$-graph $\G$ is isomorphic to $\Phi_{n+1}$, or, if
$n=3$, then $\G$ can also be isomorphic to the 3-fold antipodal cover of
$\Phi_4$.
\medskip

Next, we consider results on the McLaughlin chain of graphs $\D_1$,
$\D_2$, $\D_3$ with automorphism groups $\Aut(U_4(3))$, $\McL$ and
$2\times \Co_3$ on 112, 275 and 552 vertices, respectively.
Pasechnik~\cite{Pa:mclco3} has shown that if $\G$ is a connected
locally $\D_n$-graph then $n\le 2$ and $\G$ is isomorphic to
$\D_{n+1}$.

There is another chain of graphs, closely related to this one,
$\D'_0$, $\D'_1$, $\D'_2$, $\D'_3$ with automorphism groups
$\Aut(L_3(4))$, $U_4(3).2^2$, $\McL.2$ and $2\times \Co_3$, on 105, 162,
275 and 552 vertices, respectively. Note that $\D'_2$ is the
complement of $\D_2$, and that both $\D_3$ and $\D'_3$ are so-called
{\em Taylor graphs}, see Brouwer et al.~\cite{BCN} and Taylor~\cite{Tay}.
It may be deduced from Cuypers~\cite{Cuy:suz} and Cuypers, Kasikova
and Pasechnik~\cite{CKP} (see Chapter~\ref{chap:hex}) that if $\G$ is
a locally $\D'_n$-graph then $n\le 2$ and $\G$ is isomorphic to $\D'_{n+1}$.

\subsection{Locally 3-transposition graphs and generalizations} 
\label{rev:loc3tr}

Fischer spaces  and  3-transposition  graphs  are  important tools  in
studying centrefree  3-transposition  groups,  see  Fischer~\cite{Fi},
Buekenhout~\cite{Bue:fi} and  J.~Hall~\cite{Ha:3tr:89:1,Ha:3tr:93}.  A
3-transposition group  $G$  is  a group generated  by a conjugacy class
$C=g^G$  of  involutions  such that $(xy)^2=1$  or  $(xy)^3=1$ for any
$x,y\in C$. The  elements  of $C$ are called 3-transpositions.  Often,
$C_G(g)\cap C$ is  again a conjugacy class of $C_G(g)$.  In this  case
$C_G(g)$ is also  a 3-transposition group.  Furthermore, if $\G$ is  a
(commuting) 3-transposition graph of $G$, then the  graph $\G(x)$ is a
3-transposition graph for any vertex $x$.

This prompts a natural combinatorial generalization of the problem of
classifying 3-transposition groups, namely one can consider the
problem of classifying the {\em graphs} which are {\em locally}
3-transposition graphs of a class of 3-transposition groups $\cD$.
From this point of view, the first natural reduction is to consider
such a class $\cD$ for which the corresponding 3-transposition graphs
are nonempty and connected. The second natural restriction is to
consider groups whose 3-transposition graphs satisfy the property that
$\G(x)=\G(y)$ implies $x=y$ for any vertices $x$, $y$ of $\G$. We
shall call such graphs {\em reduced}.

Observe that the following property holds in any 3-transposition graph
$\D$. For any  pair $x$, $y$ of  nonadjacent vertices there  exists  a
unique third vertex  $z$  such  that any  $v\in \V\D$  $-\{x,y,z\}$ is
adjacent  to 0, 1, or 3 vertices in $\{x,y,z\}$.  $\G$ is called  {\em
cotriangular\/}  if, in  addition, $v$  is  adjacent to  at least  one
vertex of  $\{x,y,z\}$, cf. Hall  and Shult~\cite{HaSh}.  It turns out
that every  cotriangular graph is  in fact  a  3-transposition  graph.
Thus the class of 3-transposition graphs is a natural extension of the
class of cotriangular graphs.

Locally  cotriangular   graphs   have  been  classified  by  Hall  and
Shult~\cite{HaSh}, where  several  particular  cases were  settled  by
Hall~\cite{Ha:locpet}        and         by        Blokhuis        and
Brouwer~\cite{BlBr:locK33pet}.   The  hardest  part  was  to  classify
graphs which are locally cotriangular and reduced.

With finitely many exceptions, Pasechnik~\cite{Pa:3tr,Pa:fi,Pa:epolsp}
(see Chapters~\ref{chap:o3n}, \ref{chap:fi2n},  \ref{chap:epolsp}) has
classified the graphs which  are locally reduced  3-transposition, but
not cotriangular, graphs.  Pasechnik [in preparation] also settled one
exceptional  case not covered  in  \cite{Pa:3tr}. Only  two exceptions
remain,  namely  the 45-vertex  graph  $\Xi_1$ for  $U_4(2)$  and  the
165-vertex graph $\Xi_2$ for $U_5(2)$. There are 3 examples of locally
$\Xi_1$-graphs known, see  e.g. Cameron,  Hughes and Pasini~\cite{CHP}
for   two   of   them,    the    third   has   been    discovered   by
Pasechnik~\cite{Pa:epolsp} (see  Chapter~\ref{chap:epolsp}).   Nothing
is known about locally $\Xi_2$-graphs.

Thus, if $\cG$ denotes the class of all reduced 3-transposition
graphs with more than 165 vertices then all locally $\cG$-graphs
are known.
\bigskip

Given  a  3-transposition  group $G$ and  the  related  class  $C$  of
3-transposit\-ions, the {\em Fischer space}  $\cL=\cL(G,C)$ is defined
as a partial linear space whose points are  elements of $C$  and whose
lines  are $S_3$-subgroups such that the involutions in them belong to
$C$.  It is well-known that  the {\em planes} of $\cL$ (planes are the
minimal subspaces containing a  pair of intersecting lines) are either
affine of order 3 or dual affine of order 2.

Cuypers~\cite{Cuy:genfisp} has introduced {\em generalized Fischer
spaces}, that is, partial linear spaces whose planes are affine or
dual affine. Under certain nondegeneracy conditions, they were
classified in \cite{Cuy:genfisp} and in Cuypers and Shult~\cite{CuySh}.  The
natural analogs of 3-transposition graphs are the complements of the
point graphs of generalized Fischer spaces, and this class includes
some nontrivial examples, like graphs $\cU_n$ of perpendicular
nonisotropic points of the GF(4)-vector space of dimension $n$
equipped with a nondegenerate Hermitian form. Cuypers (personal
communication) has classified locally $\cU_n$-graphs for
$n>6$.  The case $n=4$ has been settled by Pasechnik~\cite{Pa:3tr}
(see Chapter~\ref{chap:o3n}). Pasechnik (in preparation) has also
settled the case $n=5$.

\medskip 

We conclude  this chapter with the remark that there is  evidence that
the  approach discussed here seems  to  be applicable to investigating
geometries  for  $(3,4)$-transposition groups, including the  geometry
for the sporadic simple group $\BM$.

\nocite{MPS:aq} \nocite{MPS:aut} \nocite{Pa:aff}

\bibliography{geom}
\bibliographystyle{abbrv}
\addcontentsline{toc}{chapter}{References}

\end{document}